\documentclass[10pt]{article}
\usepackage{amsmath}
\usepackage{amssymb,amsmath,mathrsfs}
\usepackage{color,cases}
\usepackage{graphics}
\let\pa=\partial
\let\al=\alpha

\let\d=\delta

\let\s=\sigma
\let\f=\frac
\let\p=\psi
\let\om=\omega
\let\G= \Gamma
\let\D=\Delta

\let\Om=\Omega

\let\eps= \epsilon

\let\si=\sigma

\def\cA{{\cal A}}

\def\cC{{\cal C}}
\def\cD{{\cal D}}

\def\cH{{\cal H}}

\def\cN{{\cal N}}

\def\cP{{\cal P}}

\def\cR{{\cal R}}
\def\cS{{\cal S}}
\def\cT{{\cal T}}

\def\mK{\mathfrak K}
\def\mJ{\mathfrak J}

\def\no{\noindent}
\def\na{\nabla}

\def\p{\partial}

\def\dive{\mathop{\rm div}\nolimits}
\def\curl{\mathop{\rm curl}\nolimits}

\def\C{\mathop{\bf C\kern 0pt}\nolimits}
\def\DD{\mathop{\bf D\kern 0pt}\nolimits}
\def\K{\mathop{\bf K\kern 0pt}\nolimits}
\def\N{\mathop{\bf N\kern 0pt}\nolimits}
\def\Q{\mathop{\bf Q\kern 0pt}\nolimits}
\def\R{\mathop{\bf R\kern 0pt}\nolimits}
\def\T{\mathop{\bf T\kern 0pt}\nolimits}

\newcommand{\ef}{ \hfill $ \blacksquare $ \vskip 3mm}

\newcommand{\beq}{\begin{equation}}
\newcommand{\eeq}{\end{equation}}
\newcommand{\ben}{\begin{eqnarray}}
\newcommand{\een}{\end{eqnarray}}
\newcommand{\beno}{\begin{eqnarray*}}
\newcommand{\eeno}{\end{eqnarray*}}

\newtheorem{thm}{Theorem}[section]
\newtheorem{lem}{Lemma}[section]
\newtheorem{rmk}{Remark}[section]

\newtheorem{prop}{Proposition}[section]
\renewcommand{\theequation}{\thesection.\arabic{equation}}

\newtheorem{theorem}{Theorem}[section]

\newtheorem{lemma}[theorem]{Lemma}
\newtheorem{proposition}[theorem]{Proposition}

\begin{document}
\title{Local well-posedness of the capillary-gravity water waves with acute contact angles}

\author{Mei Ming$^1$ and Chao Wang$^2$ \\[2mm]
{\small $ ^1$ School of Mathematics and Statistics,Yunnan University, Kunming 615000, P.R.China
  }\\[2mm]
{\small E-mail: mingmei@ynu.edu.cn}\\[2mm]
{\small $ ^2$ School of  Mathematical Science, Peking University, Beijing 100871, P.R.China}\\[2mm]
{\small E-mail: wangchao@math.pku.edu.cn}
}
\date{}

\maketitle
\begin{abstract}
We consider the  two-dimensional capillary-gravity water waves problem where the free surface $\G_t$ intersects the bottom $\G_b$ at two contact points. In our previous works \cite{MW2, MW3},  the local well-posedness for this problem has been proved with the contact angles less than $\pi/16$. In this paper, we study the case where the contact angles belong to $(0, \pi/2)$. It involves much worse singularities generated from corresponding elliptic systems, which have a strong influence on the regularities for the free surface and the velocity field.
Combining the theory of  singularity decompositions for elliptic problems with the structure of the water waves system, we  obtain a priori energy estimates. Based on these estimates, we also prove the local well-posedness  of the solutions in a geometric formulation.
\end{abstract}

\renewcommand{\theequation}{\thesection.\arabic{equation}}
\setcounter{equation}{0}

\tableofcontents

 \section{Introduction}
 \setcounter{equation}{0}

We consider the irrotational incompressible capillary-gravity water waves problem in a two-dimensional  domain $\Om_t$, where $\Om_t$ is a bounded domain  with an  upper free surface $\Gamma_t$ and a fixed  bottom $\Gamma_b$.  This moving domain contains two moving contact points $p_{l}, p_{r}$ (left and right) with the contact angles $\om_l, \om_r\in(0, \pi/2)$, which are the intersection points of $\G_t, \G_b$:
\[
\G_t\cap \G_b=\{p_l, p_r\}.
\]
 Moreover, the fixed bottom $\G_b$ is assumed to be smooth enough, and it becomes straight near the contact points $p_i$ ($i=l, r$) for the sake of simplicity.

\begin{center}
\resizebox{8cm}{!}{\includegraphics{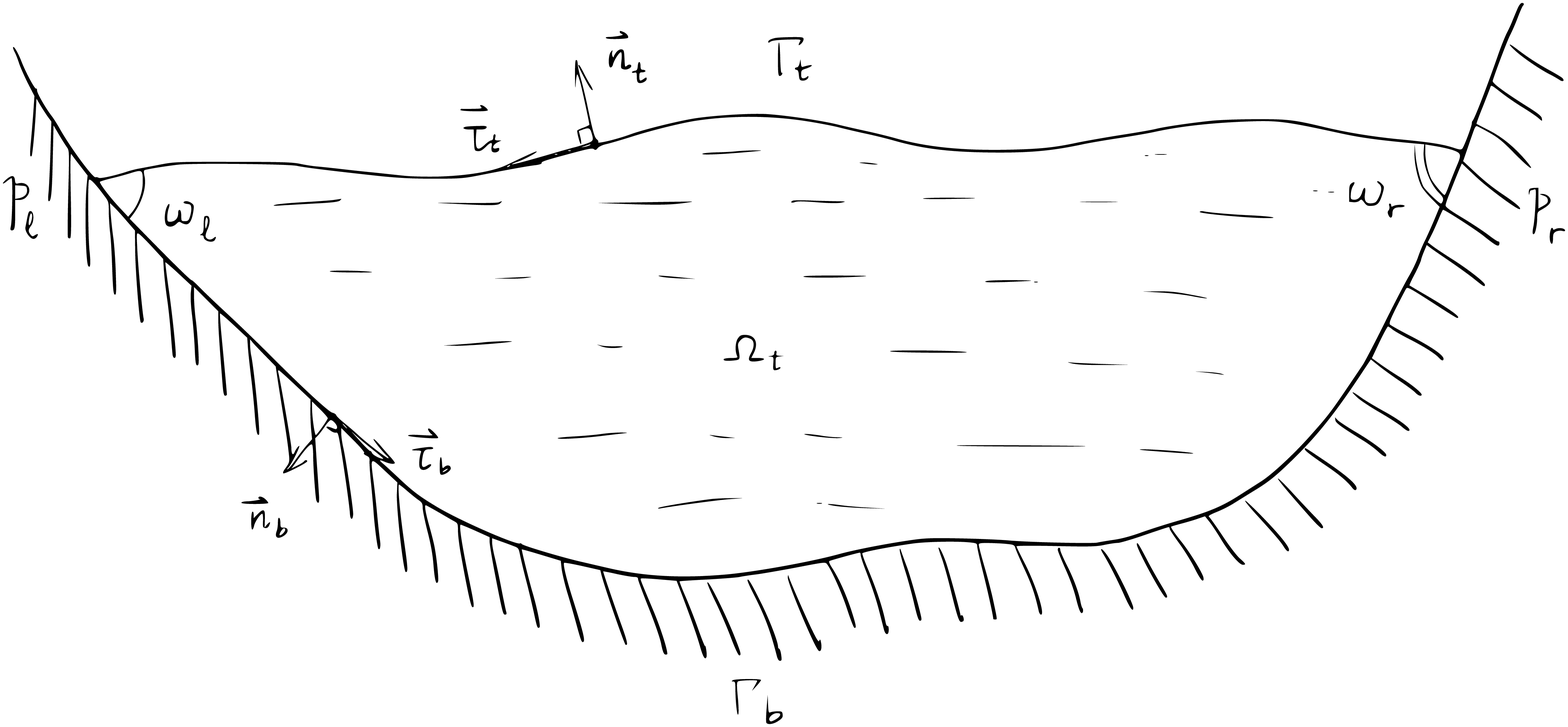}}
\end{center}

The  water waves problem has been widely studied in centuries, see for example \cite{Stokes, LC}. This problem focuses on the motion of an ideal fluid and describes the evolution of the free surface $\G_t$ as well as the velocity field $v$. Mathematically, it is described by Euler's equation with boundary conditions and initial conditions, and in our case, we also need some boundary conditions at contact points.

We express the water waves problem on the corner domain $\Om_t$ as the  following  system $\mbox{(WW)}$ of  velocity $v$ and pressure $P$:
\begin{numcases}{}
\pa_t v+v\cdot \na v=-\na P+{\bf g} \qquad \hbox{in}\quad \Om_t, \label{Euler eqn}\\
\dive v=0,\quad \curl v=0\qquad \hbox{in}\quad \Om_t, \label{div curl}\\
P|_{\Gamma_t}=\sigma \kappa, \label{surface tension}\\
D_t=\pa_t+v\cdot \na \  \textrm{is tangent to}\ \{(t, X) | X \in \pa \Omega_t\}, \label{kinematic}\\
v\cdot n_b |_{\Gamma_b}=0, \label{bottom v}\\
\beta_c v_i=\sigma(\cos{\om_s}-\cos{\om_i})\qquad\hbox{at}\quad p_i\ (i=l,r). \label{contact points}
\end{numcases}

Here \eqref{Euler eqn} from $\mbox{(WW)}$ is the Euler's equation where ${\bf g}=-g {\bf e_z}$ is the vertical gravity vector;  \eqref{div curl} describes the incompressibility and irrotationality;
\eqref{surface tension} is the condition of the pressure on the free surface in the case with surface tension, where $\si$ is the coefficient of surface tension and $\kappa$ is the mean curvature of $\G_t$ (see Section 1.2); \eqref{kinematic} is the classical kinematic condition on the free surface $\G_t$ with $D_t$ the material derivative; Meanwhile, \eqref{bottom v} describes that the velocity along the fixed bottom $\G_b$ is always tangential, where $n_b$ is the unit outward normal vector of $\G_b$. These equations and conditions are  standard in water waves, see \cite{Lannes, LannesBook}.

In particular, \eqref{contact points} gives the conditions at contact points, which come from  \cite{RE}. We denote by $v_i$ the upward tangential components of the velocity at the corner points along $\G_b$:
\[
v_l=-v\cdot\tau_b\quad\hbox{at}\  p_{l},\quad\hbox{and}\quad v_r=v\cdot\tau_b\quad\hbox{at}\  p_{r}.
\]
Here $\om_s$ is the stationary contact angle decided by the materials of the bottom and the fluid (see \cite{Young}), and $\beta_c$ denotes the  effective friction coefficient.  This condition shows that the slip velocity is dominated by the unbalanced Young stress,  and it is indeed  an effective variation of Young's law (1805) for  stationary contact angles \cite{Young}. In fact, this kind of conditions are commonly seen, see  \cite{BEIMR, CDA, SA, GL}.

\medskip
Before presenting our results, we  recall briefly earlier works on the well-posedness of classical water waves problem, where one has smooth surfaces $\G_t$ satisfying $\G_t\cap\G_b=\emptyset$. 

We recall results on the local well-posedness. When the fluid is irrotational, some early works like Nalimov \cite{Na}, Yosihara \cite{Yo1, Yo2} and Craig \cite{Craig} established  the two-dimensional local well-posedness with small initial data in Sobolev spaces. In late 1990s,  Wu \cite{Wu1, Wu2}  proved for the first time the local well-posedness with general initial data in Sobolev spaces and showed  that the Taylor sign condition
\beq\label{TC}
-\na_{n_t} P\geq c_0>0
\eeq
 held on $\G_t$ as long as $\Gamma_t$ was not self-intersecting.   Iguchi, Tanaka and Tani \cite{Ig-Ta} and  Iguchi \cite{Iguchi}   proved the local well-posedness in two-dimensional case respectively. Later on, Lannes \cite{Lannes} derived the local well-posedness of the gravity water waves under Zakharov formulation, which is convenient to link with approximate models.  Alazard, Burq and Zuily in \cite{ABZ1, ABZ4, ABZ2, ABZ3} used  paradifferential operators and Strichartz estimates to study the problem in a low-regularity space. On the other hand, when the fluid is rotational,  Christodoulou and Lindblad \cite{CL}  proved a priori estimates based on the geometry of the moving domain;  Lindblad \cite{Lin} obtained the existence of solutions using Nash-Moser iteration. In 2007, Coutand and Shkoller \cite{CS} used Lagrangian coordinates to show the local well-posedness. Shatah and Zeng \cite{SZ, SZ2} adopted a geometric point of view to reformulate the problem and prove the local well-posedness, while Beyer and G\"unther \cite{BG1, BG2} used a similar geometric approach to study the irrotational flow. Zhang and Zhang \cite{ZZ} proved the local well-posedness for rotational flow using a framework of Clifford analysis  introduced by Wu \cite{Wu2}. For more references, see  \cite{OT1, OT2, AM1, AM2, Sch, ZZ, WZZZ, AL, LannesBook, MZ} e.t.c..

For the global well-posedness of small data,   Wu \cite{Wu5} and  Germain, Masmoudi and Shatah \cite{GMS} proved the global three-dimensional existence of gravity water waves respectively using different approaches.   One can check \cite{Wu4, AD, IP, DIP, HIT1,HIT2, HIT3, Wang1, Wu6} e.t.c. and their references for more  results on   gravity or capillary-gravity water waves,.  Meanwhile, there are also some works concerning  geometric singularities on the free surfaces.  The authors in \cite{CCFG} proved the existence of a wave which is  given initially as the graph of a function and then can overturn at a later time. Later on, the authors in \cite{CCFGG}  showed the existence for some  ``splash" or ``splat'' singularities. This result was extended  to three-dimensional case and some other models in  \cite{CS2}.

\medskip

Compared to the rich literature on the well-posedness of classical water waves, the research on the well-posedness of  water waves problem with non-smooth boundaries (we call it ``non-smooth water waves'') just started several years ago and there are a lot of open questions . In general, there are two kinds of non-smooth water waves problems: The first kind of problems has contact points (or contact lines) between the free surface and the bottom, i.e. $\G_t\cap\G_b\neq\emptyset$; The other kind contains crests or cusps on the free surface, i.e. the surface is Lipschitz. We would like to mention that in the case with large crest angles,  the famous Stokes waves can  be dated back  to papers by G. Stokes \cite{Stokes, Stokes2} which obtained  traveling-wave solutions with limit crest angle $2\pi/3$. Obviously, the main difference here compared to classical water waves lies in the corners on boundaries. As a result, the analysis involving the corners (i.e. domain singularities)  becomes  the key point in the non-smooth water waves.

\bigskip
Now, we are in a position to state an informal version of the main result.
To begin with,  we introduce the following compatibility conditions at $t=0$
\beq\label{eq:com cond}
\beta_c\pa_t^k v_i(0)=\sigma\pa_t^k\big(\cos{\om_s}-\cos{\om_{i0}}\big)\quad \hbox{at}\quad p_i \,(i=l, r), \quad k=0,1,2,3.
\eeq
These conditions are needed since  we will reduce system $\mbox(WW)$ into initial-boundary value problem (see \eqref{eqn cDt2 J}, \eqref{eq: D_t^2 J-1} and \eqref{eq: IS-C}).
 
\begin{thm}\label{thm:main}
Let the initial data belong to a suitable space and the initial contact angles $\om_{i0}\in (0, \pi/2)$ for $i=l, r$. If the compatibility conditions \eqref{eq:com cond} are satisfied,  there exists a time interval depending on the initial data such that system $\mbox{(WW)}$ is locally well-posed in a suitable space. 
\end{thm}
\begin{rmk}
The suitable space is defined by $\Sigma_h$ in Section 6.3 for some good unknowns. In fact, the proof of Theorem \ref{thm:main} is divided into three parts. In Section 4 and Section 5, we obtain a priori estimates (see Theorem \ref{lower-order energy estimates} and Theorem \ref{1-order energy estimates}). In the last section, we construct approximate solutions to prove the well-posedness based on the a priori estimates. The precise statement of local well-posedness for a geometric form of $\mbox(WW)$ is given in Theorem \ref{main thm} (Section 6.3). 
\end{rmk}


\begin{rmk}
 When there is no surface tension, the authors of \cite{WuK, CEG} studied the case with contact angles. They  assume that the wall  $\G_b$ is vertical and the contact angle $\om\in(0, \pi/4)$. Then, by a symmetric extension,  they reduced the problem to the crest case with the crest angle less than $\pi/2$. As a result, \cite{WuK, CEG} do not contain the case of limit Stokes waves with the `` $2\pi/3$'' crest. In that case,  one needs a contact angle of $\pi /3$ even with the help of a symmetric extension.  In this paper, we only require that  the contact angles  $\om_i\in(0, \pi /2)$, which brings us some useful experiences to deal with  domain singularities in a more general case for  water waves. 
\end{rmk}
\begin{rmk} By continuity  argument of the contact angle in time similarly as in the last part of \cite{MW2}, one can see that as long as initial contact angles $\om_{i0}\in (0, \pi/2)$, we will have $\om_i(t)\in (0, \pi/2)$ in a short time interval.
\end{rmk}

\medskip
We recall some works concerning the local well-posedness of non-smooth water waves.
In the case where there are crests or cusps on the free surface (the ``crest'' case),  Alazard, Burq and Zuily \cite{ABZ3} study a special case (without surface tension) when the contact angle is equal to $\pi/2$ (the right angle), where they used symmetric and periodic extension  to turn this problem into a classical smooth periodic case.  A  breakthrough in this subject is made by Kinsey and Wu, see \cite{WuK, Wu3}.
They focus on  gravity water waves where the crest angle  is less than $\pi/2$.   The main difficulty is that the free surface is a non-$C^1$ interface with angled crest  and    the Taylor sign $-\na_{n_t}P$ degenerates at the crest point.  To be more precise,  they start  with reducing the water waves problem into the following equation in \cite{WuK, Wu3}:
\[
\pa_t^2 u+\mathfrak a\na_{n} u= f(u, \pa_t u)
\]
where $\mathfrak a=-\na_{n_t}P$ is the Taylor sign. When $\mathfrak a$ degenerates at the crest point, the above system will loss its hyperbolicity, and classical analysis  does not apply any more. To solve the problem, Kinsey and Wu flatten the domain with a Riemann mapping, and some singular weights appear naturally in their equation. As a result, they introduce some weighted Sobolev spaces accordingly for  energy estimates to deal with these singular weights. Based on these works, Agrawal \cite{Agra} show that these singularities are ``rigid", which means that the angle of these crests can not change in time. Very recently, 
C\'ordoba, Enciso and Grubic \cite{CEG}  study a  similar case with cusps and crests without gravity, where the angles of these crests are less than $\pi/2$ and change in time.

\medskip

For the other case where  there are  contact angles (that is $\G_t\cap\G_b\neq\emptyset$),  things become different from the ``crest'' case. First, the corners appear due to intersections of the free surface and the bottom (or wall); Second, there are different boundary conditions in corresponding elliptic systems compared to the crest case. In fact, various boundary conditions may have no big difference if we only focus on the elliptic theory, but there will be a series of consequences in water waves when boundary conditions change. For example, the evolution of the free surface is different from the crest case. Moreover, when there is surface tension, boundary conditions at the contact points as \eqref{contact points} are needed  in order to close the system, and dissipations appear at the contact points too (See Section 4 and \cite{MW2}).

In the case with contact angles,  de Poyferr\'e \cite{Poyferre} prove a priori estimates in bounded n-dimensional corner domains without surface tension. The contact angle is assumed to be small to ensure sufficient Sobolev regularity near the corner. Meanwhile, under a similar assumption of small contact angles, we  obtain the local well-posedness  in a two-dimensional corner domain (beach type) with surface tension, see \cite{MW2, MW3}.  Meanwhile, we notice that \cite{Poyferre} and \cite{MW2, MW3} use similar geometry formulations  introduced in \cite{SZ}. 

To explain why  small contact angles are needed in \cite{MW2, MW3, Poyferre},   we look at a typical mixed-boundary system in water waves:
\[
\begin{cases}{}
\D u=0\qquad in \quad \Om_t,\\
u|_{\G_t}=f,\quad \na_{n_b}u|_{\G_b}=0.
\end{cases}
\]
In fact, the elliptic theory on corner domains is well known already, see for example \cite{Kon, PG-2, KMR}. Generally, one still has variational solution $u\in H^1(\Om_t)$ if the right-side functions lie in proper Sobolev spaces. If one seeks for $H^2$ and above regularities,  singularity decompositions are  needed naturally for the solution, which decompose the solution into a singular  part $u_s$ (i.e. not good enough) near the corners and a regular part $u_r$:
\[
u=u_s+u_r\qquad\hbox{according to the required regularity.}
\]
For the mixed-boundary system above, the most singular part in $u_s$  is like $r^{\pi/2\om}$, where $r$ is the radius with respect to the corner point and $\om$ is the contact angle.   Consequently,  when  $\om$ is small enough, the singular part $u_s$ will be good enough so that we will have enough regularities from elliptic systems as in classical water waves to close the energy estimates. One can also find a singularity decomposition for $v$ in Proposition \ref{est:elliptic-N} (Section 2).

\medskip

In contrast, when the contact angle is larger or more general, the idea of taking small angles to improve  regularities in \cite{Poyferre, MW2, MW3} does not work any more. Meanwhile,  there is no obvious weighted space to use  due to the structure of the water waves problem. 

\medskip  
We want to show here the main  ingredients of this paper. 
Similarly as in our previous works \cite{MW2, MW3}, we still adopt the geometric formulation from Shatah and Zeng \cite{SZ, SZ2}. In fact, we rewrite $(\hbox{WW})$ into an equation for $\mJ=\na \mK_\cH$ with boundary conditions, where $\mK$ is the modified mean curvature on $\G_t$ and $\cH$ means the harmonic extension in $\Om_t$  (see Section 3, and we only present a simpler form for this equation here):
\beq\label{eqn simple for mJ}
D^2_t\mJ+\si\cA\mJ=R\quad\hbox{with}\quad \cA\mJ=\na \cH(-\D_{\G_t}\mJ^\perp),
\eeq
where $R$ is some remainder part.

The  trouble here  is that  singularities from the domain $\Om_t$ affect directly the regularity of the solution to $\mbox{(WW)}$, which means  singularities always appear in  related elliptic systems   even if the boundary conditions are good enough in Sobolev spaces (see for example \cite{PG-2}). More precisely,  the natural  norm $\|\cA^{k}f \|_{L^2(\Om_t)}+\|f\|_{L^2(\Om_t)}$ with $k\geq 1$ arising from this equation above is not equivalent with  $\|f \|_{H^{3k}(\Om_t)}$ due to  larger contact angles, and apparently ``singular parts" are contained in this norm. 

Due to this kind of singularities from elliptic systems, we only have limited regularities for some quantities (for example, the velocity $v$) in Sobolev spaces. 
Compared to \cite{MW2, MW3}, we must make full use of the maximal regularity for each quantity (especially for $v$) while the contact angles $\om_i\in (0, \pi/2)$.  To do this, some  delicate estimates together with singular parts from singularity decompositions (see \cite{MW1, PG-2}) are carefully used.    Meanwhile, it is also very important to gain more information from the  structure of $\mbox{(WW)}$.

The main part of the (lower-order) energy functional is  defined as
\[
E(t)=\|D_t \mJ\|^2_{L^2(\Om_t)}+\|  \mJ\cdot n_t\|^2_{H^1(\G_t)},
\]
which gives us the following estimate (see Section 4.1)
\[
\|v\|_{H^2(\Om_t)}+\|D_t v\|_{H^{3/2}(\Om_t)}+\|D_t^2 v\|_{L^2(\Om_t)}+\|\na v\|_{L^\infty(\Om_t)}+\|\kappa\|_{H^2(\G_t)} \leq P(E(t))
\]
with $P(E(t))$ the positive-coefficient polynomial of $E(t)$.
This means that the free surface $\G_t$ still has enough regularity and $v$ is Lipschitz.

In our previous work \cite{MW2, MW3}, one main part from the remainder term $R$ in \eqref{eqn simple for mJ} is  about the higher-order derivative terms of $P_{v,v}$, where $P_{v, v}$ is defined by an elliptic system with mixed-boundary conditions. According to the elliptic theory on corner domains (see for example \cite{PG-2}), the  system of $P_{v,v}$ only gives $P_{v,v}$ limited regularity around $H^2(\Om_t)$, when   the contact angles are less than $\pi/2$. To improve its regularity, we modify the definition of $P_{v, v}$ in this paper to have a Neumann-boundary system, see 
\eqref{eq: P_{vv}}. In fact, thanks to the elliptic theory,  solutions of Neumann-boundary system (or Dirichlet system) may become more regular than solutions of mixed-boundary problem in Sobolev spaces, while the right-side data have the same regularity.  As a result, we obtain a bit more regularity from the elliptic system of $P_{v,v}$ (we have $P_{v,v}\in H^3(\Om_t)$ indeed), which is important in the energy estimates.  
 
Moreover, another main part of $R$ lies in the higher-order derivative terms of the velocity  $v$. Here, we point out that compared to \cite{MW2, MW3}, $v$ loses some regularity due to the existence of corners. Fortunately,  $D_t v$ and $D_t^2 v$ have the same regularities as before. As a result, when we deal with $v$, sometimes we need to  use material derivatives $D_t$ instead of spatial derivatives.  Meanwhile, we also need to apply singularity decompositions to $v$ and its potential $\phi$ in the estimates, see for example Lemma \ref{na v estimate}.  

For the higher-order energy estimates, we use the material derivative $D_t$ instead of $\cA^{1/2}$ from \cite{SZ}. The higher-order energy is defined as 
\[
E_1(t)=\|\na_{\tau_t} D_t  (\mathfrak{J}\cdot n_t) \|^2_{L^2(\Gamma_t)} +\|D^2_t \mathfrak{J}\|^2_{L^2(\Om_t)}.
\]
In fact, one will see that $D_t$ is convenient to use when there are contact points. For example, taking $D_t$ on elliptic systems does not violate  boundary conditions, while it will change boundary conditions if one takes spatial derivatives. Moreover, one will find in the higher-order energy estimate that taking $D_t$ leads to better regularities for some quantities  than their own regularities (like $D_t v$ for the velocity $v$). However, due to the singularities of the boundary (contact points), it is not as convenient as before to turn these regularities into spatial regularities, which explains the reason why we  need to choose very carefully for ``a suitable space'' in Theorem \ref{thm:main}.  As long as the energy estimates for $E(t), E_1(t)$ are finished, we  use the equation of $\mJ$ to gain more spatial regularities.

\bigskip

In the end,  we mention some other related works. Lannes and M\'etivier \cite{LM} studied the Green-Naghdi system in a beach-type domain, which is a shallow-water model of the water waves problem. Lannes \cite{Lannes1} studied  the floating-body problem and proposed a new formulation  that can be easily  generalized in order to take into account the presence of a floating body.  Lannes and Iguchi \cite{LI} proved some sharp results for initial boundary value problem with a free boundary arising in wave-structure interaction, and it contains the floating problem in the shallow water regime. Besides, Guo and Tice \cite{GT} showed  a priori estimates for the contact line problem in the case of the stokes equations. Later on, Tice and Zheng proved the local well-posedness of the contact line problem in two-dimensional Stokes flow, see \cite{TZ}. In 2020, Guo and Tice \cite{GT2} proved a priori estimates for the contact line problem for two-dimensional Navire-Stokes flow. For  Darcy's flow, one can see \cite{KM1, KM2}.


\subsection{Organization of the paper} 
In Section 2, we present various useful lemmas including  singularity decompositions and estimates for elliptic systems. In Section 3, we derive the equation for the good unknown $\mJ$ from $(\hbox{WW})$ with modified curvature $\mK$ and modified pressure $P_{v,v}$. The lower-order energy is constructed and the energy estimate is proved in Section 4, where estimates for various quantities like $\G_t$, $v$, $P_{v,v}$ are proved. Moreover, we consider the higher-order energy estimate using $D_t$  in Section 5. In Section 6, we present the precise main theorem in our paper and  show the local well-posedness.

\subsection{Notations}

\no - $X$ stands for a point in $\Om_t\subset \R^2$. $p_l, p_r$ are the left and right contact points. $n_j (j=t, b)$ are the unit outward normal vectors on $\G_j$, and $\tau_j$ are the corresponding unit tangential vectors obeying the right-hand rule with $n_j$. \\
\no - $\si$ is the surface tension coefficient. $\beta_c$ is the effective friction coefficient determined by interfacial widths, interactions between the fluid and the bottom, and the normal stress contributions.\\
\no - $\chi_{\om}$ is a characteristic function of  contact angles: 
\[
\chi_{\om}(\theta)=
\begin{cases}
1,\quad \theta\in (\pi/3, \pi/2),\\
0,\quad \theta\in (0,\pi/3].
\end{cases}
\]
\no - $\chi_i$ ($i=l, r$) are cut-off functions near the corner points $p_i$:
\[
\chi_{i}(X)=
\begin{cases}
1,\quad \hbox{when}\ |X-p_i|\le r_0,\ X\in \Om_t;\\
0,\quad\ \hbox{otherwise}
\end{cases}
\]
with some small $r_0>0$.\\
- $f|_c=\chi_l (f|_{p_l})+\chi_r (f|_{p_r})$ stands for taking values of $f$ at the corner points.\\
- $\cS_{t,i}$ are straightened sector of $\Om_t$ with radius $r_0>0$ near the corner points $p_i$.\\
- $D_t=\p_t+\na_v$ is the material derivative.\\
- $M^*$ denotes the transport of a matrix $M$.\\
- $w^\perp$ on $\G_t$:   $w\cdot n_t$ for a vector $w\in T_X\G_t$.  \\
- $w^\top$ on $\G_t$: $(w\cdot \tau_t)\,\tau_t$. Sometimes we also use $w^\top$ on $\G_b$ with a similar definition.\\
- $\Pi$ is the second fundamental form on $\G_t$, where $\Pi(w)=\na_w n_t\in T_X\G_t$ for $w\in T_X\G_t$.\\
- $\kappa=tr \Pi=\na_{\tau_t}n_t\cdot \tau_t$ is the mean curvature of the surface $\G_t$.\\
- We define on $\G_t$ that $\cD w=\cD_{\tau_t}w=(\na_{\tau_t}w)^\top=(\na_{\tau_t}w\cdot \tau_t\big)\tau_t$ for a vector $w\in T_X\G_t$.\\
- $\cD^2f(\tau_1, \tau_2)=D^2f(\tau_1,\tau_2)-\big(\Pi(\tau_1)\cdot \tau_2\big)\na_{n_t}f$ for any two vector $\tau_1, \tau_2\in T_X\G_t$.\\
- $\D_{\G_t}$ is the Beltrami-Laplace operator on $\G_t$:
\[
\D_{\G_t}f=\cD^2f(\tau_t, \tau_t)=\na_{\tau_t}\na_{\tau_t}f-\na_{(\na_{\tau_t}\tau_t)^\top}f.
\]
\no - $\cH(f)$ or $f_\cH$ is the harmonic extension for some function $f$ on $\G_t$, which is defined by the elliptic system
\[
\left\{\begin{array}{ll}
\D \cH(f)=0\qquad\hbox{in}\quad \Om_t,\\
\cH(f)|_{\G_t}=f,\quad \na_{n_b}\cH(f)|_{\G_b}=0.
\end{array}\right.
\]
\no - $\cN=\na_{n_t}\cH$ is the Dirichlet-Neumann operator on $\G_t$.\\
\no - $\D^{-1}(h,g)$ stands for  the solution $u$ to the system
\[
\left\{\begin{array}{ll}
\D u=h\qquad \hbox{in}\quad \Om_t\\
u|_{\G_t}=0,\qquad \na_{n_b}u|_{\G_b}=g.
\end{array}
\right.
\]
\no - The Sobolev norm $H^s$ for the boundary $\G_t$ or $\G_b$ is defined by local coordinates and local graphs.\\
\no - $\tilde H^{1/2}(\G_j)$ ($j=t,b$) (see \cite{PG-2}) is a subspace of $H^{1/2}(\G_j)$ related to corner domains
\[
\tilde H^{1/2}(\G_j)=\Big\{u\in \dot{H}^{1/2}(\G_j)\Big| \,\rho_i^{-1/2}u\in L^2(\G_j),\ i=l, r\Big\}\] where $\dot H^{1/2}(\G_j)$ is the closure of $\mathscr D(\G_j)$ in $H^{1/2}(\G_j)$,  and $\rho_i=\rho_i(X)$ ($i=l, r$) is the distance (arc length) between the point $X\in \G_j$ and the end $p_i$.  We define the norm 
\[
\|u\|^2_{\tilde H^{1/2}(\G_j)}=\|u\|^2_{H^{1/2}(\G_j)}+\int_{\G_j} \rho_l^{-1}|u|^2dX+\int_{\G_j} \rho_r^{-1}|u|^2dX.
\] 
\no - $\tilde H^{-1/2}(\G_j)$ stands for the dual space of $\tilde H^{1/2}(\G_j)$. For more details, see \cite{PG-2}. \\
\no - We define $P_{w,v}$ (with $w\neq v$) by the following Neumann system:
\beq\label{eq: cP_{w,v}}
\left\{\begin{array}{ll}
\D  P_{w,v}=-tr (\na w \na v),\qquad\hbox{in}\quad\Om_t\\
\na_{n_t} P_{w,v}|_{\G_t}=C_{w,v}(t)-(w\cdot \tau_t) \na_{\tau_t} v\cdot n_t,
\qquad \na_{n_b}  P_{w,v}|_{\G_b}=w\cdot\na_v n_b.
\end{array}\right.
\eeq
Here $C_{w,v}$ is a function of $t$ satisfying the compatibility condition
\[
|\G_t|C_{w, v}(t)= -\int_{\Om_t}tr (\na w \na v)dX+\int_{\G_t}(w\cdot\tau_t) \na_{\tau_t} v\cdot n_t  ds-\int_{\G_b}w\cdot \na_v n_b ds.
\]
\no - $P_{v,v}$ is defined in \eqref{eq: P_{vv}} from Section 3.\\
\no - $C(a)$ stands for a positive constant $C$ depending on a quantity $a$. $P(E(t))$ stands for a polynomial of the energy $E(t)$ with positive constant coefficients.

 \section{Preliminaries}

To get started, we recall some useful estimates for Sobolev space.
\begin{lem}\label{embedding} (Sobolev Embeddings) One has the following inclusion:
\[
W^{s, p}(\Om)\subseteq L^q( \Om)
\] for $1/q=1/p-s/n$ and $\Om$ any bounded open subset of $\R^n$ with a Lipschitz boundary. Moreover, one also has
\[
W^{s, p}(\Om)\subseteq W^{t,q}(\bar\Om)
\]
for $t\le s, q\ge p$ such that  $s-n/p=t-n/q$.
\end{lem}

\begin{lem}\label{product} (Product estimates)
(1) For functions $f\in H^{1/2}(\Om_t)$ and $g\in H^1(\Om_t)\cup L^\infty(\Om_t)$, one has the following product estimate:
\[
\|f\,g\|_{H^{1/2}(\Om_t)}\le C\|f\|_{H^{1/2}(\Om_t)}\big(\|g\|_{H^1(\Om_t)}+\|g\|_{L^\infty(\Om_t)}\big)
\]
with a constant $C$ independent of $f, g$;

(2) For functions $f, g\in L^\infty(\G_t)\cup H^{1/2}(\G_t)$, one has
\[
\|f\,g\|_{H^{1/2}(\G_t)}\le C \big(\|f\|_{L^\infty(\G_t)}\|g\|_{H^{1/2}(\G_t)}+\|f\|_{H^{1/2}(\G_t)}\|g\|_{L^\infty(\G_t)}\big)
\]
with a constant $C$ independent of $f, g$.
\end{lem}
\no{\bf Proof.}
(1) In fact, one firstly extends $f, g$  to be defined on the full plane $\R^2$ with a control of their corresponding norms. Secondly, one can apply standard para-product analysis to prove the estimate on $\R^2$. The details are omitted here. (2) The proof is similar to the proof of (1).

\ef
We quote some Hardy inequalities here.
\begin{lem}\label{hardy} (Hardy inequalities) (1)(Corollary 2.3\cite{BK}) Let $f\in H^1(0, d)\cup C^0(0, d)$ with $d>0$ and  $f(0)=0$. Then there exists  a positive  constant $C=C(\eps, d)$ such that
\[
\int^d_0 r^{-2\epsilon}|f(r)|^2dr\le C\int^d_0r^{-2\epsilon+2}|f'(r)|^2dr
\]
for $\eps\in (1/2, 1)$;

(2) (Fractional-order version, see \cite{D, PG-1}) For a number $\epsilon\in (0, 1)$ and any function $f\in H^\epsilon(\R^+)$ with $f(0)=0$, there exists a positive constant
$C=C(\epsilon)$ such that
\[
\int_{\R^+}r^{-2\epsilon}|f(r)|^2dr\le C\int_{\R^+}\int_{\R^+}\frac{|f(r)-f(\rho)|^2}{|r-\rho|^{1+2\epsilon}}drd\rho.
\]
\end{lem}

\begin{lem}\label{trace}
(Traces on $\G_t$ or $\G_b$, Theorem 5.3 \cite{MW2})
The maps 
 \[
 u\mapsto \{u, \na_{n_j} u\}|_{\G_j},\quad \hbox{for}\quad j=t, b
 \]
 have  unique continuous extensions as  operators from
$H^{s}(\Om_t)$ onto $\Pi^1_{i=1} H^{s-i-1/2}(\Gamma_j)$ for $s>3/2$.

Moreover, one has the estimate :
\[
\|u\|_{H^{s-1/2}(\G_j)}+\| \na_{n_j} u\|_{H^{s-3/2 }(\Gamma_j )}\le    C\big(\|\Gamma_t\|_{H^{s-1/2} }\big)\|u\|_{H^{s}(\Om_t)}.
\]

\end{lem}

Next,  we present some some special trace theorems on corner domains involving $\tilde H^{1/2}(\G_j)$ and $\tilde H^{-1/2}(\G_j)$ ($j=t, b$). 
\begin{lem}\label{H1 zero trace}
\it  Assume that $u|_{\G_t}=f$ with $f|_{p_i}=0$ ($i=l,r$) for a function $u\in H^1(\Om_t)$. Then one has $f\in \tilde H^{1/2}(\G_t)$ and 
\[
\|f\|_{\tilde H^{1/2}(\G_t)}\le C(\|\Gamma_t\|_{H^{5/2} }) \|u\|_{H^1(\Om_t)}.
\]
The case on $\G_b$ holds similarly.
\end{lem}
\no{\bf Proof.} First, one has $f\in H^{1/2}(\G_t)$ immediately by Lemma \ref{trace}. Second, noticing $f|_{p_i}=0$ and applying Lemma \ref{hardy}(2) with some straightening localizations near $p_i$, one can see  that $f\in \tilde H^{1/2}(\G_t)$ with the desired estimate. Moreover, a similar lemma can be found as Lemma 5.5 \cite{MW2}.

\ef

\begin{lem}\label{G_b H-1/2} (Lemma 5.6\cite{MW2})
Let $u\in H^{1/2}(\G_j)$ ($j=t, b$), then $\na_{\tau_j} u$ belongs to  $\tilde H^{-1/2}(\G_j)$ and satisfies the estimate
\[
\|\na_{\tau_b}u\|_{\tilde H^{-1/2}(\G_j)}\le C(\|\Gamma_j\|_{H^{5/2} } ) \|u\|_{H^{1/2}(\G_j)}.
\] 
\end{lem}

\medskip

We now recall  some elliptic systems and estimates including singular decompositions in corner domains.
First, for  the  mixed-boundary system
 \[
\mbox{(MBVP)}\quad\left\{\begin{array}{ll}
\Delta u=h,\qquad \hbox{in}\quad \Omega_t\\
u\,|_{\G_t}=f,\qquad \na_{n_b}u \,|_{\G_b}=g,
\end{array}\right.
\]
we quote directly the following variational result.
\begin{lem}\label{Harmonic extension H1 estimate}(Lemma 5.9 \cite{MW2})
For a given function $f\in H^{1/2}(\G_t)$, the system 
\beq\label{H f system}
\left\{\begin{array}{ll}
\D \cH(f)=0\qquad\hbox{in}\quad \Om_t,\\
\cH(f)|_{\G_t}=f,\quad \na_{n_b}\cH(f)|_{\G_b}=0
\end{array}\right.
\eeq
admits a unique  solution  $f_\cH=\cH(f)\in H^1(\Om_t)$,  and  there holds
\[
\|\cH(f)\|_{H^1(\Om_t)}\le C(\|\Gamma_t\|_{H^{5/2} }) \|f\|_{H^{1/2}(\G_t)}.
\]
\end{lem}

The following proposition shows the existence and estimate for solutions in $H^2(\Om_t)$. Notice that when the contact angles are  below $\pi/2$, no singularity appears.
\begin{prop}\label{H2 MBVP}
  Let $h\in L^2(\Om_t)$, $f\in H^{3/2}(\G_t)$, $g\in H^{1/2}(\Gamma_b)$ and $\G_t \in H^{5/2}$ be given in $\mbox{(MBVP)}$. The contact angles $\om_i\in(0, \pi/2)$.
Then there exists a unique solution $u\in H^2(\Om_t)$ to $\mbox{(MBVP)}$.
Moreover, one  has
\[
\|u\|_{H^2(\Om_t)}\le C\big(\|h\|_{L^2(\Om)}+\|f\|_{H^{3/2}(\G_t)}+\|g\|_{H^{1/2}(\G_b)}\big)
\]
with the constant $C=C(\|\G_t\|_{H^{5/2}})$.
\end{prop}
\no{\bf Proof.}
This proposition is a direct conclusion from Proposition 5.1, Lemma 5.2 and Theorem 5.3 in \cite{MW1}.
\ef

\medskip

Second, we consider  the Neumann-boundary system:
 \ben\label{eq:elliptic-N}
\mbox{(NBVP)}\quad\left\{
\begin{array}{l}
\Delta u=h\quad \textrm{in}\quad \Om_t, \\
\na_{n_t}u|_{\Gamma_t}=f,\quad\quad \na_{n_b} u |_{\Gamma_b}=g
\end{array}
\right.
\een
satisfying the compatibility condition
\beno
\int_{\Om_t} hdX=\int_{\Gamma_t} fds+\int_{\Gamma_b}gds.
\eeno

The following results  for $H^3$-and-below case are needed in our paper.
\begin{prop}\label{elliptic estimate Neumann H3}
(1)  Let $h\in L^2(\Om_t)$, $f\in L^2(\G_t)$, $g\in L^2(\Gamma_b)$  in $\mbox{(NBVP)}$. 
Then there exists a unique  (up to an additive constant) variational solution $u\in H^1(\Om_t)$ to  $\mbox{(NBVP)}$ satisfying 
\[
\|u\|_{H^1(\Om_t)}\le C(\|\G_t\|_{H^{5/2}})\big(\|h\|_{L^2(\Om_t)}+\|f\|_{L^{2}(\G_t)}+\|g\|_{L^{2}(\G_b)}+\|u\|_{L^2(\Om_t)}\big);
\]
(2) Let $h\in H^s(\Om_t)$, $f\in H^{1/2+s}(\G_t)$, $g\in H^{1/2+s}(\Gamma_b)$  in $\mbox{(NBVP)}$. The contact angles $\om_i\in(0, \pi/2)$.
Then there exists a unique  (up to an additive constant) solution $u\in H^3(\Om_t)$  to  $\mbox{(NBVP)}$ satisfying 
\[
\|u\|_{H^{2+s}(\Om_t)}\le C(\|\G_t\|_{H^{3}})\big(\|h\|_{H^s(\Om_t)}+\|f\|_{H^{1/2+s}(\G_t)}+\|g\|_{H^{1/2+s}(\G_b)}+\|u\|_{L^2(\Om_t)}\big)
\]
for a constant $s\in [0,1]$.
\end{prop}
\no{\bf Proof.} We only write a sketch for the proof here, since one can find similar details in our previous papers \cite{MW1, MW2, MW3}.  First, the existence of the variational solution $u\in H^1(\Om_t)$ can be proved directly by a standard variation procedure based on  Lemma 4.4.3.1\cite{PG-2}. Note that here one doesn't require that the contact angles are below $\pi/2$.
Second, one needs to notice that, there is no singular part  since the most singular part for $\mbox{(NBVP)}$ behaves like $r^{\pi/\om}$ near the corners (see for example \cite{PG-2}). Therefore, one can have  directly  $H^2$ estimate from Theorem 3.2.5\cite{BR} and Proposition 5.1\cite{MW1}, Theorem 4.3.1.4\cite{PG-2}. Similarly to the proof for Proposition 5.13\cite{MW1}, one has $H^3$ estimate without singular part. As a result, the desired estimate in (2) can be proved by an interpolation.

\ef
\begin{rmk} \label{H1 neumann estimate remark}
For the first estimate in Proposition \ref{elliptic estimate Neumann H3}, when one has additionally
\[
\int_{\Om_t}udX=0,
\]
the norm $\|u\|_{L^2(\Om_t)}$ can be deleted from the right side during the proof.

\end{rmk}

Moreover, we present the $H^4$ singularity decomposition and estimate. 
\begin{prop}\label{est:elliptic-N}
  Let $h\in H^2(\Om_t)$, $f\in H^{5/2}(\G_t)$, $g\in H^{5/2}(\Gamma_b)$ and $\G_t \in H^{4}$ in $\mbox{(NBVP)}$. The contact angles $\om_i\in(0, \pi/2)$.
Then there exists a unique  (up to an additive constant) solution $u\in H^3(\Om_t)$ to  $\mbox{(NBVP)}$ such that
\[
u=u_r+ u_s\quad\hbox{with the singular part}\ u_s=\chi_{\om}(\om_l)\chi_{l}\,c_l  S_l\circ T_l+\chi_{\om}(\om_r)\chi_{r} \, c_r  S_r\circ T_r,
\]
and the regular part $u_r\in H^4(\Om_t)$. Here the cut-off functions $\chi_\om, \chi_i$ are defined in the notation part.
$T_i\in H^{4}(\Om_t)$ are  boundary-straightening diffeomorphisms from $\Om_t$ onto the sectors $\cS_{t,i}$ near the corner points $p_i$ (from \cite{MW1}), and $S_i=  r^{\pi/\om_i}s(\theta)$ with $r$ the radius with respect to $p_i$ in $\cS_{t, i}$ and $s_i(\theta)$ some fixed sine or cosine functions.

Moreover, one  has  estimates for the singular coefficients $c_i$ ($i=l,r$) and the regular part $u_r$:
\[
 |c_l|+\|u_r\|_{H^4(\Om_t)}\le C\big(\|h\|_{H^2(\Om_t)}+\|f\|_{H^{5/2}(\G_t)}+\|g\|_{H^{5/2}(\G_b)}+\|u\|_{L^2(\Om_t)}\big)
\]
with the constant $C=C(\|\G_t\|_{H^{4}})$.
\end{prop}
\no{\bf Proof.}  Similarly to the proofs of Proposition 5.17, Proposition 5.18\cite{MW1} for the mixed-boundary problem, the proof can be checked and follows line by line thanks to \cite{PG-2}, so we omit the details here.

\ef
\begin{rmk} A direct conclusion from this proposition is that when the contact angles $\om_i\in (0,\pi/3]$, the solution $u$ to $\mbox{(NBVP)}$ with the same right side will be in $H^4(\Om_t)$ with the corresponding estimate.
\end{rmk}

\begin{rmk}\label{interp Neumann} Based on the  proposition above and estimate (9.15) in \cite{DNBL}, when the contact angles $\om_i\in (\pi/3, \pi/2)$,  one 
can have a more delicate and also natural estimate for $\mbox{(NBVP)}$  with $h\in H^{1+\eps}(\Om_t)$, $f\in H^{3/2+\eps}(\G_t)$, $g\in H^{3/2+\eps}(\Gamma_b)$ and $\G_t \in H^{4}$:
\[
\|u\|_{H^{3+\eps}(\Om_t)}\le C(\|\G_t\|_{H^{4}})\big(\|h\|_{H^{1+\eps}(\Om_t)}+\|f\|_{H^{3/2+\eps}(\G_t)}+\|g\|_{H^{3/2+\eps}(\G_b)}+\|u\|_{L^2(\Om_t)}\big),
\]
where $\eps\in (0, \pi/\om-2)\subset (0,1)$.
\end{rmk}

\medskip
In the end, we recall  some useful expressions and commutators from  \cite{SZ, MW2, MW3}.
 
\noindent - $D_tn_t$ and $D_t\tau_t$. One has 
\beq\label{Dtnt Dt taut}
D_t n_t=-\big((\na v)^*n_t\big)^\top,\quad D_t\tau_t=(\na_{\tau_t}v\cdot n_t\big)n_t.
\eeq

\noindent - $[D_t,\cH]$. One has for a smooth function $f$ on $\G_t$ that
\beq\label{commutator Dt H}
[D_t,\cH]f=\D^{-1}\big(2\na v\cdot\na^2f_{\cH}+\D v\cdot\na f_{\cH},\, (\na_{N_b}v-\na_v N_b)\cdot \na f_{\cH}\big) \quad\hbox{in}\quad \Om_t,
\eeq
where $\D^{-1}(h,g)$ and $\cH$ are defined in the notation part.

\noindent - $[D_t, \cN]$.
One has
\beq\label{commutator DN}
\begin{split}
[D_t,\,\cN]f=&\na_{n_t}\D^{-1}\Big(2\na v\cdot \na^2f_{\cH}+\D v\cdot\na f_{\cH},\,(\na_{n_b}v-\na_{v}n_b)\cdot\na f_{\cH}\Big)\\
&\,-\na_{n_t}v\cdot \na f_{\cH}-\na_{(\na f_{\cH})^\top} v\cdot n_t\qquad\qquad\hbox{on}\quad \G_t.
\end{split}
\eeq

\noindent  - $[D_t,\D_{\G_t}]$.  For a smooth function $f$ on $\G_t$, there holds
\beq\label{commutator surface delta}
[D_t,\,\D_{\G_t}]f=2\cD^2f\big(\tau_t,(\na_{\tau_t}v)^\top\big)-(\na f)^\top\cdot \D_{\G_t}v+\kappa \na_{(\na f)^\top} v\cdot n_t  \qquad\hbox{on}\quad \G_t.
\eeq

\noindent - $[D_t,\,\D^{-1}]$. We have
\beq\label{commutator Dt D-1}
D_t\D^{-1}(h,g)=\D^{-1}(D_th,\,D_tg)+\D^{-1}(h_1, g_1)
\eeq with
\[
h_1=2\na v\cdot\na^2\D^{-1}(h,g)+\D v\cdot \na \D^{-1}(h,g),\quad
g_1=(\na_{N_b}v-\na_{v}N_b)\cdot \na\D^{-1}(h,g).
\]

\noindent -  $[D_t, \na_{\tau_t}]$. Direct computations lead to
\beq\label{commutator Dt na taut}
[D_t, \na_{\tau_t}]=(D_t\tau_t-\na_{\tau_t}v)\cdot\na=(n_t \na_{\tau_t}v\cdot n_t-\na_{\tau_t}v)\cdot\na
=-(\na_{\tau_t}v\cdot \tau_t)\na_{\tau_t}.
\eeq

%
%
%

\section{Reformulation of the problem}

In this section, we derive the equation for a good unknown $\mathfrak J$, which is slightly different from the quantity $J=\na\kappa_{\cH}$ introduced in \cite{SZ} and used in our previous papers \cite{MW2, MW3}.

\medskip
To begin with, we define $P_{v, v}$ by the following Neumann-boundary system:
\beq\label{eq: P_{vv}}
\left\{\begin{array}{ll}
\D P_{v, v}=-tr (\na v \na v)\qquad \hbox{in}\quad \Om_t\\
\na_{n_t}P_{v, v}|_{\G_t}=C_{v,v}(t),\qquad \na_{n_b}P_{v, v}|_{\G_b}=v\cdot \na_v n_b.
\end{array}
\right.
\eeq
where $C_{v,v}(t)$ satisfies
\[
|\G_t| C_{v,v}(t) = -\int_{\Om_t}tr (\na v \na v)dX-\int_{\G_b}v\cdot \na_v n_b ds.
\]
Moreover, due to the non-uniqueness of the variational solution to this Neumann problem,  we assume that $P_{v, v}$ satisfies
\beq\label{assumption on Pvv}
\int_{\Om_t} P_{v, v}dX=0.
\eeq
In this way, we will have  a unique solution  to \eqref{eq: P_{vv}} in Sobolev space.
\medskip

Let the pressure $P$ be decomposed into
\beq\label{decomp P}
P=\mK_{\cH}+  P_{v,v},
\eeq
where  $\mK$ is the modified curvature on $\G_t$ defined by
\beq\label{def of mK}
\mK=\sigma\kappa-P_{v,v}\qquad\mbox{on}\ {\G_t}.
\eeq
\medskip

Now,  we  can define the new quantity
\[
\mathfrak J=\na \mathfrak{K}_{\cH}
\]
and we are ready to derive its equation from $\hbox{(WW)}$. In fact, since the derivation for  the equation of $\mJ$ is much similar to
the derivation for $J$ in  \cite{MW2},  many details are omitted here.

First, recall the equation for the curvature $\kappa$ from \cite{MW2}:
 \beq\label{Dt k 1}
D_t\kappa=-\D_{\G_t}v^\perp-v^\perp |\na_{\tau_t}n_t|^2+(\na_{\tau_t}\na_{v^\top}n_t)\cdot\tau_t-(\na_{\tau_t}v^\top)^\top\cdot\na n_t\cdot\tau_t
\eeq
or equivalently
\beq\label{Dt k 2}
D_t\kappa=-\D_{\G_t}v\cdot n_t-2\Pi(\tau_t)\cdot \na_{\tau_t}v \qquad\hbox{on}\quad \G_t.
\eeq

Applying $D_t$ to \eqref{Dt k 2} and using \eqref{def of mK}, we have
\beq\label{Dt2 kappa}
D^2_t\mK=-\si \D_{\G_t}D_tv \cdot n_t +2\s \Pi(\tau_t)\cdot \na_{\tau_t} J+\si R_1- D_t^2 P_{v,v}\qquad\hbox{on}\quad \G_t
\eeq
with 
\[
\begin{split}
R_1= - [D_t, \D_{\G_t}]v\cdot n_t-\D_{\G_t}v\cdot D_tn_t+2 \Pi(\tau_t)\cdot\na_{\tau_t}\na P_{v,v}
-2 \Pi(\tau_t)\cdot[D_t, \na_{\tau_t}]v-2  D_t\big(\Pi(\tau_t)\big)\cdot\na_{\tau_t}v.
\end{split}
\]
Moreover, direct computations almost the same as in \cite{MW2} lead to
\beq\label{Dt2 J eqn1}
D^2_t\mJ=\na\cH(D^2_t\mK)+A_1+A_2+A_3
\eeq
where the remainder terms
\beq\label{A_1}
\begin{split}
A_1=\na [D_t,\cH]D_t\mK=&\na\D^{-1}\big(2\na v\cdot\na+\D v\cdot ,(\na_{n_b}v-\na_vn_b)\cdot\big)\\
&\quad
\Big(D_t\mJ-\na\D^{-1}\big(2\na v\cdot \na \mJ
+\D v\cdot \mJ,\,(\na_{n_b}v-\na_v n_b)\cdot \mJ\big)+(\na v)^*\mJ\Big),
\end{split}
\eeq
\beq\label{A_2}
\begin{split}
A_2&=\na D_t\D^{-1}\big(2\na v\cdot\na \mJ+\D v\cdot \mJ,(\na_{n_b}v-\na_vn_b)\cdot \mJ\big)\\
&=\na\D^{-1}\big(2\na v\cdot\na^2w_{A2}+\D v\cdot \na w_{A2},(\na_{n_b}v-\na_vn_b)\cdot\na w_{A2}\big) +\na\D^{-1}(h_{A2},\, g_{A2})
\end{split}
\eeq
with
\[
\begin{split}
 w_{A2}=&\D^{-1}\big(2\na v\cdot \na \mJ+\D v\cdot \mJ,\,(\na_{n_b}v-\na_v n_b)\cdot \mJ\big),\\
 h_{A2}=&2\na v\cdot(\na D_t\mJ-(\na v)^*\mJ)+2(\na D_tv-(\na v)^*\na v)\cdot \na \mJ+D_t\mJ\cdot\D v\\
 &+\mJ\cdot (\D D_tv-\D v\cdot\na v-2\na v\cdot \na^2 v),\\
g_{A2}=&(\na_vn_b-\na_{n_b}v)\cdot \na v\cdot \mJ+\na_{n_b}D_t v\cdot \mJ-(D_t v-\na_v v)\cdot \na n_b\cdot \mJ\\
& +\na_v\big((\na v)^*n_b\big)^\top\cdot \mJ+(\na_{n_b}v-\na_vn_b)\cdot D_t\mJ,
\end{split}
\]
and
\beq\label{A_3}
A_3=-2(\na v)^*D_t\mJ-(\na D_t v)^* \mJ-\big((\na v)^2\big)^*\mJ+\big((\na v)^*\na v\big)\mJ.
\eeq
Here we note that the leading-order terms in $A_1, A_2$  are like $\mJ,D_t\mJ,\na v,\na D_tv$.

Substituting \eqref{Dt2 kappa} into \eqref{Dt2 J eqn1} and applying Euler's equation, we can have after some more computations that
\beq\label{eqn for Dt2 J}
D^2_t\mJ=\s\na \cH(\D_{\G_t}\mJ^\perp)+R_0+D_t\na \cH(  D_t P_{v,v} ) .
\eeq
with
\[
R_0=-\si\na \cH(\mJ\cdot \D_{\G_t}n_t)+\si \na \cH(n_t\cdot \D_{\G_t}\na P_{v,v})+\si \na\cH(R_1 ) +A_1+A_2+A_3+[D_t, \na \cH](  D_t P_{v,v}).
\]
 Here, one observation is that since $\na_{n_t} P_{v,v}\big|_{\G_t} =C_{v,v}(t)$ from \eqref{eq: P_{vv}}, we find immediately in $R_0$ that
\[
 n_t\cdot \D_{\G_t}\na P_{v,v} =[n_t,\D_{\G_t} ] \cdot\na P_{v,v} .
\]

Moreover, we use the following Hodge decomposition from \cite{SZ, MW2}:
\beq\label{Hodge decom}
\cD_t \mJ=D_t \mJ+\na \cP_{\mJ,v}
\eeq
such that $\cD_t \mJ$ satisfies
\[
\na\cdot \cD_t\mJ=0\quad\mbox{on} \ \Om_t,\quad\hbox{and}\quad  n_b\cdot\cD_t\mJ\big|_{\G_b}=0,
\]
and
$\cP_{\mJ,v}$ satisfies the Neumann-boundary system
\beq\label{PJv system}
\left\{\begin{array}{ll}
\D   \cP_{\mJ,v}=-tr (\na\mJ\na v)\qquad\hbox{in}\quad\Om_t\\
\na_{n_t} \cP_{\mJ,v}|_{\G_t}=C_{\mJ,v}(t)-\na_{\tau_t}  \mK \na_{\tau_t} v\cdot n_t,\quad \na_{n_b}  \cP_{\mJ,v}|_{\G_b}= \mJ\cdot\na_v n_b
\end{array}\right.
\eeq
with
\[
|\G_t| C_{\mJ,v}(t) = -\int_{\Om_t}tr (\na \mJ \na v)dX+\int_{\G_t}\na_{\tau_t}  \mK \na_{\tau_t} v\cdot n_t  ds-\int_{\G_b}\mJ\cdot \na_v n_b ds.
\]
Here the assumption
\[
\int_{\Om_t} \cP_{\mJ, v}dX=0
\]
holds similarly to guarantee the uniqueness as for $P_{v,v}$ system.

This implies immediately
\[
\begin{split}
D_t\cD_t\mJ= D_t(D_t\mJ+\na \cP_{\mJ,v})=D^2_t\mJ+D_t\na \cP_{\mJ,v}.
\end{split}
\]

As a result, we finally conclude the equation for $\mJ$ in the following form based on \eqref{eqn for Dt2 J}:
\beq\label{eqn cDt2 J}
D_t\big[\cD_t \mJ-\na \cH\big(  D_t P_{v,v} -v\cdot(\na P_{v,v}|_c)\big)\big]+\s \cA \mJ=\cR
\eeq where the third-order elliptic operator $\cA$ is defined in  the same way  as in \cite{SZ, MW2}:
\[
\cA(w)=\na\cH\big(-\D_{\G_t}(w|_{\G_t})^\perp\big)
\]
for any smooth-enough function $w$ defined on $\Om_t$,
and the remainder term $\cR$ is
\[
\cR=R_0  +D_t\na \cP_{\mJ,v}+D_t\na\cH\big( v\cdot(\na P_{v,v}|_c)\big).
\]
Moreover,  one can find the definition of $\na P_{v,v}|_c$ in the notation part.

\begin{rmk}\label{explain the eqn} Instead of \eqref{eqn for Dt2 J}, we write the equation of $\mJ$ in a more complicated form \eqref{eqn cDt2 J}. This happens because of technical reasons. In fact, $\cD_t \mJ$ is more handy to use in the energy estimates, while as a price we have the remainder $\cP_{\mJ,v}$ part to deal with. Meanwhile, the part $v\cdot(\na P_{v,v}|_c)$ is added to derive better estimate for $D_tP_{v,v}$, see Lemma \ref{prop: D_t P_{v, v}-H}.

\end{rmk}

\section{Lower-order energy estimates}

The lower-order energy $E(t)$ is defined as
\[
E(t)=\|\na_{\tau_t}\mJ^\bot\|^2_{L^2(\Gamma_t)} +\|\cD_t \mJ\|^2_{L^2(\Om_t)}+\|\Gamma_t\|^2_{H^{5/2}}+\|v\|^2_{H^{3/2}(\Om_t)},
\]
where  we write
\[
E_l(t)= \|\Gamma_t\|^2_{H^{5/2}}+\|v\|^2_{H^{3/2}(\Om_t)}.
\]
Moreover, recalling from \cite{MW2, MW3} that our $\mbox{(WW)}$ system has some dissipation at corner points,  we have the following dissipation term at corner points
\[
F(t)=\sum_{i=l,r}\big|(\sin \om_i)\na_{{\bf \tau}_t}\mJ^\perp |_{p_i}\big|^2.
\]


\begin{thm}\label{lower-order energy estimates}
Let the contact angles $\om_i\in (0,\pi/2)$. Assume that  $E(t), \int_0^T F(t)dt$ are both bounded above in $[0, T]$ for some  $T>0$. Then  the following a priori estimate holds for system $\mbox{(WW)}$:
\[
\sup_{0\leq t\leq T} E(t)+\int_0^T F(t)dt\le P(E(0))+\int_0^T P(E(t))dt,
\]
where $P(\cdot)$ is some polynomial with positive constant coefficients depending on $\si, \beta_c$.
\end{thm}

One can see immediately that, our energy is defined mainly in forms of $\mJ$. As a result, we need to related $(v, P)$, $\G_t$ and other quantities to $\mJ$ firstly before we prove this energy estimate.


\medskip
\subsection{Dependence on $E(t)$}\label{estimate section}
We show in this part that all the quantities related to our problem can be controlled by $E(t)$. 
To start with, the following proposition focuses on the estimates for $P_{v,v}, n_t, \mK_\cH, \mJ$.
\begin{lem}\label{prop: eta}
 Assuming that $E(t)\in L^\infty[0,T]$ for some $T>0$, one has the following estimates
\[
\|P_{v, v}\|_{H^{2}(\Om_t)} \leq P(E_l(t))
\]
and
\[
\|n_t\|_{H^2(\Gamma_t)}+\|\mJ\|_{H^{3/2}(\Om_t)}+\|\mK_\cH\|_{H^{5/2}(\Om_t)}\leq  P(E(t)).
\]
\end{lem}
\no{\bf Proof.} (1) The first inequality. Applying Proposition \ref{elliptic estimate Neumann H3} with an interpolation to $P_{v,v}$ system \eqref{eq: P_{vv}}, we get
\[
 \|P_{v, v}\|_{H^{2}(\Om_t)} \leq C(\|\Gamma_t\|_{H^{5/2}})\big(\|tr(\na v \na v)\|_{L^2(\Om_t)}+|\G_t|^{1/2}|C_{v,v}(t)|+\|v\cdot \na_v n_b\|_{H^{1/2}(\G_b)}+\|P_{v,v}\|_{L^2(\Om_t)}\big).
\]

For the right-hand side, we have firstly  by Lemma \ref{embedding} that
\[
\|tr(\na v \na v)\|_{L^2(\Om_t)}\le C(\|\Gamma_t\|_{H^{5/2}})\|v\|^2_{H^{3/2}(\Om_t)}.
\]
Second, it's straightforward to find
\[
|C_{v,v}(t)|\le C\|v\|^2_{H^{3/2}(\Om_t)}
\]
and
\[
\begin{split}
\|P_{v,v}\|_{L^2(\Om_t)}&=\big\|P_{v,v}-\int_{\Om_t}P_{v,v}dX\big\|_{L^2(\Om_t)}\le C\|\na P_{v,v}\|_{L^2(\Om_t)}\\
&\le C(\|\Gamma_t\|_{H^{5/2}})\big( \|tr(\na v \na v)\|_{L^2(\Om_t)}+|\G_t|^{1/2}|C_{v,v}(t)|+\|v\cdot \na_v n_b\|_{L^2(\G_b)}\big),
\end{split}
\]
where \eqref{assumption on Pvv} is used and the variational estimate above is proved similarly as Proposition \ref{elliptic estimate Neumann H3}(1).

Consequently, summing up all these inequalities above, we conclude that
\[
 \|P_{v, v}\|_{H^{2}(\Om_t)} \le C(\|\Gamma_t\|_{H^{5/2}})\|v\|^2_{H^{3/2}(\Om_t)}\leq P(E_l(t)).
\]

(2) The second inequality. In fact, the key point of the proof  lies in the estimates for mean curvature $\kappa$ and $\mK=\sigma\kappa-P_{v,v}$.

To begin with, one can see immediately that $\mK_\cH$ satisfies the system
\beq\label{kH system}
\left\{\begin{array}{ll}
\D \mK_\cH=0\qquad\hbox{in}\quad\Om_t\\
\na_{n_t}\mK_\cH|_{\G_t}=\mJ^\perp|_{\G_t}\in H^1(\Gamma_t),\qquad \na_{n_b} \mK_\cH|_{\G_b}=0.
\end{array}\right.
\eeq
Applying Proposition \ref{elliptic estimate Neumann H3}(2), we find
\beq\label{est:kappa1}
\begin{split}
\|\mK_{\cH}\|_{H^{5/2}(\Om_t)}&\leq C(\|\Gamma_t\|_{H^{5/2}})(\|\mJ^\perp\|_{H^1(\Gamma_t)}+\| \mK_\cH\|_{L^2(\Om_t)})\nonumber\\
&\leq C(\|\Gamma_t\|_{H^{5/2}})( \|\na_{\tau_t} \mJ^\bot\|_{L^2(\Gamma_t)} +\| \mK_\cH\|_{L^2(\Om_t)}),
\end{split}
\eeq
 where interpolation for $\mK_\cH$ is applied to $\|\mJ^\perp\|_{L^2}$.

Moreover, from the definition of $\mK$ and Lemma \ref{Harmonic extension H1 estimate} we have
\[
\| \mK_\cH\|_{L^2(\Om_t)}\le C(\|\Gamma_t\|_{H^{5/2}})\big(\|\kappa\|_{H^{1/2}(\G_t)}+\|P_{v,v}\|_{H^{1/2}(\G_t)}\big),
\]
which together with the previous inequality leads to
\beq\label{est:kappa-1}
\|\mK_{\cH}\|_{H^{5/2}(\Om_t)}\leq C(\|\Gamma_t\|_{H^{5/2}}) (1+ \|\na_{\tau_t} \mJ^\bot\|_{L^2(\Gamma_t)}+\|  P_{v, v}\|_{H^{1/2}(\Om_t)} ).
\eeq
This implies immediately the desired $H^{3/2}$ estimate for $\mJ$.

For the estimate for $n_t$, we find by Lemma \ref{trace}, \eqref{est:kappa-1} and part (1) that
\beq\label{estimate kappa}
\begin{split}
 \|\kappa\|_{H^{1}(\Gamma_t)} &\leq \si^{-1}(\| \mK\|_{H^{1}(\Gamma_t)}+\| P_{v,v}\|_{H^{1}(\Gamma_t)})\\
 &\leq  C(\|\Gamma_t\|_{H^{\f52}}) (1+ \|\na_{\tau_t} \mJ^\bot\|_{L^2(\Gamma_t)}+\|  P_{v, v}\|_{H^{3/2}(\Om_t)})
 \le P(E(t)).
\end{split}
\eeq
As a result, the proof is finished thanks to  the fact
\[
\kappa=\na_{\tau_t}n_t\cdot \tau_t\quad\hbox{where}\ \na_{\tau_t}n_t\| \tau_t\quad\hbox{on}\ \G_t.
\]

\ef

\begin{lem}\label{prop: DtPvv}
Assuming that $E(t)\in L^\infty[0,T]$ for some $T>0$, one has the following estimates for $v$:
\[
\|v\|_{H^2(\Om_t)}+\|v^\perp\|_{H^{5/2}(\G_t)}\leq P(E(t)).
\]
Meanwhile, one also has
\[
\|D_tP_{v, v}\|_{H^{1}(\Om_t)} +\|\cP_{\mJ, v}\|_{H^1(\Om_t)}\leq P(E(t))
\]
\end{lem}
\no{\bf Proof.} \noindent - $H^1$  estimate for $\cP_{\mJ, v}$.
 From  system \eqref{PJv system} of $\cP_{\mJ, v}$, we can have the following variational estimate by Proposition \ref{elliptic estimate Neumann H3}(1) and Remark \ref{H1 neumann estimate remark}:
\[
\begin{split}
 \|\cP_{\mJ, v}\|_{H^{1}(\Om_t)} \leq &C(\|\Gamma_t\|_{H^{5/2}})\big(\|tr(\na \mJ \na v) \|_{L^2(\Om_t)}+|\G_t|^{1/2}|C_{\mJ,v}(t)|+\|\na_{\tau_t}  \mK\cdot \na_{\tau_t} v\cdot n_t\|_{L^2(\Gamma_t)}\\
 &
 +\|\mJ\cdot \na_v n_b\|_{L^2(\G_b)}\big),
\end{split}
\]
which together with Lemma \ref{prop: eta} leads to the desired estimate for $\cP_{\mJ, v}$.

\noindent - $v^\perp$ estimate. In fact, it is straightforward  to get by Lemma \ref{trace} that
\[
\|D_t \kappa\|_{H^{1/2}(\G_t)}  \le C(\|\G_t\|_{H^{5/2}})\|D_t\kappa_\cH\|_{H^1(\Om_t)}
\leq C(\|\G_t\|_{H^{5/2}})(\|\na D_t \kappa_\cH\|_{L^2(\Om_t)} +\|D_t \kappa\|_{L^2(\G_t)} ).
\]
Rewriting $\na\kappa_\cH$ in terms of $\mJ, \na P_{v,v}$ leads to
\[
\begin{split}
\|D_t \kappa\|_{H^{1/2}(\G_t)}
&\leq C(\|\G_t\|_{H^{5/2}})\big(\|D_t \mJ\|_{L^2(\Om_t)}+\|D_t\na P_{v,v}\|_{L^2(\Om_t)}+\|\na v \cdot\mJ\|_{L^2(\Om_t)}+\|\na v\cdot \na P_{v,v}\|_{L^2(\Om_t)}\\
&\qquad  +\|D_t \kappa\|_{L^2(\G_t)} \big) \\
&\leq P(E(t)) \big( 1+\|\na \cP_{\mJ, v}\|_{L^2(\Om_t)}+\|D_tP_{v,v}\|_{H^1(\Om_t)}+\|P_{v,v}\|_{H^{3/2}(\Om_t)}+\|\mJ\|_{H^{3/2}(\Om_t)} \\
&\qquad  +\|D_t \kappa\|_{L^2(\G_t)}  \big).
\end{split}
\]
Applying Lemma \ref{prop: eta}, we have
\[
\|D_t \kappa\|_{H^{1/2}(\G_t)}\le P(E(t)) \big( 1+\|\na \cP_{\mJ, v}\|_{L^2(\Om_t)}+\|D_tP_{v,v}\|_{H^1(\Om_t)}
+\|D_t \kappa\|_{L^2(\G_t)}  \big).
\]
On the other hand, we rewrite \eqref{Dt k 1} into
\beno
\D_{\G_t}v^\perp=-D_t\kappa-v^\perp |\na_{\tau_t}n_t|^2+(\na_{\tau_t}\na_{v^\top}n_t)\cdot\tau_t-(\na_{\tau_t}v^\top)^\top\cdot\na n_t\cdot\tau_t,
\eeno
which together with the inequality above for $\|D_t \kappa\|_{H^{1/2}(\G_t)}$ and the estimate for $\cP_{\mJ, v}$  imply
\ben\label{est: v-H-1}
\begin{split}
\|\D_{\G_t}v^\perp\|_{H^{1/2}(\G_t)} &\leq  C(\|\G_t\|_{H^{5/2}},\|n_t\|_{H^{5/2}(\G_t)})\big(\|D_t\kappa\|_{H^{1/2}(\G_t)}  + \|v\|_{H^2(\Om_t)}\big)\\
&\leq  P(E(t))\big(1+\|\na \cP_{\mJ, v}\|_{L^2(\Om_t)}+\|D_tP_{v,v}\|_{H^1(\Om_t)}
+\|D_t \kappa\|_{L^2(\G_t)} + \|v\|_{H^2(\Om_t)}\big)\\
&\leq  P(E(t))\big(1 +\| D_tP_{v, v}\|_{H^1(\Om_t)}  + \|v\|_{H^2(\Om_t)}\big).
\end{split}
\een
Here we use \eqref{Dt k 2} to obtain
\[
\|D_t \kappa\|_{L^2(\G_t)}\le C(\|\G_t\|_{H^{5/2}}) \|v\|_{H^2(\Om_t)}.
\]
As a result, the estimate for $v^\perp$ depends on the estimates for $\| D_tP_{v, v}\|_{H^1(\Om_t)}, \|v\|_{H^2(\Om_t)}$.

\noindent - $H^2$ estimate for $v$. Recalling  system \eqref{eq:phi} for velocity potential $\phi$ with $v=\na \phi$, we have by  Proposition \ref{elliptic estimate Neumann H3}  that
\[
\begin{split}
\|v\|_{H^2(\Om_t)} &\le \|\phi\|_{H^3(\Om_t)}\leq  C(\|\Gamma_t\|_{H^{5/2}})\big(\|v^\perp\|_{H^{3/2}(\G_t)}+\|\phi\|_{L^2(\Om_t)}\big)\\
&\le C(\|\Gamma_t\|_{H^{5/2}})\big(\|v^\perp\|_{H^{3/2}(\G_t)}+\|v\|_{L^2(\Om_t)}\big)\\
&\leq P(E(t))\big(1 +\|v^\perp\|_{H^{3/2}(\G_t)}\big).
\end{split}
\]


Bringing the above estimate into \eqref{est: v-H-1} and applying interpolations to $\|v^\perp\|_{H^{3/2}(\G_t)}$,  we arrive at
\beq\label{v H2 estimate 1}
\|\D_{\G_t}v^\perp\|_{H^{1/2}(\G_t)}+\|v\|_{H^2(\Om_t)}\le  P(E(t))\big(1+\| D_tP_{v, v}\|_{H^1(\Om_t)} \big).
\eeq
Therefore, it remains to deal with $\| D_tP_{v, v}\|_{H^1(\Om_t)}$.

\noindent - $D_tP_{v,v}$ estimates.  From the definition of $P_{v, v}$, we derive the system for $D_tP_{v,v}$:
\beq\label{DtPvv system}
\left\{\begin{array}{ll}
\D D_t P_{v, v}=-D_t tr(\na v \na v)+ 2tr(\na  v  \na^2  P_{v, v})\qquad \hbox{in}\quad  \Om_t\\
\na_{n_t} D_t P_{v, v}|_{\G_t}=C'_{v,v}(t)+\na_{n_t} v\cdot \na P_{v, v}|_{\G_t},\quad \na_{n_b} D_tP_{v, v}|_{\G_b}=D_t(v\cdot\na_v n_b)+\na_{n_b} v\cdot \na P_{v, v}|_{\G_b}.
\end{array}
\right.
\eeq
with
\[
\int_{\Om_t}D_tP_{v,v}dX=0.
\]
It is easy to see from Euler's equation that
\beq\label{Dt tr na v2}
D_t tr(\na v \na v)=-2tr (\na v\cdot \na v\, \na v)-2tr \big(\na(\na P_{v,v}+\na\mK_\cH)\na v\big),
\eeq
and similar expression can be derived for $D_t(v\cdot\na_v n_b)$ on $\G_b$.

Moreover, we have
\[
 C'_{v,v} =|\Gamma_t|^{-1}\int_{\Om_t} 2tr(D_t \na v\na v)dX
-|\G_t|^{-2}\pa_t|\G_t|\int_{\Om_t}tr(\na v\na v)dX,
\]
where $\f d{dt}|\G_t|=\int_{\G_t}\big(v^\perp \kappa+\na_{\tau_t}(v\cdot\tau_t)\big)ds$. So we obtain by Lemma \ref{prop: eta}, \eqref{est:kappa-1} and \eqref{estimate kappa} that
\beq\label{C'vv estimat}
|C'_{v,v}|\le C(\|\Gamma_t\|_{H^{5/2}}, \| v\|_{H^{3/2}(\Om_t)}, \|P_{v,v}\|_{H^2(\Om_t)},
\|\mK_\cH\|_{H^2(\Om_t)}, \|\kappa\|_{L^2(\G_t)}\big)\leq P(E(t)).
\eeq

To prove $H^1$ estimate for $D_tP_{v,v}$, we use the following form of variation equation:
\[
\begin{split}
\int_{\Om_t}\big|\na D_tP_{v,v}\big|^2dX=&-\int_{\Om_t}(\Delta D_tP_{v,v})D_tP_{v,v}dX+\int_{\G_t}(\na_{n_t}D_tP_{v,v})D_tP_{v,v}ds\\
&+\int_{\G_b}(\na_{n_b}D_tP_{v,v})D_tP_{v,v}ds.
\end{split}
\]
Substituting the expressions in \eqref{DtPvv system} into the equality above, we can deal with the integrals one by one, where Lemma \ref{embedding} is applied. For example,  for the highest-order terms of $v$ on the boundary, we use Green's Formula and Lemma \ref{embedding} to find
\[
\begin{split}
&\int_{\G_t}\na_{n_t}v\cdot P_{v,v} \,D_tP_{v,v}ds+\int_{\G_b}\na_{n_b}v\cdot P_{v,v} \,D_tP_{v,v}ds\\
&=\int_{\Om_t}tr(\na v \na^2 P_{v,v}) D_tP_{v,v}dX+\int_{\Om_t}\na v\cdot\na P_{v,v}\cdot \na D_tP_{v,v}dX\\
&\le \|\na v\|_{L^4(\Om_t)}\|\na^2P_{v,v}\|_{L^2(\Om_t)}\|D_tP_{v,v}\|_{L^4(\Om_t)}+|\na v\|_{L^4(\Om_t)}\|\na P_{v,v}\|_{L^4(\Om_t)}\|\na D_tP_{v,v}\|_{L^2(\Om_t)}\\
&\le C(\|\Gamma_t\|_{H^{5/2}})\|v\|_{H^{3/2}(\Om_t)}\|P_{v,v}\|_{H^2(\Om_t)}\|D_tP_{v,v}\|_{H^1(\Om_t)}.
\end{split}
\]
The other integrals above can be handled similarly and the details are omitted.

As a result, we  obtain by Lemma \ref{prop: eta}  the following desired estimate:
\[
\|D_t P_{v, v}\|_{H^1(\Om_t)} \leq P(E(t))\big(\|P_{v,v}\|_{H^2(\Om_t)}+\|\mK_\cH\|_{H^{5/2}(\Om_t)}+|C'_{v,v} |\big)\le P(E(t)).
\]
Consequently,  applying  \eqref{v H2 estimate 1} leads to the estimates for $v^\perp$ and $v$.

\ef

\begin{rmk}\label{est:D_t k}
We know immediately from the proof of Lemma \ref{prop: DtPvv} the following estimate
\[
\|D_t\kappa\|_{H^{1/2}(\G_t)}\le P(E(t)).
\]
Moreover, we also have
\[
\|D_t\mK_\cH\|_{H^1(\Om_t)}\le \sigma\|D_t\kappa_\cH\|_{H^1(\Om_t)}+\|D_tP_{v,v}\|_{H^1(\Om_t)}\le P(E(t)).
\]
\end{rmk}
\medskip

Based on these lemmas above, we are ready to show some higher-order estimates.
\begin{lem}\label{nt H3 estimate}
Let $E(t)\in L^\infty[0,T]$ for some $T>0$, then one has
\[
\|P_{v, v}\|_{H^{5/2}(\Om_t)}+\|n_t\|_{H^3(\G_t)}\le P(E(t)).
\]
\end{lem}
\no{\bf Proof.} For $P_{v,v}$, similar arguments as in the proof of Lemma \ref{prop: eta} lead to
\[
\begin{split}
 \|P_{v, v}\|_{H^{5/2}(\Om_t)} &\leq C(\|\Gamma_t\|_{H^{5/2}})\big(\|tr(\na v \na v)\|_{H^{1/2}(\Om_t)}+|\G_t|^{1/2}|C_{v,v}(t)|+\|v\cdot \na_v n_b\|_{H^{1}(\G_b)}\big)\\
&\leq P(E(t))(1+\|   v\|_{H^2(\Om_t)} ).
\end{split}
\]
Applying Lemma \ref{prop: DtPvv}, we have the desired estimate.

Again, similar arguments as \eqref{estimate kappa} in the proof of Lemma \ref{prop: eta} lead directly to the estimate
\beno
\|\kappa \|_{H^{2}(\Gamma_t)} \leq \si^{-1}\big(\| \mK\|_{H^{2}(\Gamma_t)}+\| P_{v,v}\|_{H^{2}(\Gamma_t)}\big)\leq  C(\|\Gamma_t\|_{H^{5/2}}) (1+ \|\na_{\tau_t} \mJ^\bot\|_{L^2(\Gamma_t)}+\|  P_{v, v}\|_{H^{5/2}(\Om_t)}),
\eeno
which implies
\[
\|n_t\|_{H^3(\Gamma_t)}\leq  P(E(t))( 1+\|  P_{v, v}\|_{H^{5/2}(\Om_t)}).
\]
Therefore, combining this with the estimate for $P_{v,v}$, the proof is finished.

\ef

We also present more estimates for $v$.
\begin{lem}\label{prop: v}
 Assuming that $E(t)\in L^\infty[0,T]$ for some $T>0$,  one has
\[
\|D_tv\|_{H^{3/2}(\Om_t)} +\|D^2_tv\|_{L^2(\Om_t)}  \leq  P(E(t))
\]
\end{lem}
\no{\bf Proof.} First, one recalls Euler's equation to get that
\[
\|D_tv\|_{H^{3/2}(\Om_t)} \leq C\big(1+\|\mJ\|_{H^{3/2}(\Om_t)}+ \|\na  P_{v, v}\|_{H^{3/2}(\Om_t)}\big),
\]
and applying Lemma \ref{nt H3 estimate} leads to the estimate for $D_tv$.

Second, taking $D_t$ on Euler's equation leads to the following estimate:
\[
\|D^2_tv\|_{L^2(\Om_t)} \leq C(\|\cD_t\mJ\|_{L^2(\Om_t)} +\|\na  \cP_{\mJ, v}\|_{L^2(\Om_t)} + \|D_t\na  P_{v, v}\|_{L^2(\Om_t)} ).
\]
Using Lemma \ref{prop: DtPvv} again, we can finish the estimate for $D^2_t v$.

\ef

We now consider the following Neumann system for velocity potential $\phi$:
\beq\label{eq:phi}
\left\{\begin{array}{ll}
\D \phi=0,\qquad\hbox{on}\quad\Om_t\\
\na_{n_t}\phi|_{\G_t}=v^\perp ,\quad \na_{n_b}\phi|_{\G_b}=0,
\end{array}\right.
\eeq
where in order to obtain the uniqueness we chose $\phi$ to satisfy
\[
\int_{\Om_t}\phi dX=0
\]
without loss of generality.

We show by singularity decompositions that $\na v$ lies in $L^\infty(\Om_t)$, which is a key ingredient in the following estimates.

\begin{lem}\label{na v estimate}
Assume that $E(t)\in L^\infty[0,T]$ for some $T>0$, then there exists a unique $\phi\in H^3(\Om_t)$ to system \eqref{eq:phi} and one finds the following singular decomposition:
\[
\phi=\phi_r+\phi_s
\]
where the regular part $\phi_r\in H^4(\Om_t)$, and the singular part is expressed in the same way as $u_s$ in Proposition \ref{est:elliptic-N}.
Moreover, one has
\[
\|\na v\|_{L^{\infty}(\Om_t)}\le P(E(t)).
\]

\end{lem}
\no{\bf Proof.}
In fact, applying Proposition \ref{est:elliptic-N}, we know immediately about the existence of $\phi\in H^3(\Om_t)$ and the singular decomposition. Therefore, we have
\beq\label{na v decomp}
 v=\na\phi_r+\na\phi_s=v_r+v_s,\qquad\hbox{with}\quad v_r\in H^3(\Om_t).
\eeq

It only remains to show the estimate for $\na v$. Since we have
\[
\phi_s=\chi_{\om}(\om_l)\chi_{l}\,c_l  r^{\pi/\om_l}\circ T_l+\chi_{\om}(\om_r)\chi_{r} \, c_r r^{\pi/\om_r}\circ T_r
\]
where the singular part exists when $\om_i\in (\pi/3, \pi/2)$. In this case, we find by Proposition \ref{est:elliptic-N} that
\[
\begin{split}
\|\na^2\phi_s\|_{L^\infty(\Om_t)}&\le C(\|\G_t\|_{H^4})\big(\|v^\perp\|_{H^{5/2}(\G_t)}+\|\phi\|_{L^2(\Om_t)}\big) \big(\|\chi_{l}r^{\pi/\om_l-2}\circ T_l\|_{L^\infty(\Om_t)}\\
&\quad+\|\chi_{r}r^{\pi/\om_r-2}\circ T_r\|_{L^\infty(\Om_t)}\big)
\end{split}
\]
with $\pi/\om_i-2\in (0,1)$. Consequently, we obtain
\[
\|\na^2\phi_s\|_{L^\infty(\Om_t)}\le C(\|\G_t\|_{H^4})\big(\|v^\perp\|_{H^{5/2}(\G_t)}+\|\phi\|_{L^2(\Om_t)}\big)
\le C(\|\G_t\|_{H^4})\|v^\perp\|_{H^{5/2}(\G_t)}
\]
thanks to a direct variational estimate
\[
\|\phi\|_{L^2(\Om_t)}\le C(\|\G_t\|_{H^{5/2}}) \|\na\phi\|_{L^2(\Om_t)}\le C(\|\G_t\|_{H^{5/2}})\|v^\perp\|_{L^2(\G_t)}.
\]

Moreover, we also have by Proposition \ref{est:elliptic-N}  the following estimate:
\[
\|\phi_r\|_{H^4(\Om_t)}\le C(\|\G_t\|_{H^4})\big(\|v^\perp\|_{H^{5/2}(\G_t)}+\|\phi\|_{L^2(\Om_t)}\big)
\le C(\|\G_t\|_{H^4})\|v^\perp\|_{H^{5/2}(\G_t)}.
\]
As a result, we find that
\[
\|\na v\|_{L^\infty(\Om_t)}\le\|\na^2 \phi_r\|_{L^\infty(\Om_t)} +\|\na^2 \phi_s\|_{L^\infty(\Om_t)}
\le C(\|\G_t\|_{H^4})\|v^\perp\|_{H^{5/2}(\G_t)},
\]
which together with Lemma \ref{prop: DtPvv} lead to the desired estimate.

\ef

As long as we have   $\|\na v\|_{L^\infty(\Om_t)}$  estimate, we are able to deal with more higher-order estimates.
\begin{lem} \label{est: P_{v, v}-H}
Let $E(t)\in L^\infty[0,T]$ for some $T>0$, one has the following estimate:
\[
 \|P_{v, v}\|_{H^{3}(\Om_t)}+ \|P_{\mJ, v}\|_{H^{5/2}(\Om_t)}+\|D_tP_{v,v}\|_{H^{3/2}(\Om_t)}\le P(E(t))
\]
\end{lem}
\no{\bf Proof.} - $P_{v,v}$ estimate. In fact, similar arguments as in the proof of Lemma \ref{prop: eta} show that
\[
 \|P_{v, v}\|_{H^{3}(\Om_t)} \leq P(E(t))(1+\|  v\|_{H^2(\Om_t)})(1+\|\na v\|_{L^\infty(\Om_t)}).
\]
Consequently, applying Lemma \ref{prop: DtPvv} and Lemma \ref{na v estimate} leads to the desired estimate.

\noindent - Estimate for $P_{\mJ, v}$. Apply Proposition \ref{elliptic estimate Neumann H3}, we have
\[
\begin{split}
 \|P_{\mJ, v}\|_{H^{5/2}(\Om_t)} \leq &C(\|\Gamma_t\|_{H^{5/2}})\big(\|tr(\na \mJ \na v) \|_{H^{1/2}(\Om_t)}+|C_{\mJ,v}(t)|+\|\na_{\tau_t}  \mK  \na_{\tau_t} v\cdot n_t\|_{H^{1}(\Gamma_t)}\\
 &\quad
 +\|\mJ\cdot \na_v n_b\|_{H^{1}(\G_b)}+\|P_{\mJ, v}\|_{L^{2}(\Om_t)}\big).
\end{split}
\]
Using Lemma \ref{embedding}, Lemma \ref{product},  Lemma \ref{prop: eta}, Lemma \ref{prop: DtPvv} and Lemma \ref{na v estimate}, we can finish the esitmate.

\noindent - $H^{3/2}$ estimate for $D_tP_{v,v}$. Following the proof of Lemma \ref{prop: DtPvv} and using Lemma \ref{prop: DtPvv} again to improve  estimates for $v$, one can easily see that we have the desired estimate. 

\ef

We now give a high-order estimate for $D_t P_{v, v}$.
 \begin{lem}\label{prop: D_t P_{v, v}-H}
  Assuming $E(t)\in L^\infty[0,T]$ for some $T>0$,  one has
\beno
 \big\|D_t P_{v, v}-v\cdot (\na P_{v,v}|_c)\big\|_{H^{5/2}(\Om_t)} \leq P(E(t)).
\eeno
 \end{lem}
\no{\bf Proof.}
 We denote by
 \[
 w=D_t P_{v, v}-v\cdot (\na P_{v,v}|_c)\quad\hbox{with}\quad \na P_{v,v}|_c=\chi_l(\na P_{v,v}|_{p_l})+\chi_r(\na P_{v,v}|_{p_r}).
 \]
  A direct computation using \eqref{DtPvv system}  leads to
 \beq\label{eq: P_{vv}-v}
\left\{\begin{array}{ll}
\D w=-trD_t(\na v \na v) +[\D,v]\cdot (\na P_{v, v}- \na P_{v,v}|_c\big)-v\cdot \D(\na P_{v,v}|_c) \quad \hbox{in}\quad \Om_t\\
\na_{n_t}w|_{\G_t}=C_{v,v}'(t)+\na_{n_t} v\cdot (\na P_{v, v}-\na P_{v,v}|_{c})-v\cdot\na_{n_t}(\na P_{v,v}|_c)\big|_{\G_t},\\
 \na_{n_b}w|_{\G_b}=D_t(v\cdot\na_v n_b)+\na_{n_b} v\cdot (\na P_{v, v}-\na P_{v,v}|_{c})-v\cdot\na_{n_b}(\na P_{v,v}|_c)\big|_{\G_b}.
\end{array}
\right.
\eeq
Applying Proposition \ref{elliptic estimate Neumann H3} to system \eqref{eq: P_{vv}-v}, we have
\[
\|w\|_{H^{5/2}(\Om_t)}\le C(\|\G_t\|_{H^{5/2}})\big(\|\Delta w\|_{H^{1/2}(\Om_t)}+\|\na_{n_t}w\|_{H^1(\G_t)}
+\|\na_{n_b}w\|_{H^1(\G_b)}+\|w\|_{L^2(\Om_t)}\big).
\]

We need to deal with the terms on the right side above. First, thanks to \eqref{Dt tr na v2}, Lemma \ref{prop: DtPvv}, Lemma \ref{nt H3 estimate} and Lemma \ref{na v estimate},  it is straightforward to check that
\beno
\|D_t tr(\na v \na v)\|_{H^{1/2}(\Om_t)}+\big\|[\D,v]\cdot (\na P_{v, v}- \na P_{v,v}|_{c}) \big\|_{H^{1/2}(\Om_t)}
+\|v\cdot \D(\na P_{v,v}|_c)\|_{H^{1/2}(\Om_t)} \leq P(E(t)).
\eeno
Notice that here we need to use Lemma \ref{product} for the product estimate of $\|\na^2\mK_\cH\na v\|_{H^{1/2}(\Om_t)}$ from $\|D_t tr(\na v \na v)\|_{H^{1/2}(\Om_t)}$:
\[
\|\na^2\mK_\cH\na v\|_{H^{1/2}(\Om_t)}\le C\|\na^2\mK_\cH\|_{H^{1/2}(\Om_t)}\big(\|\na v\|_{H^1(\Om_t)}+
\|\na v\|_{L^\infty(\Om_t)}\big)\le P(E(t)).
\]
The estimate for  $\|\na^2 P_{v,v}\na v\|_{H^{1/2}(\Om_t)}$ also follows in a similar way.

Second, for the boundary terms, we only need to take care of the following estimate from $\|\na_{n_t} v\cdot (\na P_{v, v}-\na P_{v,v}|_{c})\|_{H^1(\G_t)}$:
\beq\label{boundary term esti}
\begin{split}
&\big\|\na_{\tau_t}\na_{n_t} v\cdot (\na P_{v, v}-\na P_{v,v}|_{c})\big\|_{L^2(\G_t)} \\
&\leq
\|\na_{\tau_t}\na_{n_t} v\cdot \chi_i(\na P_{v, v}-\na P_{v,v}|_{p_i})\|_{L^2(\G_t)}
+\|\na_{\tau_t}\na_{n_t} v\cdot (1-\chi_l-\chi_r)\na P_{v, v}\|_{L^2(\G_t)}\\
&\le \|\na_{\tau_t}\na_{n_t} v\cdot \chi_i(\na P_{v, v}-\na P_{v,v}|_{p_i})\|_{L^2(\G_t)} +P(E(t)).
\end{split}
\eeq
Here the estimate for the second term in the equality above holds thanks to \eqref{na v decomp} and Proposition \ref{est:elliptic-N}.

To handle the first term in \eqref{boundary term esti}, using a straightening differeomorphism as $T_i$ from Proposition \ref{est:elliptic-N} and \eqref{na v decomp}, we know immediately
\beq\label{apply hardy}
\begin{split}
&\big\|\na_{\tau_t}\na_{n_t} v\cdot \chi_i(\na P_{v, v}-\na P_{v,v}|_{p_i})\big\|_{L^2(\G_t)}\\
 &\le
C(\|\G_t\|_{H^{5/2}})\big\|(\na_{\tau_t}\na_{n_t} v)\circ T^{-1}_i\cdot (\na P_{v,v}-\na P_{v,v}|_{p_i})\circ T^{-1}_i\big\|_{L^2(0, r_0)}\\
&\le C(\|\G_t\|_{H^{3}})\big(\big\|c_ir^{\pi/\om-3} (\na P_{v,v}-\na P_{v,v}|_{p_i})\circ T^{-1}_i\big\|_{L^2(0, r_0)}\\
&\quad+\big\|(\pa^2\na \phi_r)\circ T^{-1}_i\cdot (\na P_{v,v}-\na P_{v,v}|_{p_i})\circ T^{-1}_i\big\|_{L^2(0, r_0)}\big)\\
&\le P(E(t))\big(\|r^{\pi/\om-3+3/4}\|_{L^2(0,r_0)}\big\|r^{-3/4}(\na P_{v,v}-\na P_{v,v}|_{p_i})\circ T^{-1}_i\big\|_{L^\infty(0, r_0)}+1\big).
\end{split}
\eeq
Moreover, applying Lemma \ref{hardy} leads to
\[
\begin{split}
&\|r^{-3/4}(\na P_{v,v}-\na P_{v,v}|_{p_i})\circ T_i\|_{L^2(0, r_0)}\le C\|(\na P_{v,v}-\na P_{v,v}|_{p_i})\circ T_i\|_{H^{3/4}(0, r_0)}
\quad\hbox{and}\\
&\|r^{-3/4}(\na P_{v,v}-\na P_{v,v}|_{p_i})\circ T_i\|_{H^1(0, r_0)}\le C\|(\na P_{v,v}-\na P_{v,v}|_{p_i})\circ T_i\|_{H^{1+3/4}(0, r_0)},
\end{split}
\]
so we obtain by an interpolation the following inequality
\[
\begin{split}
\big\|r^{-3/4}(\na P_{v,v}-\na P_{v,v}|_{p_i})\circ T_i\big\|_{L^\infty(0, r_0)}&\le C(\|\G_t\|_{H^{5/2}})\big\|r^{-3/4}(\na P_{v,v}-\na P_{v,v}|_{p_i})\circ T_i\big\|_{H^{3/4}(0, r_0)}\\
&\le C(\|\G_t\|_{H^{5/2}})\|(\na P_{v,v}-\na P_{v,v}|_{p_i})\circ T_i\|_{H^{3/2}(0, r_0)}.
\end{split}
\]
As a result, combining this with \eqref{apply hardy}, we conclude that
\[
\big\|\na_{\tau_t}\na_{n_t} v\cdot \chi_i(\na P_{v, v}-\na P_{v,v}|_{p_i})\big\|_{L^2(\G_t)}\le P(E(t))(\|P_{v,v}\|_{H^3(\Om_t)}+1)\le P(E(t)).
\]
Substituting this estimate into \eqref{boundary term esti}, we finally obtain
\[
\big\|\na_{n_t} v\cdot (\na P_{v, v}-\na P_{v,v}|_{c})\big\|_{H^1(\G_t)}\le P(E(t)).
\]
In the end, the  estimate for $\|\na_{n_b} v\cdot (\na P_{v, v}-\na P_{v,v}|_{c})\|_{H^1(\G_b)}$ follows in a similar way, and the proof can be finished.

\ef

In the end, we give the estimate of $\int_{\Om_t} \na D_t\cP_{\mathfrak{J}, v}\cdot D_t \mathfrak{J}dX$.
\begin{lem}\label{lem: D_t P_{J, v}}
Let $E(t)\in L^\infty[0,T]$ for some $T>0$. Then, we have the following estimate:
\[
\Big|\int_{\Om_t} \na D_t\cP_{\mathfrak{J}, v} \cdot D_t \mathfrak{J}dX-\f d{dt}\int_{\Om_t} \na  \cP_{\mathfrak{J}, v}\cdot \na v\cdot\na \mathfrak{K}_{\cH}dX \Big|\leq P(E(t)).
\]
\end{lem}
\no{\bf Proof.} Recalling system \eqref{PJv system}, we get
\beq\label{eq:D_t P}
\left\{\begin{array}{ll}
\D  D_t\cP_{ \mathfrak{J},v}=- D_ttr(\na\mathfrak{J}\na v)+2tr(\na v \na^2 \cP_{ \mathfrak{J},v})\qquad\hbox{in}\quad\Om_t\\
\na_{n_t} D_t\cP_{\mJ,v}|_{\G_t}= C'_{ \mathfrak{J}, v}-D_t(\na_{\tau_t}\mathfrak{K}_{\cH}\cdot  \na_{\tau_t} v \cdot n_t)+\na_{n_t}v\cdot \na \cP_{ \mathfrak{J},v},\\
 \na_{n_b}  D_tP_{ \mathfrak{J},v}|_{\G_b}=D_t( \mJ\cdot\na_v n_b)+\na_{n_b}v\cdot \na P_{ \mathfrak{J},v}.
\end{array}\right.
\eeq
Therefore, we have directly
\[
\begin{split}
&\int_{\Om_t} \na D_t \cP_{\mathfrak{J}, v} \cdot D_t \mathfrak{J}dX= \int_{\Om_t} \na D_t \cP_{\mathfrak{J}, v}\cdot  D_t \na \mathfrak{K}_{\cH}dX\\
 &= \int_{\Om_t} \na D_t \cP_{\mathfrak{J}, v}\cdot \na D_t \mathfrak{K}_{\cH}dX+\int_{\Om_t} \na D_t \cP_{\mathfrak{J}, v}\cdot \na v\cdot\na \mathfrak{K}_{\cH}dX\\
 &= \int_{\Om_t} \na D_t \cP_{\mathfrak{J}, v}\cdot \na D_t \mathfrak{K}_{\cH}dX+\f d{dt}\int_{\Om_t} \na  \cP_{\mathfrak{J}, v}\cdot \na v\cdot\na \mathfrak{K}_{\cH}dX-\int_{\Om_t} \na  \cP_{\mathfrak{J}, v}\cdot D_t(\na v\cdot\na \mathfrak{K}_{\cH})dX\\
 &\quad
 +\int_{\Om_t}\na v\cdot\na \cP_{\mJ,v}\cdot\na v\cdot\na\mK_\cH dX,
\end{split}
\]
which leads to the following equality:
\[
\begin{split}
&\int_{\Om_t} \na D_t \cP_{\mathfrak{J}, v} \cdot D_t \mathfrak{J}dX-\f d{dt}\int_{\Om_t} \na  \cP_{\mathfrak{J}, v}\cdot \na v\cdot\na \mathfrak{K}_{\cH}dX\\
&=\int_{\Om_t} \na D_t \cP_{\mathfrak{J}, v}\cdot \na D_t \mathfrak{K}_{\cH}dX-\int_{\Om_t} \na  \cP_{\mathfrak{J}, v}\cdot D_t(\na v\cdot\na \mathfrak{K}_{\cH})dX+\int_{\Om_t}\na v\cdot\na \cP_{\mJ,v}\cdot\na v\cdot\na\mK_\cH dX.
\end{split}
\]
To finish the proof, the key lies in the analysis for the first integral. The remainder part is controlled by $P(E(t))$  thanks to Lemma \ref{prop: eta}, Lemma \ref{prop: DtPvv}, Remark \ref{est:D_t k} and  Lemma \ref{est: P_{v, v}-H}.
\medskip

For the first integral on the right side of the equality above, we have by Green's Formula that
\beq\label{first integral}
\begin{split}
 &\int_{\Om_t} \na D_t \cP_{\mathfrak{J}, v}\cdot \na D_t \mathfrak{K}_{\cH}dX\\
 &=\int_{\G_t}(\na_{n_t} D_t \cP_{\mathfrak{J}, v})  D_t \mathfrak{K}_{\cH}ds+\int_{\G_b}(\na_{n_b} D_t\cP_{\mathfrak{J}, v})  D_t \mathfrak{K}_{\cH}ds-\int_{\Om_t}\Delta D_t\cP_{\mJ,v} D_t \mathfrak{K}_{\cH}dX\\
&=\int_{\G_t}\big[C'_{\mJ, v}-D_t(\na_{\tau_t}\mathfrak{K}_{\cH}\cdot  \na_{\tau_t} v \cdot n_t)+\na_{n_t}v\cdot \na \cP_{ \mathfrak{J},v}\big] D_t \mathfrak{K}_{\cH}ds
+\int_{\G_b}D_t( \mJ\cdot\na_v n_b)D_t \mathfrak{K}_{\cH}ds\\
&\quad +\int_{\G_b}\na_{n_b}v\cdot \na P_{ \mathfrak{J},v}\,D_t \mathfrak{K}_{\cH}ds
+\int_{\Om_t} D_ttr(\na\mathfrak{J}\na v)\,D_t \mathfrak{K}_{\cH}dX-\int_{\Om_t}2tr(\na v \na^2 \cP_{ \mathfrak{J},v})\,D_t \mathfrak{K}_{\cH}dX\\
&=C'_{\mJ, v}\int_{\G_t}D_t \mathfrak{K}ds-\int_{\G_t}\big(\na_{\tau_t}D_t\mathfrak{K}\cdot  \na_{\tau_t} v \cdot n_t\big) D_t \mathfrak{K}ds+\int_{\G_b}(D_t\mJ\cdot\na_v n_b) D_t \mathfrak{K}_{\cH}ds\\
&\quad+\int_{\Om_t}tr(\na D_t\mJ\na v)D_t\mK_\cH dX +\hbox{l.o.t.},
\end{split}
\eeq
where the remainder lower-order terms can be controlled by $P(E(t))$ in a similar way as before.
\medskip

Now we deal with the terms in \eqref{first integral} one by one. To begin with, we write from the definition of $C_{\mJ,v}$ in system \eqref{PJv system} that
\[
\begin{split}
\f d{dt}\big(|\G_t|C_{\mJ,v}\big)&=-\int_{\Om_t}D_t tr(\na\mJ\na v)dX+\int_{\G_t}D_t(\na_{\tau_t}\mK\na_{\tau_t}v\cdot n_t ds)
-\int_{\G_b}D_t(\mJ\cdot \na_v n_b ds)\\
&=-\int_{\Om_t}tr(\na D_t\mJ\na v)dX+\int_{\G_t}\na_{\tau_t}D_t\mK\na_{\tau_t}v\cdot n_tds-\int_{\G_b}D_tJ\cdot\na_v n_bds+\hbox{l.o.t.},
\end{split}
\]
where the lower-order terms can be controlled again and hence the details are omitted. Applying Green's Formula again and using the decompositions
\[
\begin{split}
&\na D_t\mK_\cH=(\na_{\tau_t}D_t\mK)\tau_t+(\na_{n_t}D_t\mK_\cH)n_t\qquad\hbox{on}\quad \G_t,\\
&\na D_t\mK_\cH=(\na_{\tau_b}D_t\mK_\cH)\tau_b+(\na_{n_b}D_t\mK_\cH)n_b\qquad\hbox{on}\quad \G_b,
\end{split}
\]
we find
\beq\label{pt CJv}
\begin{split}
&\f d{dt}\big(|\G_t|C_{\mJ,v}\big)\\
&=-\int_{\G_t}D_t\mJ\cdot\na v\cdot n_tds-\int_{\G_b}D_t\mJ\cdot\na v\cdot n_bds+\int_{\G_t}\na_{\tau_t}D_t\mK\na_{\tau_t}v\cdot n_tds-\int_{\G_b}D_tJ\cdot\na_v n_bds+\hbox{l.o.t.}\\
&=-\int_{\G_t}\na D_t\mK_\cH\cdot\na v\cdot n_tds+\int_{\G_t}\na_{\tau_t}D_t\mK\na_{\tau_t}v\cdot n_tds-\int_{\G_b}\na D_t\mK_\cH\cdot\na v\cdot n_bds\\
&\quad
-\int_{\G_b}\na D_t\mK_\cH\cdot\na_v n_bds+\hbox{l.o.t.}\\
&=-\int_{\G_t}\na_{n_t}D_t\mK_\cH\na_{n_t}v\cdot n_tds-\int_{\G_b}\na_{n_b}D_t\mK_\cH\na_{n_b}v\cdot n_bds-\int_{\G_b}\na_{\tau_b}D_t\mK\cdot\na_{\tau_b} v\cdot n_bds
\\
&\quad-\int_{\G_b}\na D_t\mK_\cH\cdot\na_v n_bds+\hbox{l.o.t.}.
\end{split}
\eeq
Here, for the last three terms on $\G_b$, thanks to the assumption that $n_b$ is constant near the contact points $p_i$ and $v\cdot n_b|_{\G_b}=0$, we know that $\na_{\tau_b} v\cdot n_b=-v\cdot\na_{\tau_b}n_b$ and $\na_v n_b=(\na_v n_b\cdot \tau_b)\tau_b$ vanish near $p_i$. Moreover, we also know from the definition of $\mK_\cH$ that
\[
\na_{n_b}D_t\mK_\cH|_{\G_b}=[\na_{n_b}, D_t]\mK_\cH=\na_{n_b}v\cdot \na\mK_\cH.
\]
Consequently, applying Lemma \ref{H1 zero trace} and Lemma \ref{G_b H-1/2}, we have
\[
\begin{split}
&\int_{\G_b}\na_{\tau_b}D_t\mK_\cH\cdot\na_{\tau_b} v\cdot n_bds+\int_{\G_b}\na D_t\mK_\cH\cdot\na_v n_bds\\
&\le \|\na_{\tau_b}D_t\mK_\cH\|_{\tilde H^{-1/2}(\G_b)}\big(\|v\cdot\na_{\tau_b}n_b\|_{\tilde H^{1/2}(\G_b)}+\|\na_v n_b\cdot \tau_b\|_{\tilde H^{1/2}(\G_b)}\big)\\
&\le C(\|\G_t\|_{H^{5/2}})\|D_t\mK_\cH\|_{H^1(\Om_t)}\big(\|v\cdot\na_{\cH(\tau_b)}\cH(n_b)\|_{H^1(\Om_t)}+\|\na_v \cH(n_b)\cdot \cH(\tau_b)\|_{H^1(\Om_t)}\big)\le P(E(t))
\end{split}
\]
and
\[
\int_{\G_b}\na_{n_b}D_t\mK_\cH\na_{n_b}v\cdot n_bds\le P(E(t)).
\]

For the first integral in \eqref{pt CJv}, applying Green's Formula again leads to
\[
\int_{\G_t}\na_{n_t}D_t\mK_\cH\na_{n_t}v\cdot n_tds=\int_{\Om_t}[\Delta, D_t]\mK_\cH\na_{\cH(n_t)}v\cdot \cH(n_t)ds-\int_{\G_b}[\na_{n_b},D_t]\mK_\cH\na_{\cH(n_t)}v\cdot \cH(n_t)ds,
\]
so this term is  controlled by $P(E(t))$ again.

Therefore, summing up all these estimates above and going back to \eqref{pt CJv}, we arrive at
\[
\f d{dt}\big(|\G_t|C_{\mJ,v}\big)= \hbox{l.o.t.}
\]
with the lower-order terms controlled by $P(E(t))$. This leads to the following estimate
\[
C'_{\mJ,v}\le P(E(t)).
\]
\medskip

For the moment, we can go back and deal with the other terms in \eqref{first integral}. In fact, using similar arguments as above, we also conclude that
\[
\begin{split}
\int_{\G_t}\big(\na_{\tau_t}D_t\mathfrak{K}_{\cH}\cdot  \na_{\tau_t} v \cdot n_t\big) D_t \mathfrak{K}_{\cH}ds+\int_{\G_b}(D_t\mJ\cdot\na_v n_b) D_t \mathfrak{K}_{\cH}ds-\int_{\Om_t}tr(\na D_t\mJ\na v)D_t\mK_\cH dX\le P(E(t))
\end{split}
\]
Therefore,  the proof is finished.

\ef

\subsection{Boundary terms at the contact points}
 \begin{lem}\label{lem:cc}
We have the following equation on $\G_b$:
\beq\label{equ:cc}
 D_t\mJ = -\f{\sigma^2 }{\beta_c}(n_t\cdot \tau_b) (\na_{\tau_t} \mJ)^\bot
 \tau_b+R_{c1},
\eeq
and there holds at the contact points $p_i$ ($i=l, r$) that
\[
(D_t\mJ)^\perp\,( \na_{{\bf \tau}_t}\mJ)^\perp \big|_{p_l}=-\f{\sigma^2}{\beta_c}F_l(t) +R_{c2,l},\quad (D_t\mJ)^\perp\,( \na_{{\bf \tau}_t}\mJ)^\perp \big|_{p_r}=\f{\sigma^2}{\beta_c}F_r(t) + R_{c2,r}
\]
where
\[
\begin{split}
&R_{c1}=r_c\,\tau_b-  ( \mJ \cdot  D_t n_b)\,n_b,\quad
R_{c2,i}= \big(-r_c +\cot \om(\mJ \cdot  D_t n_b)\big)\, (\sin\om) \na_{{\bf \tau}_t}\mJ^\perp\big|_{p_i},\qquad\hbox{and}\\
&r_c=-\sigma  \sin\om (\na_{\tau_t}v\cdot D_t n_t)+\sigma  \tau_b\cdot D_t n_t\,(\na_{\tau_t}v\cdot n_t)-\sigma  \sin\om  \big(\na P_{v,v}\cdot\na_{\tau_t}n_t +[D_t, \na_{\tau_t} ] v\cdot n_t\big)\\
&\qquad
- \beta_c D_t \na P_{v,v}\cdot \tau_b
\end{split}
\]
satisfy
\[
|R_{c1}|\le P(E(t)),\quad |R_{c2,i}|\le P(E(t))F(t)^{1/2}.
\]
 \end{lem}
\no{\bf Proof.} The proof follows the proof of Lemma 7.1 in \cite{MW2}. For the estimates, we apply Lemma \ref{embedding}, Lemma \ref{nt H3 estimate}, Lemma \ref{na v estimate} and  Lemma \ref{prop: D_t P_{v, v}-H}. Here the differences compared to the proof of Lemma 7.1 in \cite{MW2} lie in that  we have different $\mJ$  and
\[
\na_{\tau_t}\na P_{v,v}\cdot n_t=-\na P_{v,v}\cdot\na_{\tau_t}n_t
\]
in $r_c$ due to the definition of  $P_{v,v}$.

\ef

\medskip

\subsection{Proof of Theorem \ref{lower-order energy estimates}}
We are ready to show the a priori energy estimate in Theorem \ref{lower-order energy estimates}. In fact,  taking $L^2(\Om_t)$ inner product with $\cD_t J-\na \cH(  D_t P_{v,v} - v\cdot \big(\na P_{v,v}|_{c})\big)$ on both sides of \eqref{eqn cDt2 J}, we have
\beq\label{equ:Ee2}
\begin{split}
&\int_{\Om_t} D_t \big[ \cD_t \mJ-\na \cH\big(  D_t P_{v,v} - v\cdot (\na P_{v,v}|_{c})\big)\big]\cdot\big[ \cD_t \mJ-\na \cH\big(  D_t P_{v,v} - v\cdot (\na P_{v,v}|_{c})\big)\big]dX \\
&+\sigma\int_{\Om_t} \cA \mJ\cdot\big[ \cD_t \mJ-\na \cH\big(  D_t P_{v,v} - v\cdot (\na P_{v,v}|_{c})\big)\big]dX =\int_{\Om_t}\cR\cdot \big[ \cD_t \mJ-\na \cH\big(  D_t P_{v,v} - v\cdot (\na P_{v,v}|_{c})\big)\big]dX.
\end{split}
\eeq

 \medskip
We deal with these integrals above one by one. To get started, for the first term on the left side, we rewrite it as
\[
\begin{split}
&\int_{\Om_t} D_t \big[ \cD_t \mJ-\na \cH\big(  D_t P_{v,v} - v\cdot (\na P_{v,v}|_{c})\big)\big]\cdot\big[ \cD_t \mJ-\na \cH\big(  D_t P_{v,v} - v\cdot (\na P_{v,v}|_{c})\big)\big]dX\\
&=\f12 \f d{dt}\int_{\Om_t} \big| \cD_t \mJ-\na \cH\big(  D_t P_{v,v} - v\cdot (\na P_{v,v}|_{c})\big)\big|^2dX .
\end{split}
\]
For the second term on the left side of  \eqref{equ:Ee2},  we have by Green's Formula  that
\beq\label{equ:Ee1}
\begin{split}
&\int_{\Om_t} \cA \mJ\cdot\big[ \cD_t \mJ-\na \cH\big(  D_t P_{v,v} - v\cdot (\na P_{v,v}|_{c})\big)\big]dX\\
&=-\int_{\Gamma_t} \Delta_{\Gamma_t} \mJ^\bot \,  \cD_t \mJ\cdot n_tds +\int_{\Gamma_t} \Delta_{\Gamma_t} \mJ^\bot \, \na \cH\big(  D_t P_{v,v} - v\cdot (\na P_{v,v}|_{c})\big)\cdot n_tds.
\end{split}
\eeq
For the first term on the right side of \eqref{equ:Ee1}, one deduces from Hodge decomposition \eqref{Hodge decom} and integration by parts as in \cite{MW3} that
\[
\begin{split}
&-\int_{\Gamma_t} \Delta_{\Gamma_t} \mJ^\bot \,  \cD_t \mJ\cdot n_tds\\
&=\int_{\Gamma_t} \na_{\tau_t} \mJ^\bot \, \na_{\tau_t}(D_t \mJ\cdot n_t )ds-(D_t \mJ)^\bot \na_{\tau_t} \mJ^\bot \big|^{p_l}_{p_r} -\int_{\Gamma_t} \Delta_{\Gamma_t} \mJ^\bot  (\na  \cP_{\mJ,v}\cdot n_t)ds\\
&=\f12\,\f d{dt}\int_{\Gamma_t}| \na_{\tau_t} \mJ^\bot|^2ds-\int_{\Gamma_t} \na_{\tau_t} \mJ^\bot \, [  D_t, \na_{\tau_t}] \mJ^\bot ds- \int_{\Gamma_t} \na_{\tau_t} \mJ^\bot \, \na_{\tau_t}(\mJ\cdot  D_t  n_t )ds
+\f{\sigma^2}{\beta_c}F(t) \\
&\quad -R_{c2,l}+R_{c2,r}-\int_{\Gamma_t} \Delta_{\Gamma_t} \mJ^\bot  (\na  \cP_{\mJ,v}\cdot n_t)ds,
\end{split}
\]
where we apply Lemma \ref{lem:cc} at corner points.

For the second term in \eqref{equ:Ee1}, one has in a similar way as above that
\[
\begin{split}
&\int_{\Gamma_t} \Delta_{\Gamma_t} \mJ^\bot \, \na \cH\big(  D_t P_{v,v} - v\cdot (\na P_{v,v}|_{c})\big)\cdot n_tds\\
&=-\int_{\Gamma_t}\na_{\tau_t} J^\bot  \cdot   \na_{\tau_t}\na_{n_t}\cH\big(  D_t P_{v,v} - v\cdot (\na P_{v,v}|_{c})\big)ds +\na_{\tau_t} J^\bot  \cdot    \na_{n_t}\cH\big(  D_t P_{v,v} - v\cdot (\na P_{v,v}|_{c})\big)\big|^{p_r}_{p_l}\\
&\leq P(E(t))\big(1+ F(t)^{1/2}\big),
\end{split}
\]
where Lemma \ref{embedding}  and Lemma \ref{prop: D_t P_{v, v}-H} are used.

\medskip
As a result, going back to \eqref{equ:Ee2}, we summerise that
\[
\begin{split}
&\f12 \f d{dt}\int_{\Om_t} \big| \cD_t \mJ-\na \cH\big(  D_t P_{v,v} - v\cdot (\na P_{v,v}|_{c})\big)\big|^2dX+\f\si2\,\f d{dt}\int_{\Gamma_t}| \na_{\tau_t} \mJ^\bot|^2ds
+ \f{\sigma^3}{2\beta_c}  F(t)\\
&\le P(E(t))+\si\int_{\Gamma_t} \na_{\tau_t} \mJ^\bot \, [  D_t, \na_{\tau_t}] \mJ^\bot ds+\si\int_{\Gamma_t} \na_{\tau_t} \mJ^\bot \, \na_{\tau_t}(\mJ\cdot  D_t  n_t )ds
+\si\int_{\Gamma_t} \Delta_{\Gamma_t} \mJ^\bot  (\na  \cP_{\mJ,v}\cdot n_t)ds\\
&\quad +\int_{\Om_t}\cR\cdot \big[ \cD_t \mJ-\na \cH\big(  D_t P_{v,v} - v\cdot (\na P_{v,v}|_{c})\big)\big]dX.
\end{split}
\]
Moreover, direct estimates similarly as in \cite{MW3} lead to the following estimate
\[
\begin{split}
&\si\int_{\Gamma_t} \na_{\tau_t} \mJ^\bot \, [  D_t, \na_{\tau_t}] \mJ^\bot ds+\si\int_{\Gamma_t} \na_{\tau_t} \mJ^\bot \, \na_{\tau_t}(\mJ\cdot  D_t  n_t )ds
+\si\int_{\Gamma_t} \Delta_{\Gamma_t} \mJ^\bot  (\na  \cP_{\mJ,v}\cdot n_t)ds\\
&\le P(E(t))+\f{\sigma^3}{4\beta_c}  F(t),
\end{split}
\]
so we obtain
\beq\label{energy esti-1}
\begin{split}
&\f12 \f d{dt}\Big(\int_{\Om_t} \big| \cD_t \mJ-\na \cH\big(  D_t P_{v,v} - v\cdot (\na P_{v,v}|_{c})\big)\big|^2dX
+\si\int_{\Gamma_t}| \na_{\tau_t} \mJ^\bot|^2ds\Big)
+ \f{\sigma^3}{4\beta_c}  F(t)\\
&\le P(E(t))
+\int_{\Om_t}\cR\cdot \big[ \cD_t \mJ-\na \cH\big(  D_t P_{v,v} - v\cdot (\na P_{v,v}|_{c})\big)\big]dX.
\end{split}
\eeq

In order to finish the energy estimate, we still need to deal with the integral on the right side involving the remainder term $\cR$ in \eqref{eqn cDt2 J}.

\medskip

\noindent - Estimates for the part $R_0-\si\na\cH(R_1)$ in $\cR$.
\begin{lem}\label{prop:R_0-R_1}
For the remainder term $R_0$ defined in \eqref{eqn for Dt2 J}, we have the estimate
\[
\|R_0-\si\na\cH(R_1) \|_{L^2(\Om_t)}\leq  P(E(t)) .
\]
\end{lem}
\no{\bf Proof.}
Recalling from \eqref{eqn for Dt2 J}, we know that
\[
R_0-\si\na\cH(R_1) =-\s\na \cH(J\cdot \D_{\G_t}n_t)+\si\na \cH(n_t\cdot \D_{\G_t}\na P_{v,v}) +[D_t, \na \cH](  D_t P_{v,v}  ) +A_1+A_2+A_3
\]
where $R_1$ and $A_1,A_2,A_3$ are defined in \eqref{Dt2 kappa} and \eqref{A_1}, \eqref{A_2}, \eqref{Dt2 J eqn1} respectively.

For the term $\na \cH(n_t\cdot \D_{\G_t}\na P_{v,v})$, we have
\[
\|\na \cH(n_t\cdot \D_{\G_t}\na P_{v,v})\|_{L^2(\Om_t)}=\|\na \cH([n_t, \D_{\G_t}]\na P_{v,v})\|_{L^2(\Om_t)}\le  P(E(t)),
\]
where the boundary condition $\na_{n_t}P_{v,v}|_{\G_t}=C_{v,v}(t)$ and Lemma \ref{est: P_{v, v}-H} are applied.

For the term $[D_t, \na \cH](  D_t P_{v,v}  )$, direct computations using \eqref{commutator Dt H} and similar arguments as in the proof of Lemma \ref{prop: DtPvv} lead to
\[
\|[D_t, \na \cH](  D_t P_{v,v}  )\|_{L^2(\Om_t)}\leq\|\na v\cdot \na \cH(  D_t P_{v,v}  )\|_{L^2(\Om_t)}+\| [D_t, \cH](  D_t P_{v,v}  )\|_{H^1(\Om_t)}\leq P(E(t)).
\]

The estimates for the other terms follow from lemmas in the previous section and can be done similarly as \cite{MW3}, so we omit the details here.

\ef

\noindent - The part $\si\na\cH(R_1 )$ in $\cR$. In fact,  the following integral from the right side of \eqref{energy esti-1} is rewritten as
\[
\begin{split}
&\int_{\Om_t}\na\cH(R_1 )\cdot \big[ \cD_t \mJ-\na \cH\big(  D_t P_{v,v} - v\cdot (\na P_{v,v}|_{c})\big)\big]dX\\
&= \int_{\G_t}R_1\, \cD_t \mJ\cdot n_t ds-\int_{\G_t}R_1 \na_{n_t} \cH\big(  D_t P_{v,v} - v\cdot (\na P_{v,v}|_{c})\big) ds\\
&=\int_{\G_t}R_1\, D_t \mJ\cdot n_t ds+\int_{\G_t}R_1 \,\na \cP_{\mJ,v}\cdot n_t ds-\int_{\G_t}R_1 \na_{n_t} \cH\big(  D_t P_{v,v} - v\cdot (\na P_{v,v}|_{c})\big) ds
\end{split}
\]

On the other hand, noticing  from  \eqref{Dt2 kappa} that $R_1$ contains terms like $\na^2 v ,\na n_t, \kappa$ and $\na^2 P_{v, v}$ and using \eqref{na v decomp} and lemmas in the previous subsections,  we have
\[
\int_{\G_t}R_1 \na_{n_t} \cH\big(  D_t P_{v,v} - v\cdot (\na P_{v,v}|_{c})\big) ds
\le \int_{\G_t}|R_1 |ds\, \|\na_{n_t} \cH\big(  D_t P_{v,v} - v\cdot (\na P_{v,v}|_{c})\big)\|_{L^\infty(\G_t)}\le P(E(t)).
\]
The details in the estimate above are omitted and we only note that for $\pa^2 v$ terms in $R_1$, we can have the following estimate thanks to \eqref{na v decomp}:
\[
\int_{\G_t}|\na^2 v|ds\le \int_{\G_t}|\na^3\phi_r|ds+\int_{\G_t}|\na^3\phi_s|ds\le P(E(t)).
\]
Similarly but more easily, we also  obtain
\[
\int_{\G_t}R_1 \,\na \cP_{\mJ,v}\cdot n_t ds\le P(E(t)).
\]

Besides, for the part $\int_{\G_t}R_1 D_t \mJ\cdot n_t ds$, we put $D_t$ out of the integral:
\[
\begin{split}
&\int_{\G_t}R_1 D_t \mJ\cdot n_t ds=\int_{\G_t}R_1 D_t(\mJ-\mJ|_c)ds+\int_{\G_t}R_1 D_t(\mJ|_c)ds\\
&=\f d{dt}\int_{\G_t} R_1\,(\mJ-\mJ|_c)ds-\int_{\G_t}D_t R_1\,(\mJ-\mJ|_c)ds-\int_{\G_t}R_1\,(\mJ-\mJ|_c)D_tds+\int_{\G_t}R_1 D_t(\mJ|_c)ds.
\end{split}
\]
Here we need to take care of the terms in $\int_{\G_t}D_t R_1\,(\mJ-\mJ|_c)ds$ and $\int_{\G_t}R_1\,(\mJ-\mJ|_c)D_tds$. In fact, similarly to the analysis in \eqref{apply hardy}, the key terms like $\|\pa^2 v (\mJ-\mJ|_c)\|_{L^2(\G_t)}$ can be handled as follows:
\beq\label{pa2 v estimate}
\|\pa^2 v (\mJ-\mJ|_c)\|_{L^2(\G_t)}\le
C(\|\G_t\|_{H^{5/2}})\|r^{\delta} (\pa^2 v)\circ T^{-1}_i \|_{L^\infty(\G_t)}\| r^{-\delta}(\mJ-\mJ|_{p_i})\circ T^{-1}_i\|_{L^2(\G_t)}\le P(E(t)),
\eeq
where the constant $\delta\in (0,1)$ is chosen to satisfy 
\[
\delta+\pi/\om_i-3>0,\quad\hbox{ when}\quad \om_i\in (\pi/3, \pi/2).
\]

Consequently, we derive
\[
\int_{\G_t}R_1 D_t \mJ\cdot n_t ds=\f d{dt}\int_{\G_t} R_1 (\mJ-
\mJ|_c)\cdot n_t ds+l.o.t.,
\]
where all the lower-order terms are controlled by $P(E(t))$.

As a result, we conclude that
\[
\int_{\Om_t}\na\cH(R_1 )\cdot \big[ \cD_t \mJ-\na \cH\big(  D_t P_{v,v} - v\cdot (\na P_{v,v}|_{c})\big)\big]dX
=\f d{dt}\int_{\G_t} R_1 (\mJ-
\mJ|_c)\cdot n_t ds+l.o.t.
\]
with the lower-order terms controlled by $P(E(t))$.

\bigskip

\noindent - The terms $D_t\na \cP_{\mJ,v}+D_t\na\cH\big(v\cdot(\na P_{v,v}|_c)\big)$ in $\cR$.
First,  one has directly from Lemma \ref{lem: D_t P_{J, v}} that
\[
\int_{\Om_t} \na D_t\cP_{\mathfrak{J}, v} \cdot D_t \mathfrak{J}dX=\f d{dt}\int_{\Om_t} \na  \cP_{\mathfrak{J}, v}\cdot \na v\cdot\na \mathfrak{K}_{\cH}dX +l.o.t.
\]
where the lower-order terms are all controlled by $P(E(t))$.

Second, for the integral
\[
\int_{\Om_t} \na D_t\cP_{\mathfrak{J}, v} \cdot \na \cH\big(  D_t P_{v,v} - v\cdot (\na P_{v,v}|_{c})\big)dX,
\]
applying Green's formula and similar calculations as in \eqref{first integral} show directly that it can be controlled by $P(E(t))$, and we omit the details.

Therefore, we conclude that
\[
\int_{\Om_t} \na D_t\cP_{\mathfrak{J}, v}
\cdot \big[ \cD_t \mJ-\na \cH\big(  D_t P_{v,v} - v\cdot (\na P_{v,v}|_{c})\big)\big]dX
=\f d{dt}\int_{\Om_t} \na  \cP_{\mathfrak{J}, v}\cdot \na v\cdot\na \mathfrak{K}_{\cH}dX+l.o.t.,
\]
where all the lower-order terms are controlled by $P(E(t))$.

\bigskip
Summing up all these estimates related to $\cR$ above and going back to \eqref{energy esti-1}, we obtain the following estimate
\[
\begin{split}
&\f12 \f d{dt}\Big(\int_{\Om_t} \big| \cD_t \mJ-\na \cH\big(  D_t P_{v,v} - v\cdot (\na P_{v,v}|_{c})\big)\big|^2dX
+\si\int_{\Gamma_t}| \na_{\tau_t} \mJ^\bot|^2ds\Big)
+ \f{\sigma^3}{4\beta_c}  F(t)\\
&\le P(E(t))+\f d{dt}\int_{\G_t} R_1 (\mJ-
\mJ|_c)\cdot n_t ds
+\f d{dt}\int_{\Om_t} \na  \cP_{\mathfrak{J}, v}\cdot \na v\cdot\na \mathfrak{K}_{\cH}dX.
\end{split}
\]
Integrating on both sides with respect to time  on $[0, t]$, we have
\[
\begin{split}
&\f12\int_{\Om_t} \big| \cD_t \mJ-\na \cH\big(  D_t P_{v,v} - v\cdot (\na P_{v,v}|_{c})\big)\big|^2dX
+\f\si 2\int_{\Gamma_t}| \na_{\tau_t} \mJ^\bot|^2ds
+ \f{\sigma^3}{4\beta_c} \int^t_0 F(t')dt'\\
&\le \f12\int_{\Om_0} \big| \cD_t \mJ-\na \cH\big(  D_t P_{v,v} - v\cdot (\na P_{v,v}|_{c})\big)\big|_{t'=0}\big|^2dX
+\f\si 2\int_{\Gamma_{t0}}| \na_{\tau_t} \mJ^\bot|^2\big|_{t'=0}ds\\
&\quad+\int^t_0 P(E(t'))dt'+\int_{\G_{t'}} R_1 (\mJ-
\mJ|_c)\cdot n_tds\Big|^t_0
+\int_{\Om_{t'}} \na  \cP_{\mathfrak{J}, v}\cdot \na v\cdot\na \mathfrak{K}_{\cH}dX\Big|^t_0\\
&\le P(E(0))+\int^t_0 P(E(t'))dt'+\int_{\G_{t'}} R_1 (\mJ-
\mJ|_c)\cdot n_tds\Big|^t_0
+\int_{\Om_{t'}} \na  \cP_{\mathfrak{J}, v}\cdot \na v\cdot\na \mathfrak{K}_{\cH}dX\Big|^t_0.
\end{split}
\]
Replacing the first integral with $\|\cD_t \mJ\|^2_{L^2(\Om_t)}$, we arrive at the following inequality:
\[
\begin{split}
E(t)+ \int^t_0 F(t')dt'
&\le P(E(0))+\int^t_0 P(E(t'))dt'+C\int_{\G_{t'}} R_1 (\mJ-
\mJ|_c)\cdot n_tds\Big|^t_0\\
&\quad
+C\int_{\Om_{t'}} \na  \cP_{\mathfrak{J}, v}\cdot \na v\cdot\na \mathfrak{K}_{\cH}dX\Big|^t_0
+C\int_{\Om_{t}} \big| \na \cH\big(  D_t P_{v,v} - v\cdot (\na P_{v,v}|_{c})\big)\big|^2dX,
\end{split}
\]
where the constant $C$ depends on $\si, \beta_c$.

We now deal with the last three integrals on the right side above one by one. First, for the terms in $R_1$, similarly as in \eqref{pa2 v estimate}, we can have by careful estimates with interpolations that
\[
\|r^{\delta}( \pa^2 v)\circ T^{-1}_i \|_{L^\infty(0, r_0)}\le P(E_l(t))E(t)^{1/2}
\]
and
\[
\|r^{-\delta}(\mJ-\mJ|_{p_i})\circ T^{-1}_i\|_{L^2(0, r_0)}\le C(\|\G_t\|_{H^{5/2}})\big\|\mJ-\mJ|_c\big\|_{H^\delta(\G_t)}\le C(\|\G_t\|_{H^{5/2}})\|\mJ\|_{H^{\delta+1/2}(\Om_t)},
\]
where $1>\delta>3-\pi/\om$ when $\om\in (\pi/3, \pi/2)$.

Therefore, we can prove by interpolations that there exist $\delta_0\in (1/2, 1)$ and $\epsilon$ small enough  such that
\[
\int_{\G_t} R_1 (\mJ-\mJ|_c)\cdot n_tds\le P(E_l(t))E(t)^{\delta_0}\le \epsilon E(t)+C_{\epsilon, \delta_0} P(E_l(t)).
\]

Moreover, the following estimates can also be proved in a similar but easier way:
\[
\int_{\Om_t} \na  \cP_{\mathfrak{J}, v}\cdot \na v\cdot\na \mathfrak{K}_{\cH}dX
+\int_{\Om_t} \big| \na \cH\big(  D_t P_{v,v} - v\cdot (\na P_{v,v}|_{c})\big)\big|^2dX
\le \epsilon E(t)+C_{\epsilon, \delta_0} P(E_l(t)).
\]

\bigskip
As a result, summing up the estimates above, we are able to conclude that
\[
E(t)+ \int^t_0 F(t')dt'\le P(E(0))+\int^t_0 P(E(t'))dt' +P(E_l(t)).
\]

In the end, to close the energy estimates, we  deal with $E_l(t)$.  First, we show by direct calculations using Euler's equation and lemmas from the previous section that
\[
\|v\|^2_{H^k(\Om_t)}\le C(T) \|v(0)\|^2_{H^k(\Om_0)}+C(T)\|\mJ\|^2_{L^2([0,T],H^k(\Om_t))}+\int^t_0P(E(t'))dt'
\]
for $k=1,2$.
Then we take the square root on both sides of the inequalities above and apply an interpolation between $H^1(\Om_t), H^2(\Om_t)$ as well as Lemma \ref{prop: eta} to obtain
\[
\|v\|_{H^{3/2}(\Om_t)}\le C(T) \|v(0)\|_{H^{3/2}(\Om_0)}+\int^t_0P(E(t'))dt'.
\]
Second,  the estimate for $\|\G_t\|_{H^{5/2}}$ can be derived in a similar way as above.
Consequently, we have
\[
E_l(t)\le P(E(0))+\int^t_0P(E(t'))dt'.
\]

Combining this estimate with the estimate above for $E(t)$, we  finish the lower-order energy estimates.

 \section{Higher-order time-derivative energy estimates}
 \setcounter{equation}{0}

 In this section, we prove  higher-order energy estimates with respect to $D_t$, which is needed in the local well-posedness part.  To start with, we define the  energy functional
\[
E_1(t)=\|\na_{\tau_t} D_t  \mathfrak{J}^\bot \|^2_{L^2(\Gamma_t)} +\|D^2_t \mathfrak{J}\|^2_{L^2(\Om_t)},
\]
and the dissipation
\[
F_1(t)=\sum_{i=l,r}\big|(\sin \om_i)\na_{{\bf \tau}_t}D_t \mathfrak{J}^\perp |_{p_i}\big|^2.
\]

The main result of this section is as follows:
\begin{thm}\label{1-order energy estimates}
Let the contact angles $\om_i\in (0,\pi/2)$ and  $E(t), \int_0^T F(t)dt, E_1(t), \int_0^T F_1(t)dt$ be bounded above in $[0, T]$ for some  $T>0$. Then the following higher-order a priori estimate holds
\[
\sup_{0\leq t\leq T} E_1(t)+\int_0^T F_1(t)dt\le P(E_1(0))+\int_0^T P(E_1(t))dt.
\]

\end{thm}
To prove this theorem, we begin with  a higher-order equation for $\mJ$ and more delicate estimates involving $D_t$ based on Section \ref{estimate section}. With these preparations, we are able to finish the energy estimate in the last subsection of this part.

\bigskip

 \subsection{The higher-order equation for $\mJ$}
 Firstly, we recall system \eqref{eqn for Dt2 J} and rewrite it as follows:
\beq\label{rewrite eqn for Dt2 J}
D^2_t \mathfrak{J}=\si\na \cH(\D_{\G_t}  \mathfrak{J}^\perp-h_v)+\tilde{R}_0,
\eeq
where
\[
\tilde{R}_0 =-\si\na \cH(\mathfrak{J} \cdot \D_{\G_t}n_t)+\si\na \cH([n_t, \D_{\G_t}]\cdot\na P_{v,v})+\s\na\cH(R_1+h_v ) +A_1+A_2+A_3+ \na \cH(  D^2_t P_{v,v}  ),
\]
and $h_v$ contains all the second-order terms of $v$ in $R_1$ which come from the commutator $[D_t, n_t\D_{\G_t}]\cdot v$.  More precisely, we have
\beq\label{hv def}
h_v= \D_{\G_t} v\cdot D_t n_t  +2D^2v\big(\tau_t, (\na_{\tau_t}v)^\top\big)\cdot n_t.
\eeq
As a result,  $R_1+h_v$ only contains lower-order derivatives like $\pa v, \pa n_t, \kappa$.

%
%
%

Acting $D_t$ on both sides of \eqref{eqn for Dt2 J}, we obtain
\beq\label{eq: D_t^2 J }
D_t (\cD^2_t \mathfrak{J}) +\na\cH\big(\D_{\Gamma_t}  D_t \mathfrak{J}^\bot-D_th_v\big)\\
= \tilde\cR_2,
\eeq
with the right side
\[
\begin{split}
\tilde\cR_2&=D_t \tilde{R}_0  +[\na, D_t]\cH \big(\D_{\Gamma_t} \mathfrak{J}^\bot-h_v)\big)  +\na [\cH, D_t]\big(\D_{\Gamma_t} \mathfrak{J}^\bot-h_v\big)+\na \cH [\D_{\Gamma_t}, D_t]\mathfrak{J}^\bot\\
&\quad+D_t\big(D_t\na P_{\mathfrak{J}, v}+\na \cP_{\na P_{\mathfrak{J}, v}, v}+\na P_{\cD_t \mathfrak{J}, v}\big).
\end{split}
\]
Here we use
\[
 -[\cA, D_t]=[\na, D_t]\cH \D_{\Gamma_t}+\na [\cH, D_t]\D_{\Gamma_t}+\na \cH [\D_{\Gamma_t}, D_t]
\]
with all these commutators from the end of Section 2,
and
\beq\label{Dt2 J decomp}
\cD^2_t \mathfrak{J}= D_t^2 \mathfrak{J}+D_t\na P_{\mathfrak{J}, v}+\na P_{\na P_{\mathfrak{J}, v}, v}+\na P_{\cD_t \mathfrak{J}, v},
\eeq
where $P_{w, v}$ is defined in \eqref{eq: cP_{w,v}}.

\subsection{More preliminary estimates.}\label{estimate section 2}
Based on Section \ref{estimate section}, we are going to prove higher-order estimates for different quantities using our higher-order energy.

\noindent - Higher-order estimates for $\mathfrak{K}_{\cH}$ and $D_t\mJ$.
 To get started, we recall system \eqref{kH system} to find
\[
\left\{\begin{array}{ll}
\D D_t\mathfrak{K}_\cH=2tr(\na v \na^2 \mathfrak{K}_\cH)\qquad\hbox{in}\quad\Om_t\\
\na_{n_t} D_t \mathfrak{K}_\cH|_{\G_t} = D_t \mJ^\perp +\na_{n_t} v \cdot \na \mathfrak{K}_\cH|_{\G_t},\qquad  \na_{n_b}D_t\mathfrak{K}_\cH|_{\G_b}=\na_{n_b}v\cdot \na\mK_\cH|_{\G_b}.
\end{array}\right.
\]
Applying Proposition \ref{elliptic estimate Neumann H3} and lemmas in Section \ref{estimate section}, we have
\beq\label{est:D_tk-1}
\begin{split}
\| D_t \mathfrak{K}_\cH\|_{H^2(\Om_t)}&\le C(\|\G_t\|_{H^{5/2}})\big(\|2tr(\na v \na^2 \mathfrak{K}_\cH)\|_{L^2(\Om_t)}+\| D_t \mJ^\perp +\na_{n_t} v \cdot \na \mathfrak{K}_\cH\|_{H^{1/2}(\G_t)}\\
&\quad
+\|\na_{n_b}v\cdot \na\mK_\cH\|_{H^{1/2}(\G_b)}+\|D_t\mK_\cH\|_{L^2(\Om_t)}\big)\\
&\leq P(E(t))(1+ \|D_t  \mathfrak{J}^\bot\|_{H^{1/2}(\G_t)}).
\end{split}
\eeq
 Here we notice that the regularity of $D_t\mathfrak{K}_\cH$ is constrained due to the regularity of $v$.  To get a higher-order estimate, we need to get rid of the worst part.

Similarly as in Lemma \ref{prop: D_t P_{v, v}-H}, we consider a good quantity $D_t\mathfrak{K}_\cH-v\cdot (\na \mathfrak{K}_\cH|_{c})$ which satisfies
\[
\left\{\begin{array}{ll}
\D \big(D_t\mathfrak{K}_\cH-v\cdot (\na\mathfrak{K}_\cH|_{c})\big)=2tr(\na v \na^2 \mathfrak{K}_\cH)-\D\big(v\cdot (\na\mathfrak{K}_\cH|_{c})\big) \qquad\hbox{in}\quad\Om_t\\
\na_{n_t}(D_t\mathfrak{K}_\cH-v\cdot (\na\mathfrak{K}_\cH|_{c}))\big|_{\G_t}=D_t \mJ^\perp +\na_{n_t} v \cdot \big(\na \mathfrak{K}_\cH-\na\mathfrak{K}_\cH|_{c}\big)-v\cdot \na_{n_t}(\na\mathfrak{K}_\cH|_{c})\big|_{\G_t},\\
 \na_{n_b}(D_t\mathfrak{K}_\cH-v\cdot (\na\mathfrak{K}_\cH|_{c}))\big|_{\G_b}=\na_{n_b}v\cdot \big(\na\mK_\cH-\na\mathfrak{K}_\cH|_{c}\big)
 -v\cdot \na_{n_b}(\na\mathfrak{K}_\cH|_{c})\big|_{\G_b}.
\end{array}\right.
\]
As a result, similar arguments as in \eqref{pa2 v estimate} applied for the part $\na v\cdot\big(\na\mK_\cH-\na\mathfrak{K}_\cH|_{c}\big)$ , we can  have
\beq\label{est:D_tk-v}
\| D_t\mathfrak{K}_\cH-v\cdot (\na\mathfrak{K}_\cH|_{c})\|_{H^{5/2}(\Om_t)}  \leq  P(E(t))(1+E_1(t)^{1/2}).
\eeq
where $\|D_t\mJ^\perp\|_{L^2(\G_t)}$ is handled by \eqref{est:D_tk-1}.

\medskip

On the other hand, we deal with $D_t\mJ$.  First, we know directly from \eqref{est:D_tk-1} that 
\beq\label{DtJ H1}
\|D_t\mJ\|_{H^1(\Om_t)}\le  P(E(t))(1+ \|D_t  \mathfrak{J}^\bot\|_{H^{1/2}(\G_t)}).
\eeq
Second, we notice that
\[
 D_t\mJ=D_t\na \mathfrak{K}_\cH=\na \big(D_t\mathfrak{K}_\cH-v\cdot (\mJ|_{c})\big)-\na v\cdot \big(\mJ-\mJ|_{c}\big).
\]
Similarly as in \eqref{apply hardy},  we apply Lemma \ref{product}(1) and \eqref{na v decomp} to obtain the estimate
\[
\begin{split}
&\|\na v\cdot (\mJ-\mJ|_{c})\|_{H^{3/2}(\Om_t)}\\
&\le P(E(t))\big(1+\|(\pa\na v)\cdot (\mJ-\mJ|_{c})\|_{H^{1/2}(\Om_t)}\big)\\
&\le P(E(t)\Big(1+\sum_{i=l,r}\big(\|r^{\pi/\om_i-2}\|_{L^\infty(0,r_0)}+\|r^{\pi/\om_i-2}\|_{H^1(0,r_0)}\big)
\|r^{-1}(\mJ-\mJ|_{p_i})\circ T^{-1}_i\|_{H^{1/2}(\cS_{t,i})}\Big)\\
&\le P(E(t)),
\end{split}
\]
while recall that $\cS_{t,i}$ are straightened sector of $\Om_t$ with radius $r_0$ near the corners.

Therefore, combining this estimate with \eqref{est:D_tk-v}, we have 
\beq\label{est:D_t na k}
\|D_t\mJ\|_{H^{3/2}(\Om_t)} \leq P(E(t))(1+E_1(t)^{1/2})
\eeq
and also
\[
\|D_t\mK_\cH\|_{H^2(\Om_t)}\le P(E(t))(1+E_1(t)^{1/2}).
\]

\noindent - Estimates for $D_tP_{v,v}$. Checking system \eqref{DtPvv system} carefully and applying lemmas from Section \ref{estimate section}, one has immediately 
\beq\label{higher estimate DtPvv}
\|D_tP_{v,v}\|_{H^{5/2}(\Om_t)}\le P(E(t)),
\eeq
where we use arguments similar to \eqref{boundary term esti} on boundary terms.

\medskip
\noindent - Estimates for $v$. 
To begin with,  applying Remark \ref{interp Neumann}, we obtain for a small number $\eps\in (0, \pi/\om_i-2)$ (when $\om_i\in (\pi/3, \pi/2)$) that
\beq\label{higher estimate v perp}
\begin{split} 
\|v\|_{H^{2+\eps}(\Om_t)}\le P(E(t)).
\end{split}
\eeq

Moreover, we improve the regularity of $v^\perp$. In fact, similarly as in \eqref{est: v-H-1}, we have thanks to \eqref{est:D_tk-1},\eqref{higher estimate DtPvv} and \eqref{higher estimate v perp} the estimate
\beq\label{higher estimate v perp}
\begin{split}
\|\D_{\G_t}v^\perp\|_{H^{1/2+\eps}(\G_t)} &\leq  C(\|\G_t\|_{H^{4}})\big(\|D_t\kappa\|_{H^{1/2+\eps}(\G_t)}  + \|v\|_{H^{1/2+\eps}(\Om_t)}\big)\\
&\leq  P(E(t))(1+E_1(t)^{1/2}).
\end{split}
\eeq

\noindent - Estimates for $D^2_tv, D^3_tv$.  First, acting $D_t$ on both sides of Euler's equation leads to
\[
D_t^2 v =-D_t \mathfrak{J}-D_t\na P_{v,v}.
\]
Applying Lemma \ref{prop: DtPvv}, Lemma \ref{est: P_{v, v}-H}, \eqref{DtJ H1}, and \eqref{est:D_t na k}, we find 
\beq\label{est:D_t v-2-H}
\|D_t^2v\|_{H^{3/2}(\Om_t)}
\leq  \|D_t \mathfrak{J}\|_{H^{3/2}(\Om_t)}+\|D_t\na P_{v,v}\|_{H^{3/2}(\Om_t)} \leq   P(E(t))(1+E_1(t)^{1/2}).
\eeq

\medskip

Next, we give the estimates of $D_t^3 v$.  Before that, we firstly prove the following lemma.
\begin{lemma}\label{lem:cp}
Let $P_{\mathfrak{J}, v}, P_{\na P_{\mathfrak{J}, v}, v}$ and $P_{\cD_t \mathfrak{J}, v}$ be defined by \eqref{eq: cP_{w,v}}. Then there hold
\[
\|D_t\na P_{\mathfrak{J}, v}+\na P_{\na P_{\mathfrak{J}, v}, v}+\na P_{\cD_t \mathfrak{J}, v}\|_{H^1(\Om_t)} \leq   P(E(t))\big(1+ E_1(t)^{1/2}\big)
\]
and
\beno
\|D^2_t\na P_{v,v}\|_{L^2(\Om_t)} \leq P(E(t))\big(1+ E_1(t)^{1/2}\big).
\eeno

\end{lemma}
\no {\bf Proof.}
First, we deal with $D_t P_{\mathfrak{J}, v}$. Recalling system \eqref{eq:D_t P}, we have by Proposition \ref{elliptic estimate Neumann H3},  \eqref{pt CJv}, \eqref{est:D_tk-1}, \eqref{est:D_t na k} and lemmas from Section \ref{estimate section} that
\[
\begin{split}
&\|D_t P_{\mathfrak{J}, v}\|_{H^2(\Om_t)}\\
&\leq C(\|\G_t\|_{H^{5/2}})\Big(\|- D_ttr(\na\mathfrak{J}\na v)+2tr(\na v \na^2 \cP_{ \mathfrak{J},v})\|_{L^2(\Om_t)} + \|D_t( \mJ\cdot\na_v n_b)+\na_{n_b}v\cdot \na P_{ \mathfrak{J},v}\|_{H^{1/2}(\G_b)}\\
&\quad+\|C'_{ \mathfrak{J}, v}-D_t(\na_{\tau_t}\mathfrak{K}_{\cH}\cdot  \na_{\tau_t} v \cdot n_t)+\na_{n_t}v\cdot \na P_{ \mathfrak{J},v}\|_{H^{1/2}(\G_t)} +\|D_t P_{\mathfrak{J}, v}\|_{L^2(\Om_t)}\Big)\\
&\leq P(E(t))(1+E_1(t)^{1/2}).
\end{split}
\]

Second, checking the definition \eqref{eq: cP_{w,v}} and applying Lemma \ref{product} (1) and Lemma \ref{est: P_{v, v}-H}, we can have  for $ P_{\na P_{\mathfrak{J}, v}, v}$ the following estimate immediately:
\beq\label{P PJv,v estimate}
\|P_{\na P_{\mathfrak{J}, v}, v}\|_{H^{5/2}(\Om_t)} \le  P(E(t)).
\eeq
Moreover,  a similar argument also leads to
\[
\|P_{\cD_t \mathfrak{J}, v}\|_{H^2(\Om_t)} \leq  P(E(t))\big(1+ E_1(t)^{1/2}\big).
\]

In the end, for  the part $D_t^2 P_{v,v}$, we deduce from \eqref{DtPvv system} the following system:
 \beq\label{eq: P_{vv}-v-2}
\left\{\begin{array}{ll}
\D (D^2_t P_{v, v} )=-tr D^2_t (\na v \na v) +2trD_t  (\na  v  \na^2  P_{v, v})\qquad \hbox{in}\quad \Om_t\\
\na_{n_t}( D^2_t P_{v, v} )|_{\G_t}=C''_{v,v}(t)+D_t(\na_{n_t} v\cdot \na P_{v, v} )+\na_{n_t} v\cdot \na D_t P_{v, v}\big|_{\G_t},\\
 \na_{n_b} ( D^2_t P_{v, v}) |_{\G_b}=D^2_t(v\cdot \na_v n_b)+D_t(\na_{n_b} v\cdot  \na P_{v, v} )+\na_{n_b} v\cdot  \na D_t P_{v, v} \big|_{\G_b}
\end{array}
\right.
\eeq
with 
\[
\int_{\Om_t}D^2_tP_{v,v}dX=0.
\]
Applying \eqref{est:D_t v-2}, \eqref{est:D_t v-2-H} and checking term by term, we obtain the variational estimate
\beq\label{est:D_tP-2}
\|D^2_t P_{v, v} \|_{H^1(\Om_t)} \leq  P(E(t))(1+E_1(t)^{1/2}),
\eeq
and the proof is finished.

\ef

Now, we are in a position to give the estimate for $D_t^3 v$. In fact, acting $D_t^2$ on both sides of Euler's equation, we derive by the previous lemma and lemmas from Section \ref{estimate section}  that
\beq\label{est:D_t v-2}
\|D_t^3v\|_{L^2(\Om_t)} \leq \|D^2_t \mathfrak{J}\|_{L^2(\Om_t)} +\|D^2_t\na P_{v,v}\|_{L^2(\Om_t)}  \le P(E(t))(1+E_1(t)^{1/2}).
\eeq

\medskip
\noindent - Estimate for $v^\perp$.
We obtain by the definition of $\mathfrak{J}$ that
\[
\na (D_t^2 \mathfrak{K}_\cH) =D_t^2 \mathfrak{J}+ \na v\cdot \na D_t \mathfrak{K}_\cH+ D_t(\na v\cdot \na  \mathfrak{K}_\cH),
\]
and applying \eqref{Dt k 1} and \eqref{est:D_tk-1} implies 
\beq\label{D2t kH estimate}
\|D_t^2 \mathfrak{K}_\cH\|_{H^1(\Om_t)} \leq P(E(t))(1+E_1(t)^{1/2}).
\eeq
Recalling \eqref{Dt k 1} again, we know 
\[
\si D_t\kappa=D_t (\mathfrak{K}+P_{v,v})=-\D_{\G_t}v^\perp-v^\perp|\na_{\tau_t}n_t|^2+\na_{\tau_t}\na_{v^\top}n_t\cdot \tau_t \qquad\hbox{on}\quad \G_t,
\]
so we obtain by \eqref{est:D_tP-2} and \eqref{D2t kH estimate} the estimate
\beq\label{Dt vperp estimate}
\| D_t \D_{\Gamma_t}v^\perp\|_{L^2(\Gamma_t)}   \leq P(E(t))(1+E_1(t)^{1/2}).
\eeq


%
%
%
%
%
\medskip
\noindent - Some more higher-order estimates for $\mK_\cH$, $v$ and $\mJ$.

\begin{lemma}\label{lem:J-v-h} Assuming that $E(t), E_1(t)\in L^\infty[0,T]$ for some $T>0$, we have the following estimate:
\[
\|\mathfrak{K}_\cH\|_{H^3(\Om_t)}+\|D_t v\|_{H^2(\Om_t)} +\|\na^2 \mathfrak{K}_\cH\|_{L^\infty(\Om_t)}  +\|D_t \na v\|_{L^\infty(\Om_t)}\leq P(E(t))(1+E_1(t)^{1/2}).
\]
\end{lemma}
\no{\bf Proof.} (1) Higher regularity for $\mK_\cH$.  To begin with, we use Euler's equation to rewrite \eqref{Dt2 kappa} as follows:
\[
D^2_t\mK= \big(\D_{\G_t}J^\perp-h_v\big)+[n_t,\D_{\G_t}](J+\na P_{v,v}+{\bf g})+2\s \Pi(\tau_t)\cdot \na_{\tau_t} J+(R_1+h_v)- D_t^2 P_{v,v}\quad\hbox{on}\  \G_t
\]
where we notice that $\D_{\G_t}(\na_{n_t} P_{v,v})|_{\G_t}=0$ thanks to the definition of $P_{v,v}$ and $(R_1+h_v)$ contains only $\pa v$ terms instead of $\pa^2 v$ terms.

Meanwhile, we know from \eqref{est:D_tP-2} and \eqref{D2t kH estimate} that
\[
\|D^2_t\mK\|_{H^{1/2}(\G_t)}+\| D_t^2 P_{v,v}\|_{H^{1/2}(\G_t)}\le P(E(t))(1+E_1(t)^{1/2}),
\]
so checking term by term in the equation above, we have immediately 
\beq\label{D J perp-hv estimate}
\|\D_{\G_t}J^\perp-h_v\|_{H^{1/2}(\G_t)}\le P(E(t))(1+E_1(t)^{1/2}).
\eeq

On the other hand, we know from \eqref{na v decomp} that (when $\om_i\in (\pi/3, \pi/2)$)
\[
h_v=h_{v,r}+h_{v,s}\quad\hbox{with}\quad h_{v,r}\in H^{1/2}(\G_t),\ h_s=\sum_i\big[\chi_{\om} a_{i,1}r^{2\alpha_i-3}+\chi_{\om}a_{i,2}(\pa v_r\circ T^{-1}_i)r^{\alpha_i-2}\big]\circ T_i,
\]
where we note $\al_i=\pi/\om_i-1\in (1,2)$, 
$a_{i,k}$ ($i=l, r$) contain $n_t, \tau_t, \pa n_t, \pa\tau_t$ and the singular coefficient from  \eqref{na v decomp}. Moreover, $a_{i,3}(\pa v_r\circ T^{-1}_i)$ is linear with respect to $\pa v_r\circ T^{-1}_i$, where we recall that $v_r=\na\phi_r\in H^3(\Om_t)$ and $\|v_r\|_{H^3(\Om_t)}$ is controlled by $P(E(t))$.  

As a result, we find 
\[
2\alpha_i-3\in (-1,1),\quad \alpha_i-2\in (-1,0),
\]
which implies immediately that 
\[
\|h_v\|_{L^p(\G_t)}\le P(E(t)),\qquad\hbox{for}\quad 1<p<\min\{(2-\alpha_i)^{-1}, |3-2\alpha_i|^{-1}\}.
\]
Summing up  the estimates above for $\|\D_{\G_t}J^\perp-h_v\|_{H^{1/2}(\G_t)}$ and $\|h_v\|_{L^p(\G_t)}$, we obtain
\[
\begin{split}
\|\D_{\G_t}J^\perp\|_{L^p(\G_t)}&\le \|\D_{\G_t}J^\perp-h_v\|_{L^p(\G_t)}+\|h_v\|_{L^p(\G_t)}\\
&\le C(\|\G_t\|_{H^{5/2}})\|\D_{\G_t}J^\perp-h_v\|_{H^{1/2}(\G_t)}+\|h_v\|_{L^p(\G_t)}\\
&\le P(E(t))(1+E_1(t)^{1/2}),
\end{split}
\]
and this leads to 
\[
\|\na_{\tau_t}J^\perp\|_{W^{1,p}(\G_t)}\le P(E(t))(1+E_1(t)^{1/2}).
\]

Applying Lemma \ref{embedding}, we finally show that
\beq\label{J perp higher order esti}
\|\na_{\tau_t}J^\perp\|_{L^\infty(\G_t)}+\|J^\perp\|_{H^{3/2+\epsilon}(\G_t)}\le P(E(t))(1+E_1(t)^{1/2})
\eeq
with $\epsilon=1-1/p$.

Consequently, applying Proposition \ref{elliptic estimate Neumann H3}, we have the desired estimate for $\|\mK_\cH\|_{H^{3}(\Om_t)}$. Notice that we have in fact the estimate by Remark \ref{interp Neumann}:
\beq\label{kH higher order estimate}
\|\mK_\cH\|_{H^{3+\epsilon}(\Om_t)}\le P(E(t))(1+E_1(t)^{1/2}).
\eeq
Moreover, we apply Lemma \ref{trace} to have  $\mK=\si\kappa-P_{v,v}\in H^{5/2+\epsilon}(\G_t)$ with the estimate
\[
\|\mK\|_{H^{5/2+\epsilon}(\G_t)}\le P(E(t))(1+E_1(t)^{1/2}).
\] 

(2) $H^2$ estimate for $D_tv$. Recalling the Euler's equations and $P_{v,v}$ estimate from Lemma \ref{est: P_{v, v}-H}, 
we get
\beno
\|D_t v\|_{H^2(\Om_t)}\le \|\mJ\|_{H^2(\Om_t)}+\|\na P_{v,v}\|_{H^2(\Om_t)}+\|{\bf g}\|_{H^2(\Om_t)} \leq P(E(t))(1+E_1(t)^{1/2}).
\eeno

\medskip

(3) $L^\infty$ estimates. Applying Remark \ref{interp Neumann} and \eqref{kH higher order estimate} with the same $\epsilon=1-1/p$ as above lead immediately to 
\[
\begin{split}
\|\na^2  \mathfrak{K}_\cH\|_{L^\infty} &\le C(\|\G_t\|_{H^{5/2}})\|\na^2\mK_\cH\|_{H^{1+\epsilon}(\Om_t)}
\le P(E(t))(1+E_1(t)^{1/2}).
\end{split}
\]
Applying Remark \ref{interp Neumann} to $P_{v,v}$, we obtain 
\[
\begin{split}
&\|\na^2 P_{v,v}\|_{L^\infty(\Om_t)}\le C(\|\G_t\|_{H^{5/2}})\|\na^2P_{v,v}\|_{H^{1+\eps}(\Om_t)}\\
&\le C(\|\G_t\|_{H^4})\big(\|tr(\na v\na v)\|_{H^{1+\epsilon}(\Om_t)}+\|C_{v,v}(t)\|_{H^{3/2+\eps}(\G_t)}+\|v\cdot\na_vn_b\|_{H^{1/2+\epsilon}(\G_b)} +\|P_{v,v}\|_{L^2(\Om_t)}\big)\\
&\leq P(E(t))\big(1+\|\na v\|^2_{H^{1+\epsilon}(\Om_t)}\big).
\end{split}
\]
Using \eqref{na v decomp} and a similar argument as in (1), we have
\[
\|\na v\|_{H^{1+\epsilon}(\Om_t)}\le \|\na^2\phi_s\|_{H^{1+\epsilon}(\Om_t)}+\|\na^2\phi_r\|_{H^2(\Om_t)}
\le P(E(t)),
\]
which implies 
\[
\|\na^2 P_{v,v}\|_{L^\infty(\Om_t)}\le P(E(t)).
\]

In the end,  apply Euler's equation again leads to 
\[
 \|\na D_t v\|_{L^\infty(\Om_t)}\le \|\na^2\mK_\cH\|_{L^\infty(\Om_t)}+\|\na^2 P_{v,v}\|_{L^\infty(\Om_t)}\leq P(E(t))(1+E_1(t)^{1/2}),
\]
which finishes the proof.

\ef

\noindent - Estimate for $D^2_tP_{v,v}$.
Based the above estimates, we firstly improve the estimate for $D_tP_{v,v}$. Using system \eqref{eq: P_{vv}-v} and Lemma \ref{lem:J-v-h},
we improve the estimate in Lemma \ref{prop: D_t P_{v, v}-H} and get
\beq\label{est:D_t P_{v,v}-H}
\begin{split}
&\| D_t P_{v, v}-v\cdot (\na P_{v,v}|_{c}) \|_{H^3(\Om_t)} \\
&\leq C(\|\G_t\|_{H^3})\Big(\|-tr D_t (\na v \na v) +[\D,v]\cdot (\na P_{v, v}- \na P_{v,v}|_{c})+v\cdot\D(\na P_{v,v}|_c) \|_{H^1(\Om_t)}+| C'_{v,v}(t)|\\
&\quad+\|\na_{n_t} v\cdot (\na P_{v, v}-\na P_{v,v}|_{c})-v\cdot\na_{n_t}(\na P_{v,v}|_c)\|_{H^{3/2}(\G_t)}
+\|D_t(v\cdot\na_v n_b)\|_{H^{3/2}(\G_b)}\\
&\quad +\|\na_{n_b} v\cdot (\na P_{v, v}-\na P_{v,v}|_{c})-v\cdot\na_{n_b}(\na P_{v,v}|_c)\|_{H^{3/2}(\G_b)}\Big)\\
 &\leq P(E(t))(1+E_1(t)^{1/2}).
\end{split}
\eeq
Here we use Lemma \ref{product} (2) for the boundary terms like $\|\na_{n_t} v\cdot (\na P_{v, v}-\na P_{v,v}|_{c})\|_{H^{3/2}(\G_t)}$, which are handled similarly as in the proof for \eqref{est:D_t na k}. 
 \medskip

Next, we derive the equation of $D_t(D_t P_{v, v}-v\cdot (\na P_{v,v}|_{c}))$. To simplify the notation, we define
\[
\mathcal P_{t,1}=D_t P_{v, v}-v\cdot (\na P_{v,v}|_{c}),
\]
and we rewrite \eqref{est:D_t P_{v,v}-H} as 
\beq\label{tilde P estimate}
\|\mathcal P_{t,1} \|_{H^3(\Om_t)}\le P(E(t))(1+E_1(t)^{1/2}).
\eeq
Direct computations lead to the following system for $D_t\mathcal P_{t,1}$:
 \beq\label{eq: P_{vv}-v-2}
\left\{\begin{array}{ll}
\D (D_t\mathcal P_{t,1})=-tr D^2_t (\na v \na v) +2trD_t\big(\na v\na (\na P_{v, v}- \na P_{v,v}|_{c}) \big)-D_t\big(v\cdot \D(\na P_{v,v}|_c)\big)\\
\qquad\qquad\qquad+[D_t, \D]\mathcal P_{t,1}\qquad \hbox{in}\quad\Om_t\\
\na_{n_t}( D_t \mathcal P_{t,1})\big|_{\G_t}=C''_{v,v}(t)+D_t(\na_{n_t} v\cdot \big(\na P_{v, v}-\na P_{v,v}|_{c})\big)-D_t\big(v\cdot\na_{n_t}(\na P_{v,v}|_c)\big)+\na_{n_t} v\cdot \na \mathcal P_{t,1}\big|_{\G_t},\\
\na_{n_b}( D_t \mathcal P_{t,1})\big|_{\G_b}= D^2_t(v\cdot \na_v n_b)+D_t\big(\na_{n_b} v\cdot (\na P_{v, v}-\na P_{v,v}|_{c})\big)-D_t\big(v\cdot\na_{n_b}(\na P_{v,v}|_c)\big)\\
\qquad\qquad\qquad\qquad+\na_{n_b} v\cdot \na \mathcal P_{t,1}\big|_{\G_b}.\\
\end{array}
\right.
\eeq

Moreover, we define 
\beq\label{eq:cP_1}
\mathcal P_{t,2}=D_t \mathcal P_{t,1}-v\cdot (\na \mathcal P_{t,1}|_{c})= D^2_tP_{v,v}-D_t\big(v\cdot(\na P_{v,v}|_c)\big)-v\cdot (\na \mathcal P_{t,1}|_{c})
\eeq
and we modify this system above as in Lemma \ref{prop: D_t P_{v, v}-H} into a new system for $\mathcal P_{t,2}$ below:
  \beq\label{eq: P_{vv}-v-3}
\left\{\begin{array}{ll}
\D \mathcal P_{t,2}=-tr D^2_t (\na v \na v) +2trD_t\big(\na v\na (\na P_{v, v}- \na P_{v,v}|_{c}) \big)-D_t\big(v\cdot \D(\na P_{v,v}|_c)\big)\\
\qquad\qquad+[D_t, \D]\mathcal P_{t,1}-v\cdot \D(\na \mathcal P_{t,1}|_c)\qquad \hbox{in}\quad\Om_t\\
\na_{n_t}\mathcal P_{t,2}\big|_{\G_t}=C''_{v,v}(t)+D_t(\na_{n_t} v\cdot \big(\na P_{v, v}-\na P_{v,v}|_{c})\big)-D_t\big(v\cdot\na_{n_t}(\na P_{v,v}|_c)\big)\\
\qquad\qquad\qquad
+\na_{n_t} v\cdot (\na \mathcal P_{t,1}-\na \mathcal P_{t,1}|_c)-v\cdot\na_{n_t}(\na\mathcal P_{t,1}|_c)\big|_{\G_t},\\
\na_{n_b}\mathcal P_{t,2}\big|_{\G_b}= D^2_t(v\cdot \na_v n_b)+D_t\big(\na_{n_b} v\cdot (\na P_{v, v}-\na P_{v,v}|_{c})\big)-D_t\big(v\cdot\na_{n_b}(\na P_{v,v}|_c)\big)\\
\qquad\qquad\qquad
+\na_{n_b} v\cdot (\na \mathcal P_{t,1}-\na \mathcal P_{t,1}|_c)-v\cdot\na_{n_b}(\na\mathcal P_{t,1}|_c)\big|_{\G_b}.\\
\end{array}
\right.
\eeq
Thanks to Lemma \ref{prop: D_t P_{v, v}-H}, \eqref{est:D_t v-2-H} and \eqref{tilde P estimate},  it is straightforward to show that
\beq\label{D2t P estimate}
 \|\mathcal P_{t,2}\|_{H^{5/2}(\Om_t)}\leq P(E(t))(1+E_1(t)^{1/2}).
\eeq

\medskip

\noindent - The boundary condition for $D_t^2 \mathfrak{J}$ at corner points.
\begin{lemma}\label{lem:cc-1}
 We have at the contact points $p_i (i=l,r)$ the following equations
\beq\label{equ:cc-1}
 D^2_t  \mathfrak{J} =\f{\sigma^2 }{\beta_c}(\sin\om_i) D_t(\na_{\tau_t}  \mathfrak{J}^\bot)
\tau_b+R_{c, h, i},
\eeq
 where the remainder terms
\[
R_{c, h, i}=D_tR_{c1} + \f{\sigma^2 }{\beta_c}D_t(\tau_b\sin\om_i )  (\na_{\tau_t}  \mathfrak{J})^\bot
- \f{\sigma^2 }{\beta_c}(\sin\om_i )  D_t(\na_{\tau_t} n_t\cdot \mathfrak{J})\tau_b\big|_{p_i}
\]
with $R_{c1}$ defined in Lemma \ref{lem:cc}. 
Moreover, there hold for $i=l, r$ that
\[
| R_{c, h,i}|\leq P(E(t))\big(1+E_1(t)^{1/2}\big).
\]
 \end{lemma}
\no{\bf Proof.} First, Acting $D_t$ directly on both sides of the equations in Lemma  \ref{lem:cc}, one has the desired equations immediately. Next, apply lemmas from Section \ref{estimate section},  Lemma \ref{lem:J-v-h} \eqref{kH higher order estimate} and \eqref{D2t P estimate}, the estimates follow. 
Moreover,  we notice that the higher-order term $\na^2_{\tau_t}v^\perp(\tau_t\cdot\mJ)$ from $ D_t(\na_{\tau_t} n_t\cdot \mathfrak{J})\tau_b$ in $R_{c,h,i}$ can be  controlled thanks to \eqref{higher estimate v perp}, where we use 
 \beq\label{Dtnt expression}
 D_t n_t=-\big((\na v)^*n_t\big)^ \top=-\big(n_t\cdot\na_{\tau_t}v\big)\tau_t=-(\na_{\tau_t}v^\perp)\tau_t+(\na_{\tau_t}n_t\cdot v)\tau_t.
 \eeq
\ef

\medskip

\subsection{Proof of Theorem \ref{1-order energy estimates}}
At this moment, we are finally ready to prove the higher-order energy estimate in Theorem \ref{1-order energy estimates}.  
To begin with, we rewrite equation \eqref{eq: D_t^2 J } of $\cD_t^2 \mathfrak{J}$ into:
\beq\label{eq: D_t^2 J-1}
D_t \big(\cD^2_t \mathfrak{J} -\na \cH(\cP_{t,2}) \big)+\na\cH\big(\D_{\Gamma_t}  D_t \mathfrak{J}^\bot-D_t h_v\big)= \cR_2 +D_t \cD_t\na P_{\mathfrak{J}, v}+D_t \na P_{\cD_t \mathfrak{J}, v},
\eeq
 where 
\[
\begin{split}
\cR_2&= D_t \big(\tilde{R}_0-\na\cH(D_t^2P_{v,v}) \big )+[\na, D_t]\cH (\D_{\Gamma_t} \mathfrak{J}^\bot-h_v)+\na [\cH, D_t](\D_{\Gamma_t} \mathfrak{J}^\bot-h_v)+\na \cH [\D_{\Gamma_t}, D_t]\mathfrak{J}^\bot    \\
&\quad+D_t( \na P_{\na P_{\mathfrak{J}, v}, v} )  -D_t \na \cH \big[D_t\big(v\cdot (\na P_{v,v}|_{c})\big)+ v\cdot (\na \mathcal P_{t,1}|_{c})\big].
\end{split}
\]
In fact,  this complicated form is used due to similar technical reasons as \eqref{eqn cDt2 J}, see Remark \ref{explain the eqn}, where  the terms $\na \cH(\cP_{t,2})$,  $D_t h_v$ are added on the left side to improve the estimates. 
\medskip

Next, we apply $L^2(\Om_t)$ inner product with $\cD^2_t \mathfrak{J} -\na \cH(\cP_{t,2})$ on both sides of \eqref{eq: D_t^2 J-1} to get
 \[
 \begin{split}
 &\f12\f d{dt} \|\cD^2_t \mathfrak{J} -\na \cH(\cP_{t,2})\|^2_{L^2(\Om_t)} -\int_{\Om_t} \na\cH\big(\D_{\Gamma_t} D_t \mathfrak{J}^\bot-D_t h_v\big)\cdot \big(\cD^2_t \mathfrak{J} -\na \cH(\cP_{t,2})\big) dX \\
 &=\int_{\Om_t}   \cR_2    \cdot \big(\cD^2_t \mathfrak{J} -\na \cH(\cP_{t,2})\big)dX +\int_{\Om_t}\big(  D_t \cD_t\na P_{\mJ, v}+D_t \na P_{\cD_t \mJ, v}   \big)\cdot \big(\cD^2_t \mathfrak{J} -\na \cH(\cP_{t,2})\big)dX.
 \end{split}
\]
Consequently, for any $t'\in [0,T]$ there holds
\beq\label{eq:1-order}
 \begin{split}
&\f12\big\|\cD^2_t \mathfrak{J}(t') -\na \cH(\cP_{t,2})(t')\big\|^2_{L^2(\Om_{t'})} -\int_0^{t'}\int_{\Om_t} \na\cH\big(\D_{\Gamma_t} D_t \mathfrak{J}^\bot-D_th_v\big)\cdot \big(\cD^2_t \mathfrak{J} -\na \cH(\cP_{t,2})\big) dXdt\\
&=\f12 \big\|\cD^2_t \mathfrak{J}(0) -\na \cH(\cP_{t,2})(0)\big\|^2_{L^2(\Om_0)} +\int_0^{t'}\int_{\Om_t}   \cR_2    \cdot \big(\cD^2_t \mathfrak{J} -\na \cH(\cP_{t,2})\big)dXdt\\
 &\quad+\int_0^{t'}\int_{\Om_t}\big(  D_t \cD_t\na P_{J, v}+D_t \na P_{\cD_t J, v}   \big)\cdot \big(\cD^2_t \mathfrak{J} -\na \cH(\cP_{t,2})\big)dXdt.
  \end{split}
 \eeq

 \subsubsection{Left side of \eqref{eq:1-order}.}
In this subsection, we  prove  estimates for the second term on the left side of \eqref{eq:1-order}. 
 \begin{proposition}\label{prop:est-Left} One has for some $\delta_0\in (0,1)$ the following estimate:
\[
\begin{split}
&-\int_0^T\int_{\Om_t} \na\cH\big(\D_{\Gamma_t} D_t \mathfrak{J}^\bot-D_t h_v\big)\cdot \big(\cD^2_t \mathfrak{J} -\na \cH(\cP_{t,2})\big)dXdt \\
&\ge \f12  \sup_{t\in [0, T]}\int_{\Gamma_t}|\na_{\tau_t}D_t \mathfrak{J}^\bot |^2ds + \f14\int_0^T F_1(t)dt -\int_0^TP(E(t))\big(1+E_1(t)^{3/2}\big)dt\\
&\quad -\sup_{t\in [0, T]}P(E(t))\big(1+E_1(t)^{\delta_0}\big)-P(E(0))\big(1+E_1(0)^{\delta_0}\big).
\end{split}
\]
\end{proposition}
 \no {\bf Proof.}
Applying Green's Formula, we have
 \beq\label{eq: d-1}
\begin{split}
&- \int_{\Om_t} \na\cH\big(\D_{\Gamma_t} D_t \mathfrak{J}^\bot - D_t h_v\big)\cdot \big(\cD^2_t \mathfrak{J} -\na \cH(\cP_{t,2})\big)dX\\
& =-\int_{\Gamma_t}  \big(\D_{\Gamma_t}\cD_t \mathfrak{J}^\bot-D_t h_v\big)  \big(\cD^2_t \mathfrak{J} -\na \cH(\cP_{t,2})\big)\cdot n_tds\\
&= -\int_{\Gamma_t}  \big(\D_{\Gamma_t} D_t \mathfrak{J}^\bot- D_t h_v\big)  \big(D_t^2 \mathfrak{J}+D_t\na P_{\mathfrak{J}, v}+\na P_{\na P_{\mathfrak{J}, v}, v}+\na P_{\cD_t \mathfrak{J}, v}\big)\cdot n_tds\\
&\quad+\int_{\Gamma_t}  \big(\D_{\Gamma_t} D_t \mathfrak{J}^\bot- D_t h_v\big) \na_{n_t} \cH(\cP_{t,2}) ds\\
&=-\int_{\Gamma_t}  \D_{\Gamma_t} D_t \mathfrak{J}^\bot\, D_t^2 \mathfrak{J}\cdot n_tds+\int_{\Gamma_t} D_t h_v \   D_t^2 \mathfrak{J}  \cdot n_tds +I_R,
\end{split}
 \eeq
 where 
\[
 \begin{split}
I_R&=-\int_{\Gamma_t}  \big(\D_{\Gamma_t} D_t \mathfrak{J}^\bot -D_t h_v\big) \big( D_t\na P_{\mathfrak{J}, v}+\na P_{\na P_{\mathfrak{J}, v}, v}+\na P_{\cD_t \mathfrak{J}, v}\big)\cdot n_tds\\
&\quad
+\int_{\Gamma_t}  \big(\D_{\Gamma_t} D_t \mathfrak{J}^\bot- D_t h_v\big)\cdot   \na_{n_t} \cH(\cP_{t,2}) ds.
\end{split}
\]

We deal with these integrals above one by one.  First, integrating by parts and taking one $D_t$ out of the first integral in \eqref{eq: d-1}, we  get
\beq\label{split into 3 lemma}
 \begin{split}
& -\int_{\Gamma_t}  \D_{\Gamma_t} D_t \mathfrak{J}^\bot\, D_t^2 \mathfrak{J}\cdot n_t ds\\
&=  \f12\f d{dt}\int_{\Gamma_t}|\na_{\tau_t}D_t \mathfrak{J}^\bot |^2ds-\f12 \int_{\Gamma_t}|\na_{\tau_t}D_t \mathfrak{J}^\bot |^2\,D_tds+ \na_{\tau_t} D_t \mathfrak{J}^\bot \,(D_t^2 \mathfrak{J}\cdot n_t)\big|^{p_r}_{p_l}\\
 &\quad+\int_{\Gamma_t}  \na_{\tau_t} D_t \mathfrak{J}^\bot\, [\na_{\tau_t},D_t] D_t\mathfrak{J}^\bot ds  +\int_{\Gamma_t}  \na_{\tau_t} D_t \mathfrak{J}^\bot \,\na_{\tau_t} \big([n_t, D_t^2]  \cdot\mathfrak{J}\big)ds.
\end{split}
\eeq
Here,   Lemma \ref{lem:cc-1} is applied at the corner terms:
\[
 \begin{split}
&\na_{\tau_t} D_t \mathfrak{J}^\bot \,(D_t^2 \mathfrak{J}\cdot n_t)\big|^{p_r}_{p_l}\\
&=
\na_{\tau_t} D_t \mathfrak{J}^\bot\ \f{\sigma^2 }{\beta_c}(\sin\om_l) D_t\na_{\tau_t}  \mathfrak{J}^\bot
(\tau_b\cdot n_t)\big|_{p_r}-\na_{\tau_t} D_t \mathfrak{J}^\bot\ \f{\sigma^2 }{\beta_c}(\sin\om_r) D_t\na_{\tau_t}  \mathfrak{J}^\bot(\tau_b\cdot n_t)\big|_{p_l}\\
&\quad+\na_{\tau_t} D_t \mathfrak{J}^\bot (R_{c,h,l}\cdot n_t)\big|_{p_r}-\na_{\tau_t} D_t \mathfrak{J}^\bot (R_{c,h,r}\cdot n_t)\big|_{p_l}\\
&=F_1(t)+\na_{\tau_t} D_t \mathfrak{J}^\bot\ \f{\sigma^2 }{\beta_c}(\sin\om_l) [D_t,\na_{\tau_t}]  \mathfrak{J}^\bot
(\tau_b\cdot n_t)\big|_{p_r}-\na_{\tau_t} D_t \mathfrak{J}^\bot\ \f{\sigma^2 }{\beta_c}(\sin\om_r) [D_t,\na_{\tau_t} ] \mathfrak{J}^\bot(\tau_b\cdot n_t)\big|_{p_l}\\
&\quad+\na_{\tau_t} D_t \mathfrak{J}^\bot (R_{c,h,l}\cdot n_t)\big|_{p_r}-\na_{\tau_t} D_t \mathfrak{J}^\bot (R_{c,h,r}\cdot n_t)\big|_{p_l},
\end{split}
\]
so one has by \eqref{kH higher order estimate} and lemmas from Section \ref{estimate section} that
\[
\na_{\tau_t} D_t \mathfrak{J}^\bot \,(D_t^2 \mathfrak{J}\cdot n_t)\big|^{p_r}_{p_l}\ge 
\f12 F_1(t)-P(E(t))\big(1+E_1(t)\big).
\]

On the other hand, since $D_tds=(v^\perp\kappa+\na_{\tau_t}(v\cdot\tau_t)\big)ds$, we obtain directly the estimate:
\[
\f12\int_{\Gamma_t}|\na_{\tau_t}D_t \mathfrak{J}^\bot |^2\,D_tds\le  P(E(t))E_1(t).
\]

Moreover, a direct computation using \eqref{Dtnt expression} shows that
\[
\begin{split}
[n_t, D_t^2]  \cdot\mathfrak{J}&=- D_t(D_tn_t)\cdot\mJ-2(D_tn_t)\cdot D_t\mJ\\
&=-D_t\big(-(\na_{\tau_t}v^\perp)\tau_t+(\na_{\tau_t}n_t\cdot v)\tau_t\big)\cdot\mJ-2\big(-(\na_{\tau_t}v^\perp)\tau_t+(\na_{\tau_t}n_t\cdot v)\tau_t\big)\cdot D_t\mJ.
\end{split}
\]
Therefore, thanks to \eqref{est:D_t na k}, Lemma \ref{lem:J-v-h} and checking term by term, we  have  the following estimate for the last two integrals in \eqref{split into 3 lemma}:
\[
\int_{\Gamma_t}  \na_{\tau_t} D_t \mathfrak{J}^\bot\, [\na_{\tau_t},D_t] D_t\mathfrak{J}^\bot ds  +\int_{\Gamma_t}  \na_{\tau_t} D_t \mathfrak{J}^\bot \,\na_{\tau_t} \big([n_t, D_t^2]  \cdot\mathfrak{J}\big)ds \leq P(E(t))\big(1+E_1(t)\big) .
\]

As a result,  the proof is finished as long as we have Lemma \ref{lem: f-D_tJ-1} and Lemma \ref{lem:I_R}.

 \ef

We deal with the remainder integrals in \eqref{eq: d-1} in the following two lemmas.
 \begin{lemma}\label{lem: f-D_tJ-1} One has for the second term of the right hand of \eqref{eq: d-1} that
\[
\begin{split}
\Big|\int_0^T\int_{\Gamma_t} D_t h_v \   D_t^2 \mathfrak{J}  \cdot n_tds \Big|  &\leq P(E(0))\big(1+E_1(0)^{\delta_0}\big)+\sup_{t\in[0, T]}P(E(t))\big(1+E_1(t)^{\delta_0}\big) \\
&\quad+\int_0^T P(E(t))\big(1+E_1(t)^{3/2}\big)dt+ \f{1}{8} \int_0^T F_1(t).
\end{split}
\]
\end{lemma}

 \no {\bf Proof.}
To begin with, one recalls the definition of $h_v$ in \eqref{hv def} to obtain
\[
\int_{\Gamma_t} D_t h_v \   D_t^2 \mathfrak{J}  \cdot n_tds=\int_{\Gamma_t} D_t\big(\D_{\G_t} v\cdot D_t n_t  \big) D_t^2 \mathfrak{J}  \cdot n_tds +\int_{\Gamma_t} D_t\big[2D^2v\big(\tau_t, (\na_{\tau_t}v)^\top\big)\cdot n_t\big] D_t^2 \mathfrak{J}  \cdot n_tds.
\]
We only deal with the first integral in the integral above, since the second one can be handled in a similar way.

In fact, it is straightforward to see by \eqref{Dtnt expression} that
\[
\begin{split}
&\int_{\Gamma_t} D_t\big(\D_{\G_t} v\cdot D_t n_t  \big) D_t^2 \mathfrak{J}  \cdot n_tds\\
&=\int_{\Gamma_t} D_t\big[\D_{\G_t} v\cdot \tau_t(-\na_{\tau_t}v^\perp+\na_{\tau_t}n_t\cdot v )\big] D_t^2 \mathfrak{J}  \cdot n_tds\\
&= \int_{\Gamma_t} D_t\D_{\G_t} v\cdot \tau_t(-\na_{\tau_t}v^\perp+\na_{\tau_t}n_t\cdot v) D_t^2 \mathfrak{J}  \cdot n_tds 
+ \int_{\Gamma_t}\D_{\G_t} v\cdot D_t\big[\tau_t(-\na_{\tau_t}v^\perp+\na_{\tau_t}n_t\cdot v)\big] D_t^2 \mathfrak{J}  \cdot n_tds \\
  &\triangleq I_1+I_2.
\end{split}
\]
The estimates for $I_1, I_2$ are proved  in the following lines.
\medskip

\noindent - Estimates of $I_1$.  First, direct computations lead to
\[
\begin{split}
  I_1 =&\int_{\Gamma_t} D_t\D_{\G_t} v\cdot \tau_t(-\na_{\tau_t}v^\perp+\na_{\tau_t}n_t\cdot v) D_t^2 (\mathfrak{J}-\mJ|_c)  \cdot n_tds
\\
  &  +\int_{\Gamma_t} D_t\D_{\G_t} v\cdot \tau_t(-\na_{\tau_t}v^\perp+\na_{\tau_t}n_t\cdot v) D_t^2 (\mathfrak{J}|_c)  \cdot n_tds  \\
= & \f d{dt} \int_{\Gamma_t} D_t\D_{\G_t} v\cdot \tau_t(-\na_{\tau_t}v^\perp+\na_{\tau_t}n_t\cdot v) D_t (\mathfrak{J}-\mJ|_c)  \cdot n_tds\\
&-\int_{\Gamma_t} D^2_t\D_{\G_t} v\cdot \tau_t(-\na_{\tau_t}v^\perp+\na_{\tau_t}n_t\cdot v) D_t(\mathfrak{J}-\mJ|_c)  \cdot n_tds\\
&-\int_{\Gamma_t} D_t\D_{\G_t} v\cdot D_t\big[\tau_t(-\na_{\tau_t}v^\perp+\na_{\tau_t}n_t\cdot v)\big] D_t (\mathfrak{J}-\mJ|_c)  \cdot n_tds\\
&-\int_{\Gamma_t} D_t\D_{\G_t} v\cdot \tau_t(-\na_{\tau_t}v^\perp+\na_{\tau_t}n_t\cdot v) \,D_t (\mathfrak{J}-\mJ|_c)  \cdot D_t(n_tds)\\
&  +\int_{\Gamma_t} D_t\D_{\G_t} v\cdot \tau_t(-\na_{\tau_t}v^\perp+\na_{\tau_t}n_t\cdot v) D_t^2 (\mathfrak{J}|_c)  \cdot n_tds  \\
\triangleq & \f d{dt} \int_{\Gamma_t} D_t\D_{\G_t} v\cdot \tau_t(-\na_{\tau_t}v^\perp+\na_{\tau_t}n_t\cdot v) D_t (\mathfrak{J}-\mJ|_c)  \cdot n_tds+\sum_{i=1}^4I_{1i}.
\end{split}
\]
We deal with $I_{1i}$ one by one.  In fact, applying integration by parts, we have  
\[
\begin{split}
I_{14}
&=D_t\na_{\tau_t} v\cdot \tau_t(-\na_{\tau_t}v^\perp+\na_{\tau_t}n_t\cdot v) \,D_t^2 (\mathfrak{J}|_c)  \cdot n_t\big|^{p_l}_{p_r}\\
&\quad-\int_{\Gamma_t} D_t\na_{\tau_t} v\cdot \na_{\tau_t}\big[ \tau_t(-\na_{\tau_t}v^\perp+\na_{\tau_t}n_t\cdot v) D_t^2 (\mathfrak{J}|_c)  \cdot n_t\big]ds\\
&\quad + \int_{\Gamma_t} [D_t,\na_{\tau_t}]\na_{\tau_t} v\cdot \tau_t(-\na_{\tau_t}v^\perp+\na_{\tau_t}n_t\cdot v) D_t^2 (\mathfrak{J}|_c)  \cdot n_tds\\
&\quad - \int_{\Gamma_t} D_t\big((\na_{\tau_t}\tau_t)^\top\cdot\na v\big)\cdot \tau_t(-\na_{\tau_t}v^\perp+\na_{\tau_t}n_t\cdot v) D_t^2 (\mathfrak{J}|_c)  \cdot n_tds.
\end{split}
\]
Using \eqref{higher estimate v perp} and similar arguments as in \eqref{pa2 v estimate} for $\pa^2v$ and $v^\perp$ terms and applying Lemma \ref{lem:cc-1} to $D^2\mJ|_{p_i}$, we derive 
\[
|I_{14}|\le \f1{16}F_1(t)+P(E(t))\big(1+E(t)\big).
\]

For $I_{11}$, since
\[
\begin{split}
&D^2_t\D_{\G_t} v=D^2_t\big(\na_{\tau_t}\na_{\tau_t}v-(\na_{\tau_t}\tau_t\cdot\tau_t)\na_{\tau_t} v\big)\\
&=\na_{\tau_t}D^2_t\na_{\tau_t}v+[D^2_t, \na_{\tau_t}]\na_{\tau_t}v-D^2_t \big((\na_{\tau_t}\tau_t\cdot\tau_t)\na_{\tau_t} v\big)\\
&=\na_{\tau_t}D^2_t\na_{\tau_t}v+D_t\big((D_t\tau_t-\na_{\tau_t}v)\cdot\na\na_{\tau_t}v\big)
+(D_t\tau_t-\na_{\tau_t}v)\cdot\na D_t\na_{\tau_t}v-D^2_t \big((\na_{\tau_t}\tau_t\cdot\tau_t)\na_{\tau_t} v\big)\\
&\triangleq h_{v,2}+l.o.t.
\end{split}
\]
where 
\[
h_{v,2}=\na_{\tau_t}D^2_t\na_{\tau_t}v+ 2\na_{\tau_t}\big((D_t\tau_t-\na_{\tau_t}v)\cdot\na D_t v\big)
-\na^2_{\tau_t}D_tv\cdot n_t(n_t\cdot\tau_t)\na_{\tau_t}v
-(\na_{\tau_t}\tau_t\cdot\tau_t)\na_{\tau_t}D^2_t v
\] contains higher-order terms of $v$, and 
the remainder part contains products like $\pa D_tv\,\pa^2v$ and other lower-order terms and can be controlled by \eqref{est:D_t na k}, \eqref{est:D_t v-2-H}, Lemma \ref{lem:J-v-h} and similar arguments as in \eqref{pa2 v estimate}.

Therefore, we arrive at 
\[
I_{11}=\int_{\Gamma_t} h_{v,2}\cdot \tau_t(-\na_{\tau_t}v^\perp+\na_{\tau_t}n_t\cdot v) D_t(\mathfrak{J}-\mJ|_c)  \cdot n_tds+l.o.t..
\]
To finish the estimate for $I_{11}$, we firstly take care of the first term $\na_{\tau_t}D^2_t\na_{\tau_t}v$ in $h_{v,2}$.  In fact,  we use integration by parts on $\G_t$ to obtain
\[
\begin{split}
&\int_{\Gamma_t} \na_{\tau_t}D^2_t\na_{\tau_t}v\cdot \tau_t(-\na_{\tau_t}v^\perp+\na_{\tau_t}n_t\cdot v) D_t(\mathfrak{J}-\mJ|_c)  \cdot n_tds\\
&=-\int_{\Gamma_t}D^2_t\na_{\tau_t}v\cdot  \na_{\tau_t}\big(\tau_t(-\na_{\tau_t}v^\perp+\na_{\tau_t}n_t\cdot v)\big) D_t(\mathfrak{J}-\mJ|_c)  \cdot n_tds\\
&\quad-\int_{\Gamma_t} D^2_t\na_{\tau_t}v\cdot \tau_t(-\na_{\tau_t}v^\perp+\na_{\tau_t}n_t\cdot v) \na_{\tau_t}\big(D_t(\mathfrak{J}-\mJ|_c)  \cdot n_t\big)ds.
\end{split}
\]
Similar arguments as in\eqref{pa2 v estimate} and  applying \eqref{est:D_t na k}, \eqref{est:D_t v-2-H}, we derive
\[
\int_{\Gamma_t} \na_{\tau_t}D^2_t\na_{\tau_t}v\cdot \tau_t(-\na_{\tau_t}v^\perp+\na_{\tau_t}n_t\cdot v) D_t(\mathfrak{J}-\mJ|_c)  \cdot n_tds\le P(E(t))\big(1+E_1(t)\big).
\]
Meanwhile, the other terms in $I_{11}$  can also be handled similarly, so we conclude
\[
|I_{11}|\le P(E(t))\big(1+E_1(t)\big).
\]

 \medskip

 For $I_{12}$, we have
 \[
 \begin{split}
I_{12}=&-\int_{\Gamma_t} \big([D_t,\na_{\tau_t}]\na_{\tau_t}v-D_t(\na_{\tau_t}\tau_t\cdot\tau_t\na_{\tau_t} v)\big)\cdot D_t\big[\tau_t(-\na_{\tau_t}v^\perp+\na_{\tau_t}n_t\cdot v)\big] D_t (\mathfrak{J}-\mJ|_c)  \cdot n_tds \\
&-\int_{\Gamma_t} \na_{\tau_t} D_t\na_{\tau_t}v\cdot D_t\big[\tau_t(-\na_{\tau_t}v^\perp+\na_{\tau_t}n_t\cdot v)\big] D_t (\mathfrak{J}-\mJ|_c)  \cdot n_tds,
\end{split}
\]
where the first integral can be handled as above. For the second integral, we know from Lemma \ref{H1 zero trace} that
\[
D_t (\mathfrak{J}-\mJ|_c)\in \tilde H^{1/2}(\G_t),
\]
so we apply Lemma \ref{H1 zero trace}, Lemma \ref{G_b H-1/2}, \eqref{est:D_t na k} and Lemma \ref{lem:J-v-h} to find
 \beq\label{tilde argument}
 \begin{split}
&-\int_{\Gamma_t} \na_{\tau_t} D_t\na_{\tau_t}v\cdot D_t\big[\tau_t(-\na_{\tau_t}v^\perp+\na_{\tau_t}n_t\cdot v)\big] D_t (\mathfrak{J}-\mJ|_c)  \cdot n_tds\\
&\le \big\|D_t\big[\tau_t(-\na_{\tau_t}v^\perp+\na_{\tau_t}n_t\cdot v)\big]\big\|_{L^\infty(\G_t)}\|\na_{\tau_t} D_t\na_{\tau_t}v\|_{\tilde H^{-1/2}(\G_t)}\|\|D_t (\mathfrak{J}-\mJ|_c)  \cdot n_t\|_{\tilde H^{1/2}(\G_t)}\\
&\le P(E(t))\|D_t\na_{\tau_t}v\|_{H^{1/2}(\G_t)}\|D_t\mJ\|_{H^1(\Om_t)}\\
&\le P(E(t))\big(1+E_1(t)\big).
\end{split}
\eeq
As a result, we summarize that
\[
|I_{12}|\le P(E(t))\big(1+E_1(t)^{3/2}\big).
\]

Moreover, similar arguments as for $I_{12}$, we also have
\[
 | I_{13}| \leq P(E(t))\big(1+E_1(t)\big).
\]
 
Together with all these estimates above for $I_{11}$ to $I_{14}$, we can go back to $I_1$ expression and integrate on both sides with respect to time $t$ on $[0,T]$ to find
\[
\begin{split}
\int^T_0 I_1 dt&\le \int_{\Gamma_t} D_t\D_{\G_t} v\cdot \tau_t(-\na_{\tau_t}v^\perp+\na_{\tau_t}n_t\cdot v) D_t (\mathfrak{J}-\mJ|_c)  \cdot n_tds\big|^T_0\\
&\quad+\int^T_0P(E(t))\big(1+E_1(t)^{3/2}\big)dt+\f{1}{16} \int_0^T F_1(t).
\end{split}\]

Moreover, we have for the first integral on the right side above the following estimate
\[
\begin{split}
&\int_{\Gamma_t} D_t\D_{\G_t} v\cdot \tau_t(-\na_{\tau_t}v^\perp+\na_{\tau_t}n_t\cdot v) D_t (\mathfrak{J}-\mJ|_c)  \cdot n_tds\\
&= \int_{\Gamma_t} \big([D_t,\na_{\tau_t}]\na_{\tau_t}v-D_t(\na_{\tau_t}\tau_t\cdot\tau_t\na_{\tau_t} v)\big)\cdot \tau_t(-\na_{\tau_t}v^\perp+\na_{\tau_t}n_t\cdot v) D_t (\mathfrak{J}-\mJ|_c)  \cdot n_tds\\
&\quad +\int_{\Gamma_t} \na_{\tau_t}D_t\na_{\tau_t} v\cdot \tau_t(-\na_{\tau_t}v^\perp+\na_{\tau_t}n_t\cdot v) D_t (\mathfrak{J}-\mJ|_c)  \cdot n_tds.
\end{split}
\]
Similar arguments as in \eqref{tilde argument} and checking carefully on the highest-order terms, we  derive
\[
\int_{\Gamma_t} D_t\D_{\G_t} v\cdot \tau_t(-\na_{\tau_t}v^\perp+\na_{\tau_t}n_t\cdot v) D_t (\mathfrak{J}-\mJ|_c)  \cdot n_tds
\le P(E(t))\big(1+E_1(t)^{\delta_0}\big)
\]
with the number $\delta_0\in (0,1)$ as above.

Consequently, we have the following estimate for $I_1$:
 \beq\label{est:I_1}
\begin{split}
 \Big|\int_0^T I_1 dt  \Big|  
 &\le P(E(0))\big(1+E_1(0)^{\delta_0}\big)+\sup_{t\in[0, T]}P(E(t))\big(1+E_1(t)^{\delta_0}\big) \\
 &\quad+\int_0^T P(E(t))\big(1+E_1(t)^{3/2}\big)dt+ \f{1}{16} \int_0^T F_1(t).
\end{split}
\eeq

\medskip

\noindent - Estimates of $I_2$. The integral $I_2$ can be  handled in the same way as $I_1$. We simply conclude that 
\beq\label{est:I_2}
\begin{split}
 \Big|\int_0^T I_2 dt  \Big|  
 &\le P(E(0))\big(1+E_1(0)^{\delta_0}\big)+\sup_{t\in[0, T]}P(E(t)\big(1+E_1(t)^{\delta_0}\big)\\
 &\quad +\int_0^T P(E(t))\big(1+E_1(t)^{3/2}\big)dt+ \f{1}{16} \int_0^T F_1(t).
\end{split}\eeq

In the end, combing \eqref{est:I_1} with \eqref{est:I_2},  the proof is finished.

 \ef

 \medskip

 Next, we deal with $I_R$.
 \begin{lemma}\label{lem:I_R} One has 
 \ben\label{est:I_R} 
|I_R| &\leq& \f1{8} F_1(t)+P(E_1(t)).
\een
 \end{lemma}
 \no{\bf Proof.} First, to simplify the notations in the first integral of $I_R$, we denote by 
 \[
 C_P =(D_t\na \cP_{\mathfrak{J}, v}+\na \cP_{\na \cP_{\mathfrak{J}, v}, v}+\na \cP_{\cD_t \mathfrak{J}, v})\cdot n_t.
 \]
Thanks to the definition of $P_{w, v}$ (see \eqref{eq: cP_{w,v}}), we have
\[
C_P=C'_{\mathfrak{J}, v}+\cC_{\na P_{\mathfrak{J}, v}, v}+C_{\cD_t \mathfrak{J}, v}   -D_t  n_t\cdot\na P_{\mathfrak{J}, v}- D_t \big((\mJ\cdot{\tau_t})\na_{\tau_t} v\cdot n_t\big)-(\na_{\tau_t}P_{\mathfrak{J}, v}+\cD_t\mathfrak{J}\cdot \tau_t) \na_{\tau_t} v\cdot n_t.\nonumber
\]
By \eqref{est:D_tk-1}, \eqref{est:D_t na k}, Lemma \ref{lem:J-v-h} and \eqref{D2t P estimate}, we  show immediately
\beq\label{C P estimate}
\|C_P\|_{H^{1}(\G_t)}+\|\na_{n_t}\cH(\cP_{t,2})\|_{H^{1}(\G_t)}\leq P(E(t))\big(1+E_1(t)^{1/2}\big).
\eeq
Now we are ready to have estimate of $I_R$. For the first integral in $I_R$, we  integrate by parts to derive
\[
\begin{split}
\int_{\Gamma_t} \D_{\Gamma_t} D_t \mathfrak{J}^\bot\, C_P ds
 &=\na_{\tau_t} D_t \mathfrak{J}^\bot  C_P\big|^{p_l}_{p_r}-\int_{\Gamma_t}  \na_{\tau_t} D_t \mathfrak{J}^\bot \cdot \na_{\tau_t}C_Pds \\
&\leq \f1{8} F_1(t)+P(E(t))\big(1+E_1(t)\big)
\end{split}
\]
where \eqref{est:D_t na k}, \eqref{higher estimate v perp}, \eqref{Dtnt expression}, Lemma \ref{lem:cc-1}  and \eqref{C P estimate} are used.

Moreover, we have from \eqref{hv def} and \eqref{C P estimate} that
\[
\begin{split}
\int_{\G_t} D_t h_v \, C_Pds
&=\int_{\G_t}\big[\D_{\G_t} v\cdot D_t n_t  +2D^2v\big(\tau_t, (\na_{\tau_t}v)^\top\big)\cdot n_t\big]C_Pds\\
&\leq P(E(t))\big(1+E_1(t)\big),
\end{split}
\]
where similar analysis as in the proof of the previous lemma is applied.

On the other hand,  the second integral in $I_R$
 can also be handled by a similar argument as above.
As a result, the proof is finished.

\ef

\medskip

\subsubsection{Right side of \eqref{eq:1-order}.} We firstly deal with the integral involving $\cR_2$.
  \begin{lemma}\label{lem: R^1} One has for some $\delta_0\in (0,1)$ the following estimate:
\[
\begin{split}
\Big|\int_0^T\int_{\Om_t}\cR_2 \cdot \big(\cD_t^2 \mathfrak{J}-\na\cH(\cP_{t,2})\big)dXdt
\Big|  &\le P(E(0))\big(1+E_1(0)^{\delta_0}\big)+\sup_{t\in[0, T]}P(E(t))\big(1+E_1(t)^{\delta_0}\big) \\
 &\quad+\int_0^T P(E(t))\big(1+E_1(t)^{3/2}\big)dt+ \f{1}{16} \int_0^T F_1(t).
\end{split}
\]
\end{lemma}
\no {\bf Proof.} Recalling from \eqref{eq: D_t^2 J-1}, we deal with the terms in $\cR_2$ one by one.
  \medskip

\noindent - Estimates of $D_t   \big(\tilde{R}_0-\na \cH(  D^2_t P_{v,v}  ) \big)  $. 
In fact, one knows directly from \eqref{rewrite eqn for Dt2 J} that
\[
\tilde{R}_0- \na \cH(  D^2_t P_{v,v}  ) =-\si\na \cH(J\cdot \D_{\G_t}n_t)+\si\na \cH\big([n_t,\D_{\G_t}]\cdot \na P_{v,v}\big)+\si\na\cH(R_1+h_v ) +A_1+A_2+A_3 ,
\]
where recall that $R_1+h_v$ only contains lower-order derivatives like $\pa v, \pa n_t, \kappa$ and $A_1$ to $A_3$ are defined in \eqref{A_1}-\eqref{A_3}.

As a result, checking term by term on $\tilde{R}_0- \na \cH(  D^2_t P_{v,v}  ) $, one finds that it contains $\pa D_t v$, $D_tJ$, $\pa v$ and other lower-order terms, so acting $D_t$ on $\tilde{R}_0- \na \cH(  D^2_t P_{v,v}  ) $ and applying lemmas in Section \ref{estimate section} and Section \ref{estimate section 2} lead to the following estimate:
\[
  \big\|D_t   \big(\tilde{R}_0-\na \cH(  D^2_t P_{v,v}  ) \big) \big\|_{L^2(\Om_t)}\leq P(E(t))\big(1+E_1(t)^{1/2}\big).
\]

\no - Estimates of $[\na, D_t]\cH (\D_{\Gamma_t}\mathfrak{J}^\bot-h_v) $.  It is straightforward to show by \eqref{D J perp-hv estimate} and Lemma \ref{Harmonic extension H1 estimate} that
\[
 \big\| [\na, D_t]\cH (\D_{\Gamma_t} \mathfrak{J}^\bot-h_v) \big\|_{L^2(\Om_t)} \leq P(E(t))\|\D_{\Gamma_t}\mathfrak{J}^\bot-h_v \|_{H^{1/2}(\G_t)}\leq P(E(t))\big(1+E_1(t)^{1/2}\big).
\]

\no - Estimates of $\na [\cH, D_t](\D_{\Gamma_t}\mathfrak{J}^\bot-h_v) $.  Recalling \eqref{commutator Dt H}
and using Lemma \ref{Harmonic extension H1 estimate}, \eqref{D J perp-hv estimate} and applying variational estimates as in page 33 \cite{MW2}  imply that
\[
 \big\|\na [\cH, D_t](\D_{\Gamma_t}\mathfrak{J}^\bot-h_v) \big\|_{L^2(\Om_t)} \leq P(E(t))\|\cH(\D_{\Gamma_t}\mathfrak{J}^\bot-h_v)\|_{H^1(\Om_t)}\leq P(E(t))\big(1+E_1(t)^{1/2}\big).
\]

\no - Estimates of $\na \cH [\D_{\Gamma_t}, D_t]\mathfrak{J}^\bot$.  To estimate this term, we use similar arguments as in the proof of Lemma \ref{lem: f-D_tJ-1}. 

First, using Green's Formula, we  obtain
\[
\begin{split}
 &\int_{\Om_t} \na \cH [\D_{\Gamma_t}, D_t]\mathfrak{J}^\bot \cdot \cD_t^2 \mathfrak{J}dX= \int_{\Gamma_t}[\D_{\Gamma_t}, D_t]\mathfrak{J}^\bot \, \cD_t^2 \mathfrak{J}\cdot n_tds\\
&=\int_{\G_t} [\big(\na_{\tau_t}\na_{\tau_t}-(\na_{\tau_t}\tau_t\cdot\tau_t)\na_{\tau_t}\big), D_t]\mJ^\perp \ \cD_t^2\mJ\cdot n_tds\\
 &= \int_{\Gamma_t}
2( \na_{\tau_t}v\cdot\tau_t)\na^2_{\tau_t}\mJ^\perp \ \cD_t^2\mJ\cdot n_tds\\
&\quad + \int_{\Gamma_t}\big[\na_{\tau_t}( \na_{\tau_t}v\cdot\tau_t)+D_t(\na_{\tau_t}\tau_t\cdot \tau_t)-( \na_{\tau_t}\tau_t\cdot\tau_t)( \na_{\tau_t}v\cdot\tau_t)
 \big] \na_{\tau_t}\mJ^\perp \  \cD_t^2 \mathfrak{J}\cdot n_tds.
\end{split}
\]
 These two integrals can be handled in a similar way as before, and here we only give the details for the first part in the second integral above. Moreover, the analysis for $\int_{\Om_t} \na \cH [\D_{\Gamma_t}, D_t]\mathfrak{J}^\bot \cdot  \na\cH(\cP_{t,2})dX$ can be done in a similar and easier way thanks to \eqref{D2t P estimate}, and  hence estimates for all these remainder parts are omitted.

Recalling \eqref{Dt2 J decomp} and the definition of $C_P$, we rewrite the first part in the second integral above  as follows:
\[
\begin{split}
 & \int_{\Gamma_t}\na_{\tau_t}( \na_{\tau_t}v\cdot\tau_t) \na_{\tau_t} \mathfrak{J}^\bot \, \cD_t^2 \mathfrak{J}\cdot n_t ds\\
  &=  \int_{\Gamma_t}\na_{\tau_t}( \na_{\tau_t}v\cdot\tau_t) \na_{\tau_t} \mathfrak{J}^\bot \, D_t^2 \mathfrak{J}\cdot n_t ds+  \int_{\Gamma_t}\na_{\tau_t}( \na_{\tau_t}v\cdot\tau_t) \na_{\tau_t} \mathfrak{J}^\bot \, C_P ds \triangleq I_3+I_4.
\end{split}
\]
Similarly as the estimates for $I_1$ in the proof of Lemma \ref{lem: f-D_tJ-1},  one has   for $I_3$ that
\[
\begin{split}
 I_3 &=\int_{\Gamma_t}\na_{\tau_t}( \na_{\tau_t}v\cdot\tau_t) \na_{\tau_t} \mathfrak{J}^\bot \, D_t^2 (\mathfrak{J}-\mJ|_c)\cdot n_t ds+\int_{\Gamma_t}\na_{\tau_t}( \na_{\tau_t}v\cdot\tau_t) \na_{\tau_t} \mathfrak{J}^\bot \, D_t^2 (\mathfrak{J}|_c)\cdot n_t ds\\
&\leq   \f1{8} F_1(t)+P(E(t))\big(1+E_1(t)\big)+\f d{dt}\int_{\Gamma_t}\na_{\tau_t}( \na_{\tau_t}v\cdot\tau_t) \na_{\tau_t} \mathfrak{J}^\bot \, D_t (\mathfrak{J}-\mJ|_c)\cdot n_t ds\\
&\quad -\int_{\Gamma_t}D_t\na_{\tau_t}( \na_{\tau_t}v\cdot\tau_t) \na_{\tau_t} \mathfrak{J}^\bot \, D_t (\mathfrak{J}-\mJ|_c)\cdot n_t ds
-\int_{\Gamma_t}\na_{\tau_t}( \na_{\tau_t}v\cdot\tau_t) D_t\na_{\tau_t} \mathfrak{J}^\bot \, D_t (\mathfrak{J}-\mJ|_c)\cdot n_t ds\\
&\le \f1{8} F_1(t)+P(E(t))\big(1+E_1(t)^{3/2}\big)+\f d{dt}\int_{\Gamma_t}\na_{\tau_t}( \na_{\tau_t}v\cdot\tau_t) \na_{\tau_t} \mathfrak{J}^\bot \, D_t (\mathfrak{J}-\mJ|_c)\cdot n_t ds,
\end{split}
\]
  where we use \eqref{J perp higher order esti}, Lemma \ref{lem:J-v-h} and  Lemma \ref{lem:cc-1} to have
\[
\begin{split}
&\int_{\Gamma_t}\na_{\tau_t}( \na_{\tau_t}v\cdot\tau_t) D_t\na_{\tau_t} \mathfrak{J}^\bot \, D_t (\mathfrak{J}-\mJ|_c)\cdot n_t ds\\
&\le P(E(t))\Big(1+\sum_i\big\|r^{\delta}(\na^2_{\tau_t} v)\circ T^{-1}_i\big\|_{L^\infty(0, r_0)}
\big\|\|D_t\na_{\tau_t}\mJ^\perp\|_{L^2(\G_t)}\big\|r^{-\delta}D_t (\mathfrak{J}-\mJ|_c)\circ T^{-1}_i\big\|_{L^2(0, r_0)}
\Big)
\end{split}
\]
with $1>\delta>3-\pi/\om_i$ for $\om_i\in (\pi/3,\pi/2)$  and 
\[
\begin{split}
\|D_t\na_{\tau_t}\mJ^\perp\|_{L^2(\G_t)}
&\le \big\|(\na_{\tau_t}v\cdot \tau_t)\na_{\tau_t}\mJ^\perp\big\|_{L^2(\G_t)}+\big\|\na_{\tau_t}D_t(\mJ\cdot n_t)\big\|_{L^2(\G_t)}\\
&\le P(E(t))\big(1+\|\na_{\tau_t}D_t\mJ\|_{L^2(\G_t)}+\|\na^2_{\tau_t}v^\perp\|_{L^2(\G_t)}\big)\le 
P(E(t))\big(1+E_1(t)\big).
\end{split}
\]

As a result, we have
\[
\begin{split}
 \Big|\int_0^T I_3 dt  \Big|  
 &\le P(E(0))\big(1+E_1(0)^{\delta_0}\big)+\sup_{t\in[0, T]}P(E(t))\big(1+E_1(t)^{\delta_0}\big) \\
 &\quad+\int_0^T P(E(t))\big(1+E_1(t)^{3/2}\big)dt+ \f{1}{64} \int_0^T F_1(t).
\end{split}
\]

Since similar arguments can be applied to $I_4$ and other integrals, we  conclude directly 
\[
\begin{split}
 \Big|\int_0^T\int_{\Om_t} \na \cH [\D_{\Gamma_t}, D_t]\mathfrak{J}^\bot \cdot \cD_t^2 \mathfrak{J}dXdt \Big|&\le P(E(0))\big(1+E_1(0)^{\delta_0}\big)+\sup_{t\in[0, T]}P(E(t))\big(1+E_1(t)^{\delta_0}\big) \\
 &\quad+\int_0^T P(E(t))\big(1+E_1(t)^{3/2}\big)dt+ \f{1}{64} \int_0^T F_1(t).
\end{split}
\]

\no - Estimates of $D_t( \na P_{\na P_{\mathfrak{J}, v}, v} )$.  Recalling the definition of $P_{\mathfrak{J}, v}$ and $P_{\na P_{\mathfrak{J}, v}, v} $ by \eqref{eq: cP_{w,v}}, one has firstly the estimate \eqref{P PJv,v estimate}.
Moreover, one can  take $D_t$ on the system of $P_{\na P_{\mathfrak{J}, v}, v}$ as in \eqref{DtPvv system} to obtain the system for $D_tP_{\na P_{\mathfrak{J}, v}, v}$. 

Consequently, checking term by term and applying Lemma \ref{est: P_{v, v}-H},  Lemma \ref{lem:cp}, one derives
\[
   \|D_t( \na \cP_{\na \cP_{\mathfrak{J}, v}, v} )\|_{L^2(\Om_t)} \leq P(E(t))\big(1+E_1(t)^{1/2}\big).
\]

\no - Estimates of the last term in $\cR_2$.  Thanks to \eqref{commutator Dt H}, \eqref{est:D_t v-2-H}, \eqref{D2t P estimate} and Lemma \ref{lem:cp}, we have 
\[
\big\|D_t \na \cH \big[D_t\big(v\cdot (\na P_{v,v}|_{c})\big)+ v\cdot (\na \mathcal P_{t,1}|_{c})\big]\big\|_{L^2(\Om_t)} \leq P(E(t))\big(1+E_1(t)^{1/2}\big).
\]

In the end, summing up all these estimates above, the proof is finished.

\ef

At this moment,   it remains to  handle the last integral on the right side of \eqref{eq:1-order}.   In fact, similar arguments as in  the proof of  Lemma \ref{lem: D_t P_{J, v}}, we can conclude the following estimate:
\[
\begin{split}
&\int_0^{t'}\int_{\Om_t}\big(  D_t \cD_t\na P_{J, v}+D_t \na P_{\cD_t J, v}   \big)\cdot \big(\cD^2_t \mathfrak{J} -\na \cH(\cP_{t,2})\big)dXdt\\
&\le P(E(0))\big(1+E_1(0)^{\delta_0}\big)+\sup_{t\in[0, T]}P(E(t))\big(1+E_1(t)^{\delta_0}\big) +\int_0^T P(E(t))\big(1+E_1(t)^{3/2}\big)dt+ \f{1}{16} \int_0^T F_1(t).
\end{split}
\]

 \subsubsection{The end of the higher-order energy estimate.}
Summing up all the estimates in the previous two subsections, we finally arrive at:
\[
\begin{split}
&
\sup_{t\in [0, T]}\Big(\big\|\cD^2_t \mathfrak{J} -\na \cH(\cP_{t,2})\big\|^2_{L^2(\Om_t)}+\int_{\Gamma_t}|\na_{\tau_t}D_t \mathfrak{J}^\bot |^2\Big)+\f1{4}\int_0^T F_1(t)\\
&\leq P(E(0))\big(1+E_1(0)^{\delta_0}\big)+\sup_{t\in[0, T]}P(E(t))\big(1+E_1(t)^{\delta_0}\big) +\int_0^T P(E(t))\big(1+E_1(t)^{3/2}\big)dt.
\end{split}
\]
Thanks to \eqref{D2t P estimate} and Lemma \ref{lem:cp}, we have for a  number $\delta_1>0$ small enough such  that
  \beno
\sup_{t\in [0, T]} E_1(t)+ \int_0^T F_1(t)\leq P(E(0))\big(1+E_1(0)^{\delta_0}\big)+\delta_1\sup_{t\in [0, T]} E_1(t)+\int_0^TP(E(t))\big(1+E_1(t)^{3/2}\big)dt.
 \eeno
Therefore, we've finished the proof of Theorem \ref{1-order energy estimates}.

\section{Well-posedness of  system $\mbox{(WW)}$}

In this section, we use  Picard iteration to prove the existence of solutions to $\mbox{(WW)}$. The main idea is the same as \cite{MW3}, although some necessary modifications are needed and presented here. As a result, we only show the skecth of the proof. For more details, see \cite{MW3}.

\subsection{Definitions of surfaces and domains}
        
Following \cite{SZ2, MW3}, we introduce a map $\Phi_{S_t}$ on the boundary $S_t=\G_t\cup \G_b$ to fix the moving domain $\Om_t$. To start with, we choose a reference domain $\Om_*$ with upper surface $\G_{t*}$ and bottom $\G_{b*}$, which can be taken as the initial domain $\Om_0$ without loss of generality. The contact points of $\Om_*$ are denoted by $p_{i*}$ ($i=l, r$) and the other notations follow similarly.

we define a unit upward vector field $\mu\in H^s(\G_{*}, \cS^1)$ with some large $s$ satisfying
\[
\mu\cdot n_{t*}\ge c_0\quad\hbox{on}\quad \G_{t*},\qquad\hbox{and}\quad\mu|_{p_{l*}}=-\tau_{b*}|_{p_{l*}},
\  \mu|_{p_{r*}}=\tau_{b*}|_{p_{r*}}
\]
for some fixed constant $c_0\in(0,1)$. Here we notice that the conditions above hold at $p_{l*}, p_{r*}$.

Applying the Implicit Function Theorem, there exists a small constant $d_0>0$  such that the map
 \[
 \Phi: \G_{t*}\times[-d_0,d_0]\rightarrow \R^2\quad\hbox{with}\quad
\Phi(p, d) \triangleq p+d\,\mu(p)
 \]
is an $H^s$ diffeomorphism  from its domain to a neighborhood of $\G_*$. 

Consequently, this map identifies each upper surface $\G_t$ near $\G_{t*}$ with a unique function 
\[
d_{\G_t}: \G_{t*}\rightarrow  \R
\]
 and we can define the following map
\[
\Phi_{S_t}: \ \G_{t*}\rightarrow \G_t\subset \R^2\quad\hbox{with}\quad \Phi_{S_t}(p)=p+d_{\G_t}(p)\mu(p).
\]
Meanwhile,  we can use the function $d_{\G_t}(p)$ as the expression of the upper surface $\G_t$, and we have at the corner points that 
\[
d_{\G_t}(p_{i*})=p_i,\quad i=l, r.
\] 

Moreover, $\Phi_{S_t}$ can be extended to  the entire boundary $S_*=\G_{t*}\cup \G_{b*}$. 
Consequently,  we obtain the map on $S_*$: 
\[
\Phi_{S_t}: \ \ S_* \rightarrow  {S_t}.
\]
Using the harmonic extension, we define the following map on $\Om_*$:
\[
\cT_{S_t}:\ \Om_*\rightarrow \Om\quad\hbox{with}\quad \cT_{S_t}=\cH_*\big(\Phi_{S_t}-Id_{S_*}\big)+Id.
\]
Here $\cH_*\big(\Phi_{S_t}-Id_{S_*}\big)$ is the harmonic extension of $\Phi_{S_t}-Id_{S_*}$ satisfying 
\[
\left\{\begin{array}{ll}
\D \cH_*\big(\Phi_{S_t}-Id_{S_*}\big)=0\qquad\hbox{in}\quad \Om_*,\\
\cH_*\big(\Phi_{S_t}-Id_{S_*}\big)\big|_{\G_{t*}}=d_{\G_t}\mu,
\quad \cH_*\big(\Phi_{S_t}-Id_{S_*}\big)\big|_{\G_{b*}}=\Phi_{S_t}|_{\G_{b*}}-Id_{\G_{b*}}.
\end{array}\right.
\]



\subsection{Recovery of the velocity}

When the domain $\Om_t$ is defined by $\cT_S$, we can define the velocity $v$ by the free surface function $d_{\G_t}$.
To begin with,  the kinematic condition on $\G_t$ in $\mbox{(WW)}$ is rewritten into 
\[
\pa_t \Phi_{S_t}\cdot (n_t\circ \Phi_{S_t})=(v\cdot N_t)\circ \Phi_{S_t}\quad\hbox{with}\quad\pa_t\Phi_{S_t}=(\pa_t d_{\G_t})\mu\quad\hbox{on}\quad \G_*.
\]
So we obtain 
\beq\label{normal derivative of v}
\pa_t d_{\G_t}=\f{(v\cdot n_t)\circ \Phi_{S_t}}{\mu\cdot (n_t\circ \Phi_{S_t})}\quad\hbox{i.e.} \quad
v\cdot n_t=(\pa_t d_{\G_t}\mu)\circ \Phi^{-1}_{S_t}\cdot n_t.
\eeq

Due to the  assumption that the velocity $v$ is  irrotational,  we define $v$ by  
\beq\label{v potential expression}
v=\na \phi
\eeq
with $\phi$ satisfying
\beq\label{phi system}
 \left\{\begin{array}{ll}
\Delta \phi=\xi \gamma\qquad\hbox{in}\quad \Om_t,\\
\na_{n_t} \phi|_{\Gamma_t}= (\pa_t d_{\G_t} \,\mu)\circ \Phi^{-1}_{S_t}\cdot n_t, \quad \na_{n_b} \phi|_{\Gamma_b}= 0,
\end{array}\right.
\eeq
where
\[
\gamma=|\Om_t|^{-1}\quad\hbox{ and }\quad \xi=\int_{\Gamma_t} v\cdot n_t\,ds=\int_{\G_t}(\pa_t d_{\G_t}\mu)\circ \Phi^{-1}_{S_t}\cdot n_tds.
\]
 
Moreover, we define as in \cite{MW3} that
\[
D_{t^*}= \pa_t+\na_{v^*},
\]
where $v^{*}= D\Phi_{S_t}^{-1}(v^\top\circ \Phi_{S_t} -\pa_t d_{\G_t} \mu^{\top} )$. A direct computation shows that
\[
(D_t f)\circ \Phi_{S_t}=D_{t^*}(f\circ \Phi_{S_t}) ,
\]
for a function $f$ on $\G_t$.

\subsection{The modified formulation  and the precise form of the main theorem}
%
%
Before we construct the approximate solutions, we derive a new equation modified from the previous sections, which turns out to be more convenient in this part. 

First, we define $\mK_a$ based on  the definition of $\mathfrak{K}$:
\beq\label{Ka def}
\mathfrak{K}_a=\mK+\si a \,d_{\G_t}\circ \Phi_{S_t}^{-1}= \si(\kappa+a \,d_{\G_t}\circ \Phi_{S_t}^{-1} )-P_{v, v} \quad\hbox{on} \quad\G_t
\eeq
for some constant $a>0$, where $\kappa$ can be expressed by $d_{\G_t}$ and $P_{v,v}$ is defined by \eqref{eq: P_{vv}}. 
So we have
\[
\mJ_a= \na \cH (\mathfrak{K}_a)=\mJ+\si a\na \cH(d_{\G_t}\circ \Phi_{S_t}^{-1}).
\]

Applying  \eqref{Dt2 kappa} and using $\mK_a$ instead of $\mK$, we obtain
\[
    D_t^2 \mathfrak{K}_a+\sigma (a-\D_{\G_t})\cN(\mathfrak{K}_a)=R_{a, 0},
\]
    where $R_{a,0}$ is defined by
\[
\begin{split}
R_{a, 0}&=\si R_1 -D_t^2P_{v,v}+\si[n_t, \D_{\G_t}]\cdot \mJ+\si \D_{\G_t}\na P_{v,v}\cdot n_t+2\si\na_{\tau_t}n_t\cdot \na_{\tau_t}\mJ+\si a D^2_t(d_{\G_t}\circ \Phi_{S_t}^{-1})\\
&\quad -\si^2 a \D_{\G_t}\cN(d_{\G_t}\circ \Phi_{S_t}^{-1})+\si a \cN\mK_a.
\end{split}
\]

    Acting $\na \cH$ on both sides of the above equation, we  get the equation of $\mJ_a$:
\[
D^2_t\mJ_a+ \sigma\na \cH\big[(a-\D_{\G_t})\mJ_a^\perp  +h_v\big]=R_a,
\]
where  
\[
R_a= \na \cH(R_{a, 0}+\si h_v)-[\na\cH, D_t^2]\mathfrak{K}_a,
\]
and recall that $h_v$ is defined in \eqref{hv def} and comes from $R_1$.

Meanwhile, the condition \eqref{equ:cc} at the contact points are rewritten as 
\[
 D_t\mJ_a = \pm\f{\sigma^2 }{\beta_c}\sin\om_i (\na_{\tau_t} \mJ_a)^\bot
 \tau_b+R_{c1, a}\quad\hbox{at}\quad p_i (i=l, r)
\]
with 
\[
R_{c1, a}=R_{c1}+\si aD_t\na \cH(d_{\G_t}\circ \Phi_{S_t}^{-1})+\si a \f{\sigma^2 }{\beta_c}(n_t\cdot \tau_b) \big(\na_{\tau_t} \na \cH(d_{\G_t}\circ \Phi_{S_t}^{-1})\big)^\perp.
\]
%
%
%

%
%

\medskip
    
Next, we consider how to recover  the free surface and the domain from  $\mathfrak{K}_a$, which is slightly different from \cite{MW3}. In \cite{MW3}, we use the equation of $\mathfrak N_a=\cN(\kappa+a \,d_{\G_t}\circ \Phi_{S_t}^{-1})$ together with the boundary information $d_i=d_{\G_t}|_{p_{i*}}$ of $d_{\G_t}$ to recover the free surface $d_{\G_t}$, so the system of $(\mathfrak{K}_a, d_l, d_r)$ is needed; In this paper, we use the quantity $\mK_a$ instead of $\mathfrak N_a$, where an extra $P_{v,v}$ is added here in \eqref{Ka def}. 
Moreover, we will need to use the equation and norms of $\mJ_a$ in the  iteration scheme, where  $\mK_a$ can be retrieved. In fact, we have $\mJ_a^\perp=\mJ_a\cdot n_t=\cN \mK_a$ on $\G_t$. Therefore,  to identify $\mK_a$, we look at  the Neumann-boundary elliptic system of $\cH(\mK_a)$ with the compatibility condition $\int_{\G_t}\mJ_a^\perp ds=0$. We know immediately that there exists a unique solution $\cH(\mK_a)$ up to an additive constant to this system. 
As a result,  as long as  we have $\mJ_a$, we obtain $\mathfrak{K}_a$. (One can also check Lemma 2.5 in \cite{MW3}.)

Consequently, to recover $d_{\G_t}$ (which is the key to recover the water-waves system), we need the system of $(\mK_a, P_{v,v}, d_l, d_r)$. As long as we have proved the existence of the solution to this system, we obtain immediately the following quantity 
\[
\kappa+a \,d_{\G_t}\circ \Phi_{S_t}^{-1}=\si^{-1}(\mK_a+P_{v,v})\quad\hbox{with the boundary information} \ d_l, d_r.
\]
As a result,  the the desired function $d_{\G_t}$ can be solved directly from these quantities above  in a similar way as in Proposition 4.2 \cite{MW3}, and then we can finally recover our water-waves system $\mbox{(WW)}$.

\medskip   
Based on the analysis above,  we need to give the boundary condition of $d_{\G_t}$ which is deduce from \eqref{normal derivative of v} (for more details, see (4.27) in \cite{MW3}).  In fact, one has the following evolution equations for $d_i(t)=d_{\G_t}(p_{i*})$ ($ i=l, r$):
\beq\label{d c eqn}
 d_i''(t)=\mathfrak{B}_i,\quad i=l, r,
\eeq
where
\[
\begin{split}
\mathfrak{B}_i&=-\f{1}{\mu\cdot (n_t \circ\Phi_{S_t})} \Big(\mu\cdot (n_t \circ\Phi_{S_t})\na_{v^*}\pa_td_{\G_t}+\na_{v^*}\mu\cdot (n_t\circ\Phi_{S_t})\pa_t d_{\G_t}+\si\mathfrak{K}_a\circ \Phi_{S_t}\\
&\quad -\si^2 a \cN (d_{\G_t}\circ \Phi^{-1}_{S_t})\circ \Phi_{S_t}+(\na P_{v,v}+{\bf g})\circ \Phi_{S_t}\cdot (n_t\circ \Phi_{S_t})\Big)\big|_{p_{i*}}.
\end{split}
\]

We rewrite the equation for $P_{v,v}$ by \eqref{DtPvv system}:
\[
D_tP_{v,v}=\D^{-1}_N(h_p, f_p, g_p),
\]
where $\D^{-1}_N$ means solving the Neumann-boundary system \eqref{DtPvv system} with $\int_{\Om_t}D_tP_{v,v}dX=0$, and the right-side functions are
\[
h_p=2tr (\na v\cdot \na v\, \na v)+2tr \big(\na(\na P_{v,v}+\mJ)\na v\big)+ 2tr(\na  v  \na^2  P_{v, v}),\quad  f_p=C'_{v,v}(t)+\na_{n_t} v\cdot \na P_{v, v}\big|_{\G_t}\]
and 
\[
g_p=-(mJ+\na P_{v,v}+{\bf g})\cdot\big (\na n_b+(\na n_b)^*\big)\cdot v+v\cdot D_t(\na n_b)\cdot v+\na_{n_b} v\cdot \na P_{v, v}.
\]
\medskip

As a result, we sum up the system of $(\mK_a, P_{v,v}, d_l, d_r)$ as follows:
\beq\label{eq: IS-C}
\left\{\begin{array}{ll}
    D_{t}^2 \mathfrak{K}_a +\sigma (a-\D_{\G_t})\cN(\mathfrak{K}_a)=R_{a, 0},\\
 D_{t}\mJ_a  = \pm\sigma^2 \beta^{-1}_c\sin\om (\na_{\tau_t} \mJ_a)^\bot 
 \tau_b +R_{c1, a} , \quad \textrm{at} \quad  p_{i} ( i=l, r),\\
 D_t P_{v,v}=\D^{-1}_N(h_p, f_p, g_p),\\
 \f {d^2}{dt^2}d_i(t)=\mathfrak{B}_i,\quad i=l, r.
\end{array}
\right.
\eeq

Based on these preparations above, we are finally ready to state our precise form of Theorem \ref{thm:main}. We start with introducing the space $\Sigma$ for given $T, L>0$:
\[
\Sigma=\big\{(\mK_a, P_{v,v}, d_l, d_r)\big| \|(\mK_a, P_{v,v}, d_l, d_r)\|_\Sigma\le L\big\}
\] 
where the  norm
\[
\begin{split}
\|(\mK_a, P_{v,v}, d_l, d_r)\|_\Sigma&\triangleq \|\pa_t(\mJ_a\circ \Phi_{S_t})  \|_{C([0, T]; L^2(\Om_*))}
+\| (\mJ_a\circ \Phi_{S_t})\cdot n_{t*}\|_{C([0, T]; H^1(\Gamma_*))} \\
&\quad +\|P_{v,v}\circ  \Phi_{S_t}\|_{C([0, T]; H^{5/2}(\Om_*))} +\sum_{i=l, r}\Big(|d_i|_{C([0, T])}+| d(d_i)/ dt |_{C([0, T])}\Big).
\end{split}
\] 
Meanwhile, according to the higher-order energy $E_h(t)$, we also define the space (for given $L_1>0$)
\[
\Sigma_h=\big\{(\mK_a, P_{v,v}, d_l, d_r)\big| \|(\mK_a, P_{v,v}, d_l, d_r)\|_{\Sigma_h}\le L_1\big\}
\]
by the norm 
\[
\begin{split}
\|(\mK_a, P_{v,v}, d_l, d_r)\|_{\Sigma_h}&\triangleq \|(\mK_a, P_{v,v}, d_l, d_r)\|_\Sigma+\|\pa_t\big((\mJ_a\circ \Phi_{S_t}) \cdot n_{t*}\big)  \|_{L^\infty([0, T]; H^1(\G_*))}\\
&\quad +\|\pa^2_t(\mJ_a\circ \Phi_{S_t})  \|_{L^\infty([0, T]; L^2(\Om_*))}.
\end{split}
\]
The initial data is given by
\[
\left\{\begin{array}{ll}
(\mJ_a\circ \Phi_{S_t})\cdot n_{t*}\big|_{t=0}=\bar \mJ_{a, 0},  \quad \pa_t(\mJ_a\circ \Phi_{S_t})\big|_{t=0} =\bar\mJ_{a, 1},\\
P_{v,v}\circ \Phi_{S_t}\big|_{t=0}=\bar P_{v, v, 0}, \quad d_i(0)=d_{i, 0},\quad  \f d{dt}d_i(0)=d_{i, 1}.
\end{array}
\right.
\]

Now we can  present  the local well-posedness theorem.
\begin{thm}\label{main thm} Assume that the  initial data $\big(\bar \mJ_{a, 0}, \bar\mJ_{a, 1}, \bar P_{v, v, 0}, d_{i,0}, d_{i, 1}\big)\in H^1(\G_*)\times L^2(\Om_*)\times  H^{5/2}(\Om_*)\times \R^2$ and initial contact angles $\om_{i0}\in (0, \pi/2) $ for $i=l, r$.  When the  compatibility conditions  \eqref{eq:com cond} at $t=0$  are satisfied for $k=0,1,2,3$, there exists a unique solution $(\mK_a, P_{v,v}, d_l, d_r)\in \Sigma_h$ to system \eqref{eq: IS-C}. Moreover, system \eqref{eq: IS-C} is locally well-posed with  $(\mK_a, P_{v,v}, d_l, d_r)$ depending continuously on the initial data in  $\Sigma$.

\end{thm}

   \subsection{Iteration scheme}
In this subsection, we present the iteration scheme. First of all, we set the initial boundary $S_0=S_*$ without loss of generality.  To simplify the notations, we denote by
\[
 D_{t*}=\pa_t+v^k_*\cdot \na,\quad D_t=\pa_t+v^k\cdot \na,
\]
when no confusion will be made.
 
When we have $(\mK^{k}_a, P^{k}_{v,v}, d^{k}_l, d^{k}_r)$, the linear system of $(\mK^{k+1}_a, P^{k+1}_{v,v}, d^{k+1}_l, d^{k+1}_r)$ for the iteration scheme is set to be
 \beq\label{eq: IS}
\left\{\begin{array}{ll}
    D_{t*}^2 (\mathfrak{K}^{k+1}_a\circ \Phi^k_{S_t})+\sigma (a-\D_{\G_t})\cN(\mathfrak{K}^{k+1}_a)\circ \Phi^k_{S_t}=R^k_{a, 0}\circ \Phi^k_{S_t},\\
 D_{t*}(\mJ^{k+1}_a\circ \Phi^k_{S_t}) = \pm\sigma^2\beta^{-1}_c\sin\om^k_i (\na_{\tau^k_t} \mJ_a^{k+1})^\bot 
 \tau_b\circ \Phi^k_{S_t}+R^k_{c1, a}\circ \Phi^k_{S_t}\quad \textrm{at} \quad  p_{i*},\quad i=l, r,\\
 D_{t*}(P^{k+1}_{v,v}\circ \Phi^k_{S_t})=\D^{-1}_N(h^k_p, f^k_p, g^k_p)\circ \Phi^k_{S_t}\\
 \f {d^2}{dt^2}d^{k+1}_i(t)=\mathfrak{B}^k_i,\qquad i=l, r.
\end{array}
\right.
\eeq
Here we use the superscript  $k$ on $R_{a, 0}$ (for example) to denote that all the quantities there are obtained using $(\mK^{k}_a, P^{k}_{v,v}, d^{k}_l, d^{k}_r)$. 

Moreover, we point out that the velocity $v$ in the definition of $R^k_{a, 0}$ is given by
\[
D_{t}\tilde{v}^k =-J^k-\na P_{v^k, v^k}-{\bf g},
\]
while in the other quantities we use $v^k$ defined by \eqref{normal derivative of v}-\eqref{phi system}. This happens due to the difference of regularities using these two definitions, which can be seen already in the previous sections.

Besides, the initial data is given by
\beq\label{eq: ID-IC}
\left\{\begin{array}{ll}
(\mJ^{k+1}_a\circ \Phi^k_{S_t})\cdot n^k_{t*}\big|_{t=0}=\bar \mJ_{a, 0},  \quad \pa_t(\mJ^{k+1}_a\circ \Phi^k_{S_t})\big|_{t=0} =\bar\mJ_{a, 1},\\
P^{k+1}_{v,v}\circ \Phi^k_{S_t}\big|_{t=0}=\bar P_{v, v, 0}, \quad d^{k+1}_i(0)=d_{i, 0},\quad  \f d{dt}d^{k+1}_i(0)=d_{i, 1}.
\end{array}
\right.
\eeq
 
As a result, a similar proof as the proof of Proposition 5.1 in \cite{MW3}, we show the existence of the solution $(\mK^{k+1}_a, P^{k+1}_{v,v}, d^{k+1}_l, d^{k+1}_r)$ to  the linear system \eqref{eq: IS}--\eqref{eq: ID-IC}. The details for the proof are omitted. 
 \begin{prop}\label{existence of linear system} 
Let $\big(\bar \mJ_{a, 0}, \bar\mJ_{a, 1}, \bar P_{v, v, 0}, d_{i,0}, d_{i, 1}\big)\in H^1(\G_*)\times L^2(\Om_*)\times  H^{5/2}(\Om_*)\times \R^2$, $\om_{i0}\in (0, \pi/2)$ and $(\mK^{k}_a, P^{k}_{v,v}, d^{k}_l, d^{k}_r)$ be given correspondingly. Moreover we assume that the conditions  for the corner points from \eqref{eq: IS} hold at $t=0$.   Then there exists a small $T>0$ such that the system \eqref{eq: IS}--\eqref{eq: ID-IC} has a unique solution on $[0, T]$.
\end{prop}

      \subsection{Uniform estimates}

Now  we are ready to give the uniform estimates for the linear system \eqref{eq: IS}--\eqref{eq: ID-IC}. To begin with, we define the energy functional for $k\in \N$ as below:
\[
 E^{k+1}(t)= E^{k+1}_{low}(t)+ E^{k+1}_{high}(t)
\]
 where $E^{k+1}_{low}(t)$ and $E^{k+1}_{high}(t)$ are defined by
\[
\begin{split}
E^{k+1}_{low}(t)&=a\| (\mJ^{k+1})^\bot\|^2_{L^2(\Gamma^k_t)}+\|\na_{\tau^k_t}  (\mJ^{k+1})^\bot\|^2_{L^2(\Gamma^k_t)} +\| D_t  \mJ^{k+1}\|^2_{L^2(\Om^k_t)} +\|P^{k+1}_{v,v}\|^2_{H^{5/2}(\Om^k_t)}\\
 &\quad +\sum_{i=l, r}\big(|d^{k+1}_i(t)|^2+\Big|\f d{dt}d^{k+1}_i(t)\Big|^2\big),
\end{split}
\]
 and
\[
E_{high}(t)=a\|D_t(\mJ^k)^\bot\|^2_{L^2(\Gamma^k_t)}+\|\na_{\tau^k_t}D_t(\mJ^{k+1})^\bot\|^2_{L^2(\Gamma^k_t)} +\| D_t^2  \mJ^{k+1} \|^2_{L^2(\Om^k_t)}.
\]
Moreover,  the dissipation $F^{k+1}(t)$ is 
\[
F^{k+1}(t)= F^{k+1}_{low}(t)+ F^{k+1}_{high}(t),
\]
where 
\[
F^{k+1}_{low}(t)=\sum_{i=l,r}\big|(\sin \om_i^k)\na_{ \tau^k_t}  (\mathfrak{J}^{k+1})^\perp |_{p_i}\big|^2
,\quad
F^{k+1}_{high}(t)=\sum_{i=l,r}\big|(\sin \om_i^k)\na_{ \tau^k_t}D_t (\mathfrak{J}^{k+1})^\perp |_{p_i}\big|^2.
\]

Meanwhile, we define some more auxiliary functions. Recalling from \eqref{eq:cP_1}, here $\cP^k_{t,1}$ and $\cP^k_{t,2}$ are defined by
\[
\cP^k_{t,1}=D_t P_{\tilde{v}^k, \tilde{v}^k}-v^k\cdot (\na P_{\tilde{v}^k, \tilde{v}^k}|_{c}),\quad
\cP^k_{t,2}=D_t\cP^k_{t,1}-v^k\cdot (\cP^k_{t,1}|_{c}).
\]

The following proposition is our main result on the uniform estimates:
\begin{prop}
Let $\big(\bar \mJ_{a, 0}, \bar\mJ_{a, 1}, \bar P_{v, v, 0}, d_{i,0}, d_{i, 1}\big)$ and $(\mK^{k}_a, P^{k}_{v,v}, d^{k}_l, d^{k}_r)$ be given as in Proposition \ref{existence of linear system}. Then there exists constants $T>0$ small enough and $A>0$ large enough such that when $a\geq A$, the inequality below holds
\beno
\sup_{t\in [0, T]} E^{k+1}(t)+\int_0^T F^{k+1}(t) dt \leq   P(E(0)).
\eeno
\end{prop}
\no{\bf Proof.}
Since the main steps of the proof follow Theorem \ref{lower-order energy estimates} and Theorem \ref{1-order energy estimates}, we only present the sketch of the  proof here.

First, we consider the basic energy estimates $E^{k+1}_l$, where we only focus on the estimates for $\mK_a$ or $\mJ_a$ and the other estimates follow from Lemma \ref{est: P_{v, v}-H} and \cite{MW3}. 

To begin with, acting $(\Phi^k_{S_t})^{-1}$ on both sides of \eqref{eq: IS}, one has
\[
D_{t}^2 \mathfrak{K}^{k+1}_a +\sigma (a-\D_{\G_t})\cN(\mathfrak{K}^{k+1}_a) =R^k_{a, 0} ,
\]
which implies that
  \beq\label{eq: IS-1}
D_t^2 \mJ^{k+1}_a  +\sigma\na \cH\big[(a-\D_{\G_t})(\mJ^{k+1}_a)^\perp+h_{\tilde{v}^k}\big]  = \na \cH(R^k_{a, 0}+\si h_{\tilde{v}^k})-[\na\cH, D_t^2]\mathfrak{K}^{k+1}_a.
\eeq
Using the same arguments as in Section 4 and Section 5, we have 
\[
 \big\|\na \cH(R^k_{a, 0}+\si h_{\tilde{v}^k})-[\na\cH, D_t^2]\mathfrak{K}^{k+1}_a\big\|_{L^2(\Om_t^k)}\leq P(E^{k}(t))(1+E^{k+1}_{low}(t)).
\]

Next,  we define  $P_{J^{k+1}, v^k}$ on $\Om_t^k$ by  system \eqref{PJv system}. Taking the $L^2(\Om^k_t)$ inner product of  \eqref{eq: IS-1} with $D_t\mJ^{k+1}_a +\na P_{J^{k+1}, v^k}$,  one has 
\[
\begin{split}
&\int_{\Om_t^k} \Big( D_t ^2\mJ^{k+1}_a  + \na\cH\big[(a-\D_{\G^k_t})(\mJ^{k+1}_a)^\perp+h_{\tilde{v}^k}\big]\Big)\cdot(D_t \mJ^{k+1}_a +\na P_{J^{k+1}, v^k})dX\\
&=\int_{\Om_t^k} \big(\na \cH(R^k_{a, 0}+\si h_{\tilde{v}^k})-[\na\cH, D_t^2]\mathfrak{K}^{k+1}_a\big) \cdot(D_t \mJ^{k+1}_a +\na P_{J^{k+1}, v^k})dX.
\end{split}
\]
Following the energy estimates in Section 4, there exists a constant  $A$ large enough such that when $a\geq A$, one concludes
\[
\sup_{t\in[0, T]} E^{k+1}_{low}(t)+\int_0^T F^{k+1}_{low}(t) \leq   P(E(0))+\int_0^TP(E^k(t)) .
\]

On the other hand, for the higher-order energy $E_h^{k+1}(t)$, we get similarly as \eqref{eq: D_t^2 J-1} the equation
\[
D_t \big[D^2_t \mathfrak{J}_a^{k+1} -\na \cH(\cP^k_{t,2}) +[\na\cH, D_t^2]\mathfrak{K}^{k+1}_a\big]+\si \na\cH\big((a-\D_{\Gamma^k_t})  D_t( \mathfrak{J}_a^{k+1})^\bot+D_t h_{\tilde{v}^k} \big)= R^k_{a, 1}  ,
\]
where $R^k_{a, 1}$ is given by 
\[
\begin{split}
R^k_{a, 1}&=D_t \na \cH (R^k_{a,0}+\si h_{\tilde v_k})-D_t\na \cH(\cP^k_{t,2}) +\si [\na\cH, D_t]\cH \big((a-\D_{\G_t})(\mJ^{k+1}_a)^\perp+\si h_{\tilde{v}^k}\big) 
\\
&\quad -\si\na \cH [\D_{\Gamma_t}, D_t] (\mJ^{k+1}_a)^\perp.
\end{split}
\]

Compared to the energy estimate of Theorem \ref{1-order energy estimates}, the main difference lies in  the term $D_t ([\na\cH, D_t^2]\mathfrak{K}^{k+1}_a)$ which contains a  term like  $D_t(\pa D_t v^k)\pa \cH(\mathfrak{K}^{k+1}_a)$. In Theorem \ref{1-order energy estimates}, since $v$ satisfies Euler's equation, we have the estimate for $D_t^2\na v$. But here in the iteration scheme,  $D_t^2\pa v^k$ acts like $\pa\pa_t^3 d_{\G_t}$, which can not be controlled by the energy. Therefore, we put $D_t ([\na\cH, D_t^2]\mathfrak{K}^{k+1}_a)$ together with $D^2_t \mathfrak{J}_a^{k+1}$ to go back to the form $D_t^2 \mathfrak{K}^{k+1}_a$:
\[
D^2_t \mathfrak{J}_a^{k+1}+[\na\cH, D_t^2]\mathfrak{K}^{k+1}_a=\na \cH(D_t^2 \mathfrak{K}^{k+1}_a). 
\]
Taking the $L^2(\Om^k_t)$ inner product of the above equation with $\na \cH(D_t^2 \mathfrak{K}^{k+1}_a) -\na \cH(\cP^k_{t,2})  $, we derive
\[
\begin{split}
&\int_{\Om_t^k}D_t\big(\na \cH(D_t^2 \mathfrak{K}^{k+1}_a) -\na \cH(\cP^k_{t,2})\big)\cdot \big(\na \cH(D_t^2 \mathfrak{K}^{k+1}_a) -\na \cH(\cP^k_{t,2})\big)dX\\
&+\si\int_{\Om_t^k}\na\cH\big((a-\D_{\Gamma^k_t})  D_t( \mathfrak{J}_a^{k+1})^\bot-D_t h_{\tilde{v}^k} \big)\cdot \big(\na \cH(D_t^2 \mathfrak{K}^{k+1}_a) -\na \cH(\cP^k_{t,2})\big)dX\\
&=\int_{\Om_t^k}R^k_{a, 1}  \cdot \big(\na \cH(D_t^2 \mathfrak{K}^{k+1}_a) -\na \cH(\cP^k_{t,2})\big)dX.
\end{split}
\]
 
  For the second integral on the left side of the above equation, we have by Green's Formula that
\[
\begin{split}
 &\int_{\Om_t^k}\na\cH\big((a-\D_{\Gamma^k_t})  D_t( \mathfrak{J}_a^{k+1})^\bot+ D_t h_{\tilde{v}^k} \big)\cdot \big(\na \cH(D_t^2 \mathfrak{K}^{k+1}_a) -\na \cH(\cP^k_{t,2})\big)dX\\
 &=\int_{\G^k_t} \big((a-\D_{\Gamma^k_t})  D_t( \mathfrak{J}_a^{k+1})^\bot+D_t h_{\tilde{v}^k} \big)\cdot \cN(D_t^2 \mathfrak{K}^{k+1}_a-\cP^k_{t,2})ds\\
 &=\int_{\G^k_t} \big((a-\D_{\Gamma^k_t})  D_t( \mathfrak{J}_a^{k+1})^\bot+ D_t h_{\tilde{v}^k} \big)\cdot  D_t^2 (\mathfrak{J}_a^{k+1})^\bot ds\\
 &\quad+\int_{\G^k_t} \big((a-\D_{\Gamma^k_t})  D_t( \mathfrak{J}_a^{k+1})^\bot+ D_t h_{\tilde{v}^k} \big)\cdot \big([\cN,D_t^2] \mathfrak{K}^{k+1}_a-\cN(\cP^k_{t,2})\big)ds,
\end{split}
\]
and proceeding in a similar way as before, we obtain
\[
\begin{split}
 &\int_{\Om_t^k}\na\cH\big((a-\D_{\Gamma^k_t})  D_t( \mathfrak{J}_a^{k+1})^\bot+ D_t h_{\tilde{v}^k} \big)\cdot \big(\na \cH(D_t^2 \mathfrak{K}^{k+1}_a) -\na \cH(\cP^k_{t,2})\big)dX\\
  &=\f d{dt}\Big(a\|D_t(\mJ^k)^\bot\|^2_{L^2(\Gamma^k_t)}+\big\|\na_{\tau^k_t}D_t(\mJ^{k+1})^\bot\big\|^2_{L^2(\Gamma^k_t)} \Big)+\na_{\tau^k_t}D_t(\mJ^{k+1}_a)^\perp\cdot D^2_t( \mathfrak{J}_a^{k+1})^\perp\big|^{p_r}_{p_l} \\
 &\quad+\int_{\G^k_t} (a-\D_{\Gamma^k_t})  D_t( \mathfrak{J}_a^{k+1})^\bot  \cdot \big([\cN,D_t^2] \mathfrak{K}^{k+1}_a-\cN(\cP^k_{t,2})\big)ds\\
 &\quad+\int_{\G^k_t}   D_t h_{\tilde{v}^k} \cdot \big(D_t^2 (\mathfrak{J}_a^{k+1})^\bot+[\cN, D_t^2] \mathfrak{K}^{k+1}_a-\cN(\cP^k_{t,2})\big)ds.
 \end{split}
\]
 On one hand, since the commutator $[D_t,\,\cN]$ is already expressed in \eqref{commutator DN}, 
one can conclude that
\[
 \|[\cN,D_t^2] \mathfrak{K}^{k+1}_a\|_{H^1(\G_t^k)}\leq  P(E^{k+1}(t)).
\]
Consequently, we have
\[
 \int_{\G^k_t} \big((a-\D_{\Gamma^k_t})  D_t( \mathfrak{J}_a^{k+1})^\bot \big)\cdot \big([\cN,D_t^2] \mathfrak{K}^{k+1}_a-\cN(\cP^k_{t,2})\big)ds\leq \f18 F_h^{k+1}(t)+ P(E^{k+1}(t)).
\]
On the other hand,  using similar arguments as in  Lemma \ref{lem: f-D_tJ-1}, we get
\[
\begin{split}
&\int_0^T \int_{\G^k_t}  D_t h_{\tilde{v}^k}  \cdot \big(D_t^2 (\mathfrak{J}_a^{k+1})^\bot+[\cN, D_t^2] \mathfrak{K}^{k+1}_a-\cN(\cP^k_{t,2})\big)dsdt\\
&\leq P(E(0))+\sup_{t\in[0, T]}\big(\f18 E_{high}^{k+1}(t)+P(E_{low}^{k+1}(t))\big) +\int_0^T P(E^{k}(t))dt+ \f{1}{8} \int_0^T F_{high}^{k+1}(t)dt.
\end{split}
\]

In the end, using similar arguments as in Section 5.3, we can show  that when $T$ small enough, the following estimate holds
\[
\sup_{t\in[0, T]} E^{k+1}_{high}(t)+\int_0^T F^{k+1}_{high}(t)dt \leq   P(E(0))+\int_0^TP(E^k(t) )dt +\sup_{t\in[0, T]} P(E_{low}^{k+1}(t)) .
\]
Therefore, by a bootstrap argument, we can prove the desired result. 

\ef

\subsection{Cauchy sequence and going back to $\mbox{(WW)}$}
      
 In this part, we are finally in a position to prove that the sequence of $(\mathfrak{K}^{k}_a, P^k_{v,v}, d^{k}_l, d^{k}_r)$ is indeed a Cauchy sequence.   In fact, we know from the previous subsection that
\beq\label{bound of Ek}
\sup_{t\in[0, T]} E^{k+1}(t) +\int_0^T F^{k+1}(t)dt \leq C\quad\hbox{for all} \ k\in\mathbb N,
\eeq
where $C>0$ is a constant depending on $E(0)$.

 
 To simplify the notation, we denote by
\[
\overline{f}^{k+1}=f^{k+1}\circ \Phi^k_{S_t}  ,\qquad \d_{\overline{f}^k} = \overline{f}^{k+1}- \overline{f}^k,\]
and
\[
(a-\D_{\G_t})g^k=h_{\tilde{v}^k}\quad\hbox{with}\quad g^k|_{p_i}=0.
\]
     
Using \eqref {eq: IS-1} and rewriting it similarly as \eqref{eqn cDt2 J}, we have the equation of $\d_{\overline{\mJ}^{k}_a} $:
\beq\label{eq:d}
D_{t*}(D_{t*}\d_{\overline{\mJ}^{k}_a}-\delta_{\overline{\na \cH(\cP^k_{t,1})}})+\si \cA(d_{\G^k_t})(\d_{\overline{\mJ}^{k}_a}-\d_{\overline{g}^k}) = D_{R^k},
\eeq
where we note
\[
\cA(d_{\G_t})f=\big(\cN(a-\D_{\G_t})(f\circ \Phi^{-1}_{S_t})\big)\circ \Phi_{S_t},
\]
and 
\[
\begin{split}
D_{R^k}&=\big((\pa_t+v^k_*\cdot\na)^2-(\pa_t+v^{k-1}_*\cdot\na)^2\big)\overline \mJ^k_a+\big(\cA(d_{\G^k_t})-\cA(d_{\G^{k-1}_t})\big)(\overline{\mJ}^{k}_a-\overline{g}^k) +\d_{ \overline{\na \cH(R^{k-1}_{a, 0}+\si h_{\tilde{v}^{k-1}})}}\\
&\quad-\d_{ \overline{[\na\cH, D_t^2]\mathfrak{K}^{k}_a}}-D_{t*}\delta_{\overline{\na \cH(\cP^k_{t,1})}}.
\end{split}\]
Besides, similar but simpler equations for $(\d_{P^k_{v,v}}, \d_{d^k_l}, \d_{d^k_r})$ can be derived, and we omit the details here.

Moreover, we define the energy of the difference according to the definition of $\Sigma$ as  
\[
\begin{split}
E^k_\d(t)&= \|\pa_t\d_{\overline{\mJ}^{k}_a}  \|_{C([0, T]; L^2(\Om_*))}+a\| \d_{\overline{\mJ}^{k}_a}\cdot n_{t*}\|_{C([0, T]; L^2(\Gamma_*))}
+\| \d_{\overline{\mJ}^{k}_a}\cdot n_{t*}\|_{C([0, T]; H^1(\Gamma_*))} \\
&\quad +\|\d_{P^k_{v,v}}\|_{C([0, T]; H^{5/2}(\Om_*))} +\sum_{i=l, r}\big(|\d_{d^{k}_i}|_{C([0, T])}+\Big|\d_{d (d^{k}_i)/dt}\Big|_{C([0, T])}\big).
\end{split}
\]

Before we consider the convergence, we need to deal with $D_{R^k}$ first.
\begin{lemma} 
The right side of \eqref{eq:d} satisfies the following estimate:
\[
\|D_{R^k}\|_{L^2(\Om_*)}\leq CE^k_\d(t)
\]
with the positive constant $C$ depending on $E(0)$.
\end{lemma}
\no{\bf Proof.}
Here, we only give the outline of the proof, and one can see similar details in \cite{MW3}. 
First, we consider  the estimate for $\big((\pa_t+v^k_*\cdot\na)^2-(\pa_t+v^{k-1}_*\cdot\na)^2\big)\overline \mJ^k_a$. Using similar arguments as  Proposition 4.12 and Corollary 4.13 in \cite{MW3} and applying \eqref{eq: IS} and  \eqref{bound of Ek}, we can have 
\[
\begin{split}
\big\|\big((\pa_t+v^k_*\cdot\na)^2-(\pa_t+v^{k-1}_*\cdot\na)^2\big)\overline \mJ^k_a\big\|_{L^2(\Om_*)} &\leq C\|\pa_t^2\d_{d_{\G^{k-1}_t}}\|_{L^2(\G_*)}+C\|\pa_t\d_{d_{\G^{k-1}_t}}\|_{H^1(\G_*)} \\
 &\leq C \big(a^{-1} \|\d_{\overline{\mJ}^{k}_a}\cdot n_*\|_{H^1(\G_*)}+C \|\pa_t\d_{\overline{\mJ}^{k}_a} \|_{L^2(\Om_*)}\big)\leq C E^k_\d(t).
\end{split}
\]

Second, for the term $(\cA(d_{\G_t})-\cA(d_{\G^{k-1}_t}))(\overline{\mJ}^{k}_a-\overline{g}^k)$, we notice by \eqref{D J perp-hv estimate} that
\[
\big\| \overline{\mJ}^{k}_a-\overline{g}^k\big\|_{H^{5/2}(\G_*)}\leq C.
\]
As a result, by similar arguments as in the proof of  Proposition 6.3 \cite{MW3}, we have
\[
\big\|\big(\cA(d_{\G_t})-\cA(d_{\G^{k-1}_t})\big)(\overline{\mJ}^{k}_a-\overline{g}^k)\big\|_{L^2(\Om_*)}\leq CE^k_\d(t).
\]

In the end, similar arguments as  in Lemma 4.6 \cite{MW3}, we can have
\[
\big\|\d_{ \overline{\na \cH(R^{k-1}_{a, 0}+\si h_{\tilde{v}^{k-1}})}}-\d_{ \overline{[\na\cH, D_t^2]\mathfrak{K}^{k}_a}}-D_{t*}\delta_{\overline{\na \cH(\cP^k_{t,1})}}\big\|_{L^2(\Om_*)}\leq CE^k_\d(t).
\]
Combining all these estimates above, the proof is finished.

\ef

Now, we are able to conclude about the convergence result.
\begin{prop}\label{prop:C-S}
The sequence $(\mathfrak{K}^{k}_a, P^k_{v, v}, d^{k}_l, d^{k}_r)$ is a Cauchy sequence.
\end{prop}
 \no{\bf Proof.}
We follow the steps in Theorem \ref{lower-order energy estimates} and Section 6.3 in \cite{MW3} to conclude that there exists a constant $T$ small enough and $A$ large enough such that when $a\geq A$, we have
\[
E^k_\d(t)\leq C\int_0^TE^k_\d(t)dt.
\]
As a result, this implies immediately that $(\mathfrak{K}^{k}_a, P^k_{v, v}, d^{k}_l, d^{k}_r)$ is convergent.

 \ef

We are finally in a position to finish the proof for Theorem \ref{main thm}.

\noindent{\bf Proof of Theorem \ref{main thm}.} We only present the sketch for the proof here. In fact, since we have proved in Proposition \ref{prop:C-S} that  $(\mathfrak{K}^{k}_a, P^k_{v, v}, d^{k}_l, d^{k}_r)\in \Sigma$ is a Cauchy sequence, we know immediately that there exists $(\mathfrak{K}_a, P_{v, v}, d_l, d_r)\in \Sigma$ satisfying 
\[
(\mathfrak{K}^{k}_a, P^k_{v, v}, d^{k}_l, d^{k}_r)\to(\mathfrak{K}_a, P_{v, v}, d_l, d_r)\qquad\hbox{in}\quad \Sigma.
\]
As a result, one can show in a standard way that $(\mathfrak{K}_a, P_{v, v}, d_l, d_r)$ satisfies  system \eqref{eq: IS-C}. Moreover, one also has $(\mathfrak{K}_a, P_{v, v}, d_l, d_r)\in \Sigma_h$ in the proof of Proposition \ref{prop:C-S}.
 \ef

In the end, we go back to our water-waves system $\mbox{(WW)}$. In fact, thanks to system \eqref{eq: IS-C} for $(\mathfrak{K}_a, P_{v, v}, d_l, d_r)\in \Sigma$, we derive the mean curvature $\kappa$, which defines the  free surface. Based on the knowledge of $\G_t$, we also obtain $v$ by Section 6.2.
Therefore, using similar arguments as in  Section 6.4 \cite{MW3} and thanks to discussions in Section 6.3, we can finally retrieve the solution $(v, P)$ to the water-waves system $\mbox{(WW)}$.

\bigskip

\noindent{\bf Acknowledgement}.  The authors would like to thank Chongchun Zeng for very fruitful discussions. The author Mei Ming is supported by NSFC no.12071415. The author Chao Wang is supported by NSFC no.12071008.

\end{document}